\documentclass[11pt,reqno]{amsart}
\usepackage{graphicx,amsmath,amsthm,amsfonts,amssymb,mathrsfs}
\usepackage{mhequ}
\usepackage{wrapfig}
\usepackage[usenames,dvipsnames]{xcolor}
\usepackage{hyperref}
\hypersetup{%
	colorlinks=true, linkcolor=blue,
	citecolor=ForestGreen
}
\usepackage[paper=letterpaper,margin=1.1in]{geometry}
\usepackage[inline]{enumitem} 
\usepackage{caption,subcaption}
\usepackage{tikz}
\usetikzlibrary{decorations.pathreplacing}
\usepackage{marginnote}
\usepackage{siunitx}
\usetikzlibrary{calc} \tikzset{>=latex}
\usetikzlibrary{calc,decorations.markings}
\usetikzlibrary{positioning,shapes,decorations.text}

\theoremstyle{plain}
\newtheorem{thm}{Theorem}[section]

\theoremstyle{definition}

\newtheorem{rmk}[thm]{Remark}
\numberwithin{equation}{section}
\numberwithin{table}{section}
\allowdisplaybreaks

\makeatletter
\newcommand{\myitem}[1]{%
\item[#1]\protected@edef\@currentlabel{#1}%
}
\makeatother

\usepackage{acronym}
\acrodef{SHE}{Stochastic Heat Equation}
\acrodef{SBE}{Stochastic Burgers Equation}
\acrodef{SPDE}{Stochastic Partial Differential Equation}
\acrodef{KPZ}{Kardar--Parisi--Zhang}
\acrodef{6V}{Six Vertex}
\acrodef{ASEP}{Asymmetric Simple Exclusion Process}

\newcommand{\eps}{\varepsilon}
\newcommand{\E}{\mathbf E}
\newcommand{\T}{\mathbf T}

\newcommand{\Z}{\mathbf{Z}}

\newcommand{\R}{\mathbf{R}}
\newcommand{\e}{\varepsilon}

\renewcommand{\bar}[1]{\overline{#1}}

\renewcommand{\tilde}[1]{\widetilde{#1}}

\usepackage{graphicx}
\graphicspath{{./figs/}}

\begin{document}
\title{Some recent progress in singular stochastic PDEs}
\author[I.\ Corwin]{Ivan Corwin}
\address{I.\ Corwin,
	Departments of Mathematics, Columbia University,
	\newline\hphantom{\hspace{15pt}I.\ Corwin}
	2990 Broadway, New York, NY 10027}
\email{corwin@math.columbia.edu}
\author[H.\ Shen]{Hao Shen}
\address{H.\ Shen,
	Departments of Mathematics, University of Wisconsin - Madison,
	\newline\hphantom{\hspace{15pt}H.\ Shen}
	480 Lincoln Drive, Madison, WI 53706}
\email{pkushenhao@gmail.com}
%


%
\begin{abstract}
\noindent
Stochastic PDEs are ubiquitous in mathematical modeling. Yet, many such equations are too singular to admit classical treatment. In this article we review some recent progress in defining, approximating and studying the properties of a few examples of such equations. We focus mainly on the dynamical $\Phi^4$ equation, KPZ equation and Parabolic Anderson Model, as well as touch on a few other equations which arise mainly in physics.
\end{abstract}
\maketitle
\tableofcontents

\section{Introduction}
\label{sec:intro}

Partial differential equations (PDEs) and randomness are ubiquitous constructions used to model both mathematical and physical phenomena.
For instance, PDEs have been used for centuries to describe the building block laws of physics, and model aggregate macroscopic phenomena such as heat conduction, diffusion,
electro-magnetic dynamics, interface and fluid dynamics. Randomness has become a default paradigm for modeling systems with uncertainty or with many complicated or chaotic microscopic interactions.

Combining these two approaches leads to the study of stochastic PDEs (SPDEs) in which the coefficients or forcing terms in PDEs are described via certain random processes.
While SPDEs have become increasingly important in applications, there remain many fundamental mathematical challenges in their study---in particular, showing how they arise from microscopic particle based models remains a major source of research problems and has seen some radical progress in the past decade.

The purpose of this article is to introduce a few important classes of SPDEs and describe how they arise and the mathematical challenges that go along with demonstrating that.
Though this article will mainly focus on nonlinear systems, we will start our investigation in Section \ref{sec:linear} in the simpler and more classical setting of linear SPDEs which are very well-understood. In Section \ref{sec:NSPDE} we turn our attention to nonlinear SPDEs, and introduce
our two main examples (the dynamical $\Phi^4$ equation and the KPZ equation) along with a host of other important SPDEs which arise in physics. Our discussion in this section is heuristic and ignores some of the serious mathematical challenges which arise when one tries to make sense of what it
means to ``solve'' an SPDE. This challenge is addressed in Section \ref{sec:Well-posedness}. In the course of making sense of SPDEs, there are often ``renormalizations'' which arise (effectively changing the equation). Section \ref{sec:Renormalize-physical} describes how these renormalizations have physical
meaning and arise in certain discrete approximation schemes for the continuum equations. Finally, Section \ref{sec:universal} seeks to demonstrate how these SPDEs (in particular, the KPZ equation) arise as universal limits from microscopic systems.

Before proceeding to our main text, one disclaimer. Our aim is to make this material approachable to non-experts. As such, we will not state precise theorems or give proofs, but rather will attempt to provide some intuition behind results and the challenges which accompany proving them.
An interested reader can find much more detail and precision in the works cited; or can consult other survey articles such as \cite{MR3828162,Gubinelli2018panorama}, \cite{ChandraWeber}, \cite{MR3336866} and \cite{HairerICM2014}.

\subsubsection*{Acknowledgements}
Ivan Corwin was partially supported by the Packard Fellowship for Science and Engineering, and by the NSF through DMS-1811143 and DMS-1664650. Hao Shen was partially supported by the NSF through DMS-1712684 and DMS-1909525. We are grateful to Weinan E, Massimiliano Gubinelli, Martin Hairer, Konstantin Matetski, Nicolas Perkowski and Li-Cheng Tsai for providing helpful comments while we complete this survey.

\section{A first (linear) SPDE}
\label{sec:linear}

We will start out discussion on linear SPDEs with the ``stochastic heat equation'' which is driven by a random additive noise term $\xi$:
\begin{equation}\label{e:LSHE-intro}
\partial_t u(t,x) = \partial_x^2 u(t,x)
+\xi(t,x)
\end{equation}
where $\xi$ is the so called space-time white noise. It will take a bit of work to define this noise and make sense of what it means to solve this equation. However, before going down that route, we will first address the question of what sort of physical system does this model? In particular,
we will explain heuristically how this equation arises from a simple microscopic model of polymers in liquid.

Consider modeling a polymer chain (e.g. composed of DNA or proteins) in a liquid. A simple model involves describing the polymer by a string of $N$ beads that are linked together sequentially by springs and subject to kicking by noise,
as shown in the following figure\footnote{The dots represent the locations of the beads connected by zigzag edges. The arrows represent the forces acting to move the beads---part due to spring force with the previous and subsequent beads, and part due to some random kicking force (yet to be specified).}
where $N=13$:
\begin{center}
\includegraphics[scale=.8]{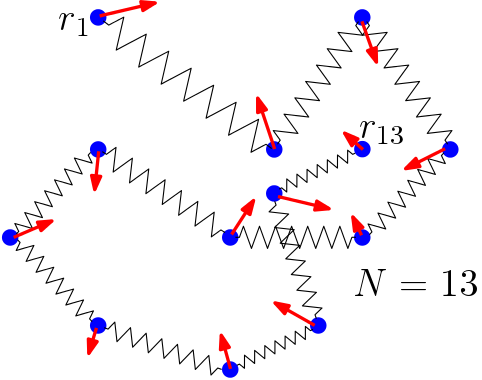}
\end{center}
Imagine that each bead of the polymer is ``kicked'' by the surrounding liquid molecules. In our simplified model\footnote{As usual one always has to make various simplifying assumptions in order to describe a complicated physical system via a mathematically analyzable model.
It is natural to ask whether having random kicking leads to a reasonable microscopic model. After all, the liquid itself is governed by certain physical laws of motion for its particles. Such concerns arose  early in the development of Brownian motion as the model for a single
tracer particle moving in a liquid---see \cite{Brush} for a nice historical review. We do not provide further justification for this as a reasonable microscopic model here.},
we describe such a system via the following equations of motion
for the position $r_i \in \R^3$ of the $i$-th bead:
\begin{equation} \label{e:r-i}
dr_i (t)=  (r_{i+1}  (t)- r_i (t))dt +   (r_{i-1} (t) - r_i (t)) dt+
\tikz[baseline]{
            \node[fill=red!20,anchor=base] (t1)
            {$ w_i (t)$};
            \path[->]  node[rectangle,draw,fill=blue!20,rounded corners=.8ex,text width=1.6cm,align=center] at (1.5,0.6) { \tiny  \baselineskip=7pt random kicking at time $t$\par} edge [bend right,thick] (t1);}
(t\in \R_+,i=1,\cdots,N)
\end{equation}
where a boundary condition is given by fixing $r_0-r_1\equiv 0$ and $r_{N+1}-r_N\equiv 0$ for all time.
Eq~\eqref{e:r-i} means the following
\begin{itemize}[leftmargin=.3in]
\item
The linear drift terms $(r_{i+1} - r_i)dt $ and   $ (r_{i-1} - r_i) dt$ arise from assuming a linear spring force between the $i$-th bead with its neighboring (in the sense of label number) beads.  Without the kicking term $w_i(t)$, Eq. \eqref{e:r-i} would simply be a coupled system of $N$ ordinary differential equations.
 \item The term $w_i(t)$ represents the {\it random kicking} that is experienced by the $i$-th bead at time $t$. We make the simplifying assumption that the kicking is ``overdamped''\footnote{Essentially, this means that the kicks occur instantaneously in time and do not result in any inertia. This effectively decouples the various kicks.} and model the kicks in terms of random jumps in the location of the $r_i(t)$. Namely, for each particle $r_i$ there is a random sequence of
 ``kicking'' times\footnote{It is natural to assume the gaps between times $t_{i,j}-t_{i,j-1}$ are chosen according to independent exponential random variables of mean 1. In this case, the times are distributed as a ``Poisson point process'' of intensity 1.} $t_{i,1}<t_{i,2}<\cdots$. At the kicking time $t_{i,j}$, we update $r_i \mapsto r_i + w_i(t_{i,j})$ where  $w_i(t_{i,j})$ is an $\R^3$-valued random variable.
We assume that the $w_i$ are statistically isotropic (i.e. their distribution is invariant under rotation) and all independent and identically distributed. Note that the resulting $r_i(t)$ process is piecewise continuous, with jumps occurring at the kicking times.
\end{itemize}

The question with which  we are concerned is what happens to the polymer when its length grows, and possibly space and time are scaled accordingly. By default one might expect that as $N$ increases, the complexity of studying this system goes likewise. However, it turns out that there is a very tractable continuum limit for the evolution of our polymer model. That is to say, in the scaling limit, things simplify! In fact, this limit is quite robust and (up to some scaling constants) is not affected by various changes in the microscopic model, such as how we model the kicking (e.g. different distribution on the $w_i$ or on the kicking times). This robustness can, itself, be seen as evidence that the microscopic model may be reasonable.

With the aim of demonstrating a continuum limit of our model, think of $N$ as large and define 
\begin{equ}[e:def-uN]
u_N(t,x) := \frac{1}{\sqrt{N}} r_{[xN]} (N^2 t)\in \R^3
\end{equ}
where $x\in [0,1]$ encodes the label $i$ via $i=[xN]$ (the closest integer to $xN$). The linear drift term in \eqref{e:r-i} is in fact a discrete Laplacian, so under our diffusive scaling (that is, scaling $x$ by $N$ and $t$ by $N^2$)
it approximates a continuum Laplacian $\partial_x^2$. Thus, for large $N$ one can expect that the following relation approximately holds
\begin{equ}[e:uN-equa]
\partial_t u_N(t,x) \approx \partial_x^2 u_N(t,x)
+ N^{\frac32} w_{[Nx]}(N^2 t)\;.
\end{equ}
In the scaled coordinates, the kicking starts to add up. Namely, in a time-space region of size $dt \times dx$, there are roughly $N^2 dt \times N dx$ kicks.
By the Central Limit Theorem the sum of $N^2\times N$ independent identically distributed random variables divided by $N^{\frac32}$ converges to a Gaussian random variable.
Thus, on each time-space region, the kicking adds up to a Gaussian random variable with variance equal to the area of the region.
Different regions have covariance given by the area of their overlap.
This limit is called space-time (or time-space, given our ordering of variables) white noise and is denoted by $\xi(t,x)$.
This heuristic lead us to the following type\footnote{This type of convergence result  was first  proved by Funaki \cite{MR692348} in the slightly different setting where the $r_i$ are driven by Brownian motions. In the present setting, we do not know if a precise result
of this sort has been proved (though have no doubt that it can be). We are suppressing  coefficients which may (depending on the nature of the discrete noise) arise in the limiting equation.} of limit as $N\to \infty$,
\begin{equation}\label{e:LSHE}
u_N(t,x)\to u(t,x) \qquad \textrm{where} \qquad
\partial_t u(t,x) = \partial_x^2 u(t,x)
+\tikz[baseline]{
            \node[fill=red!20,anchor=base] (t1)
            {$ \xi(t,x)$.};
            \path[->]  node[rectangle,draw,fill=blue!20,rounded corners=.8ex,text width=1.5cm,align=center] at (1.5,0.8) {\scriptsize  \baselineskip=7pt space-time white noise \par} edge [bend right,thick] (t1);
        }
\end{equation}

Eq. \eqref{e:LSHE} is our first example of an SPDE---it is called the linear stochastic heat equation with additive noise. Even in this simple linear example we encounter an equation which requires some work to make sense of because of the noise.
For convenience of our exposition, from this point, we will think of $x$ as a {\it spatial} variable (although in our example it actually stands for the parametrization of the limiting polymer length);
and although $u$ and $\xi$ are $\R^3$-valued,
the three components are completely independent (decoupled), so
it will be convenient to simply consider Eq.~\eqref{e:LSHE} as $\R$-valued instead of $\R^3$-valued in the sequel.

\vspace{2ex}

Let us look at Eq. \eqref{e:LSHE} more closely, with spatial dimension $d$  now being arbitrary (recall $d=1$ corresponds with the above polymer example)
\begin{equation}\label{e:LSHEd}
\partial_t u(t,x) = \Delta u(t,x)  + \xi(t,x)
\qquad x\in [0,1]^d \subset \R^d\;,
\end{equation}
with, for instance, periodic boundary condition.

Consider the case $d=0$, where \eqref{e:LSHEd} becomes the stochastic (ordinary) differential equation $\partial_t u(t) = \xi(t)$.  The $d=0$ white noise $\xi(t)$ is defined to be the derivative of Brownian motion so that $u(t) = \int_0^t \xi(s)$ is a Brownian motion.
Of course, Brownian motion is famously almost nowhere differentiable so $\xi$ is not defined as a function. Rather, $\xi(t)$ can be defined as a random {\it distribution} in a suitable {\it negative} regularity space (such as a negative Sobolev space).
Since the H\"{o}lder regularity of Brownian motion is $1/2-$ (meaning any exponent below $1/2$), its derivative is said to have regularity $-1/2-$. There are other ways to define $\xi(t)$.
For instance, if we restrict to a periodic $t\in [0,T]$, then $\xi(t) = \sum_{k\in \Z} N_k e_k(t)$ where the $N_k$ are independent identically distributed Gaussian random variables, and the $\{e_k\}_{k\in \Z}$ constitute an orthonormal basis of $L^2([0,T])$.
Alternatively, one can define $\xi(t)$ via the machinery of Gaussian processes (see, for instance \cite{janson_1997}) wherein it suffices to specify its mean and variance. 
As definition $\xi$ is mean zero,
and since $\xi$ is a {\it distribution} as mentioned above, its covariance
$$
\mathbf E\big(\xi(t)\xi(t')\big) =\delta(t-t')
$$
(here $\mathbf E$ represents the expectation value operator and $\delta$ is a Dirac delta function)
must be interpreted in a distributional sense as well. For a smooth test functions $f$, one defines the stochastic integral $\xi(f)=\int f(t)\xi(t)$. Then $\xi$ is defined by the property that $\mathbf E\big(\xi(f)\big)=0$ for all test functions $f$, and for test functions $f$ and $g$,
\begin{equ}[e:xi-cov-fg]
\mathbf E \big(\xi(f)\xi(g)\big)
=\langle f,g \rangle_{L^2}.
\end{equ}
The random distribution $\xi$ is defined from this information using Kolmogorov's continuity theorem.

For general dimension $d\geq 1$, space-time white noise can be defined via analogous methods. As a Gaussian process, $\xi(t,x)$ is a random distribution  with covariance
\begin{equ}[e:xi-cov-delta]
\mathbf E\big(\xi(t,x)\xi(t',x')\big) =\delta(t-t')\delta(x-x'),
\end{equ}
where the last $\delta$ is the Dirac delta function on $d$-dimensional space.
Its action on space-time test functions $f$ and $g$ has covariance
given by \eqref{e:xi-cov-fg} where the $\langle \cdot,\cdot \rangle_{L^2}$ is $L^2$ product over space-time.

As the dimension $d$ increases, the regularity of $\xi$ {\it decreases}. We will work with spaces  of space-time distributions (for $\alpha< 0$) or functions (for $\alpha\geq 0$) denoted by $\mathcal{C}^{\alpha}$. These are essentially equivalent to the Besov spaces $\mathcal B^{\alpha}_{\infty,\infty}$ in harmonic analysis, and their precise definitions can be given via wavelets in \cite[Eq.~(3.2)]{Hairer14}. The smaller $\alpha$ corresponds to less regular (or more {\it singular}) functions or distributions. A well-known result is that 
 that space-time white noise $\xi \in \mathcal C^{\alpha}$ for  $\alpha <-\frac{d+2}{2}$. \footnote{As we will primarily work with parabolic equations, these spaces $\mathcal{C}^{\alpha}$ have a built-in {\it parabolic scaling} between time and space wherein time regularity is doubled. For instance a $\mathcal C^{2}$ function has second continuous spatial derivatives and first continuous time derivative. Extending the situation discussed earlier for ($d=0$) white noise, we have that space-time white noise $\xi \in \mathcal C^{\alpha}$ for  $\alpha <-\frac{d+2}{2}$. The case $d=0$ have $\xi\in \mathcal{C}^{-1-}$ which, given the doubling of time regularity corresponds to the $-\frac12-$ regularity discussed above.}

Having made sense of $\xi$, it remains to understand what it means to ``solve'' Eq. \eqref{e:LSHEd}. For a linear equation as \eqref{e:LSHEd}, the meaning of ``solution''
is not hard to define; essentially one only needs to give a suitable meaning to the inverted linear differential operator acting on the noise $\xi$:
given an  initial data $u(0,x)=u_0(x)$,
the solution to \eqref{e:LSHEd}  is defined by
\begin{equ} [e:solve-SHE]
u= (\partial_t - \Delta)^{-1} \xi  + e^{t\Delta} u_0\;.
\end{equ}
Here $ e^{t\Delta}$ is the heat semi-group so that
$e^{t\Delta} u_0 (x):=  \int_{[0,1]^d} P(t,x-y) u_0(y)dy$
solves  the classical (deterministic) heat equation starting from $u_0$, with heat kernel $P(t,x):=(4\pi t)^{-\frac{d}{2}} e^{-\frac{|x|^2}{4t}}$.
The expression $ (\partial_t - \Delta)^{-1} \xi $ also acts as an integral operator on $\xi$ via
\begin{equ} [e:heat-kernel-conv]
(\partial_t - \Delta)^{-1} \xi (t,x):=\int_0^t \int_{[0,1]^d} P(t-s,x-y)\xi(ds,dy)
\end{equ}
which is the space-time convolution of the heat kernel $P(t,x)$ with $\xi$.
Just like $\xi$ itself, $ (\partial_t - \Delta)^{-1} \xi $ is also a {\it well-defined}
random distribution.

To get a bit more flavor of ``solution theories'' of stochastic PDE,
we list some well-known  properties  for \eqref{e:LSHEd}:
   \begin{enumerate}[leftmargin=.3in]
   \setlength\itemsep{1ex}
 \item[(P1)]
 The solution, in the above sense, is obviously {\it unique},
 since the difference  of two solutions
 would solve a deterministic heat equation with zero initial condition
which must be zero.
\item[(P2)]
With the aforementioned regularity of $\xi$,
by standard parabolic PDE theory, in particular the Schauder estimate which states that the operator $(\partial_t - \Delta)^{-1}$ increases regularity by $2$,
one  has
$u\in \mathcal C^{\alpha}$ for any $\alpha< -\frac{d+2}{2} +2 =\frac{2-d}{2} $.
In particular,
$u$ is (almost surely) a random continuous function in $d=1$, and a random distribution in $d\ge 2$.
So the limiting polymer parametrized by $x\in [0,1]$ in the above example is a random continuous curve.
\item[(P3)]\label{P3Gaussian}
The random distribution $ (\partial_t - \Delta)^{-1} \xi $ has {\it Gaussian probability law}. This is because $\xi$ is Gaussian and any {\it linear} combination of Gaussian random variables  is still Gaussian.
\item[(P4)]
Eq. \eqref{e:LSHEd} has an {\it invariant measure} called the Gaussian free field. This is a Gaussian random field on $[0,1]^d$ with covariance given by the Green's function of the Laplace $\Delta$.
Being invariant means that if the initial data $u_0$ is random and
distributed as Gaussian free field, then $u(t,\cdot)$
has the same law of  Gaussian free field  for all $t>0$.
On the other hand starting from arbitrary $u_0$,
the law of $u(t,\cdot)$ will approach that of the Gaussian free field as $t\to \infty$. We refer to \cite{MR2322706} for a nice review about the Gaussian free field.
\item[(P5)]
Eq. \eqref{e:LSHEd} is scaling invariant in any dimension $d$, namely, $\tilde u(t,x):= \lambda^{\frac{d-2}{2}} u(\lambda^2 t, \lambda x)$ satisfies $\partial_t \tilde u(t,x) = \Delta \tilde u(t,x)  +\tilde\xi(t,x)$
where $\tilde\xi(t,x):= \lambda^{\frac{d+2}{2}} \xi(\lambda^2 t,\lambda x) \stackrel{law}{=}\xi(t,x)$.
The last scaling relation of the white noise can be seen from its covariance \eqref{e:xi-cov-delta} recalling that the Dirac $\delta$ on $n$-dimensional space has scaling dimension $-n$.
Note that the scaling taken in \eqref{e:def-uN}, \eqref{e:uN-equa} was precisely the one
which leaves the limit equation invariant.
%
\end{enumerate}

So far a reader who is new to the area of SPDEs should have acquired the following message: the ``solution theories'' of SPDEs share some of the same fundamental challenges as in the study of classical PDEs. These include showing that solutions exist (or can be defined) both locally or globally, are unique
within certain regularity classes, and arise as a scaling limits for various approximation schemes. The rest of this article will focus on recent progress on these challenges for nonlinear SPDEs. Before doing so, let us briefly remark on another important challenge present both for PDEs and SPDEs---explicitly representing solutions via formulas.

Linear PDEs always admit explicit solutions. Linear SPDEs, as we saw above, have solutions which are random Gaussian processes with explicit mean and covariance (computable from the equation explicitly). Most nonlinear PDEs {\it do not} admit explicit solutions; those that do are generally related to
the area of ``integrable systems''. Likewise, most nonlinear SPDEs do not admit explicit descriptions for the probability distribution of their solutions. There are, however, a few special SPDEs (such as the KPZ equation discussed below) which can be explicitly solved in this sense. The study of such SPDEs
fall under the area of ``integrable probability'' or ``exactly solvable systems'', see for instance \cite{CorwinICM2014, BP16b} and references therein. We will not pursue this direction further in this article.

\section{Nonlinear SPDEs}
\label{sec:NSPDE}
Linear systems are often insufficient to effectively model many interesting phenomena. Indeed, as we will now see, nonlinear SPDEs arise in a number of important areas of physics (and many other directions which we will not discuss). Such nonlinear systems, however, are generally much more challenging to work with. Before coming to that, let us start with a few examples.

Consider a piece of magnet which is being heated up:

\begin{center}
\includegraphics[scale=0.3]{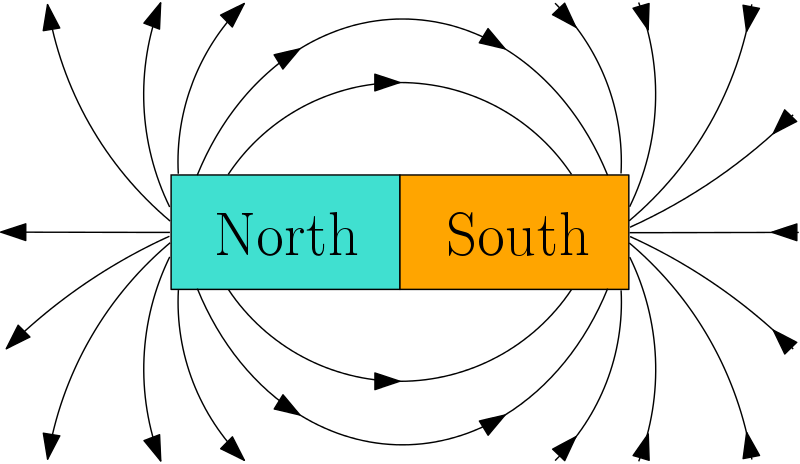}
$\qquad\qquad$
\includegraphics[scale=.35]{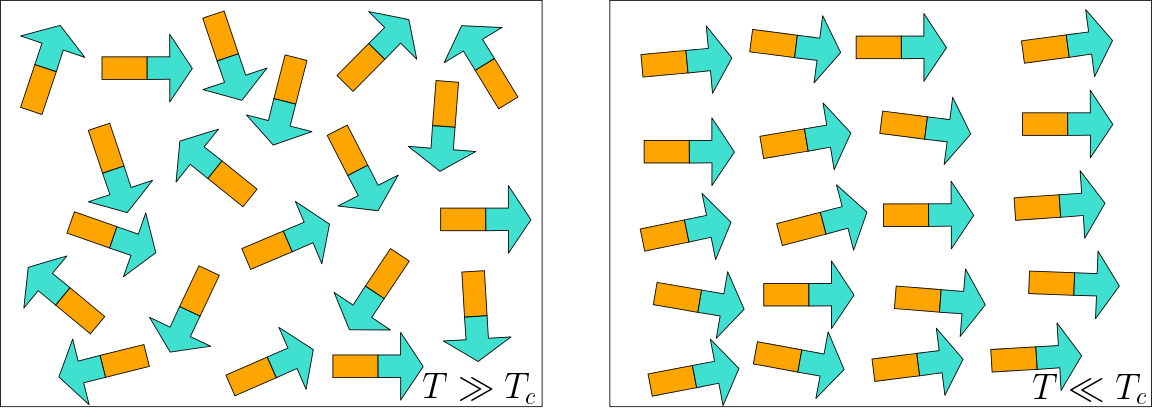}
\end{center}
%

As the temperature $T$ increases, the
magnetic field produced by the magnet weakens,
and
at a critical temperature $T_c$ known as the ``Curie temperature''\footnote{The Curie temperature is named after Pierre Curie who first experimentally demonstrated that certain magnets lost their magnetism entirely at a critical temperature.}
the magnetic field
disappears. Though various  magnet materials have different microscopic structures, a common physical explanation for magnetism is that it comes from the alignment of the magnetic moments of many of the atoms in the material.
As a simplified mathematical model one can  imagine that a magnet is made up of millions of tiny arrows (or spins),
with directions oscillating over time.
Below the Curie temperature, i.e. $T<T_c$, the spins tend to align
in order to minimize an interaction energy (which energetically prefers alignment), which
causes a macroscopic magnetization (shown in the right picture);
above the Curie temperature  $T>T_c$, the spin configurations
are much more disordered due to strong thermal fluctuation\footnote{In statistical mechanics, thermal fluctuations are random deviations of a system from its low energy state. All thermal fluctuations become larger and more frequent as the temperature increases, and likewise they decrease as temperature approaches absolute zero.
Thermal fluctuations are a basic manifestation of the temperature of systems.}
and as a result the magnetic fields cancel out (shown in the left picture).


%


A general mantra in statistical physics holds that ``interesting scaling limits arise at critical points''. In particular, here we would like to understand what happens to the spin system when the temperature $T$ approaches $T_c$ while time and space are tuned accordingly.
Near criticality, the spins start to oscillate more and more drastically and the small scale disorder starts to propagate to larger and larger scales. The resulting magnetic field fluctuations are believed to be described by the following nonlinear SPDE:
\[
  \partial_t \Phi = \Delta \Phi - \Phi^3  + \xi
  \tag{$\Phi^4$ Equation}
\]
when $\Phi=\Phi(t,x)$ and the spatial dimension of $x$ is $1,2$ or $3$. This is called the ``dynamical $\Phi^4$ equation'' since the deterministic part arises from the gradient of an energy $\int \frac12 |\nabla\Phi|^2+ \frac14 \Phi^4\,dx$.
We will return to this equation later in Section~\ref{sec:phy-Phi4} and describe how it arises from a particular model of magnets.

As another example,  we consider a model for interface growth,
where each point of the interface randomly grows up or drops down over time, with a trend to locally smooth the interface out (like the spring force in the polymer example in Section~\ref{sec:linear}). Such systems are ubiquitously found in nature---for instance, the left image in the figure
shows the end result of an interface grown in the ocean from a volcanic eruption.
\begin{center}
 \includegraphics[scale=0.045,trim={0 5cm  0 4.5cm},clip]{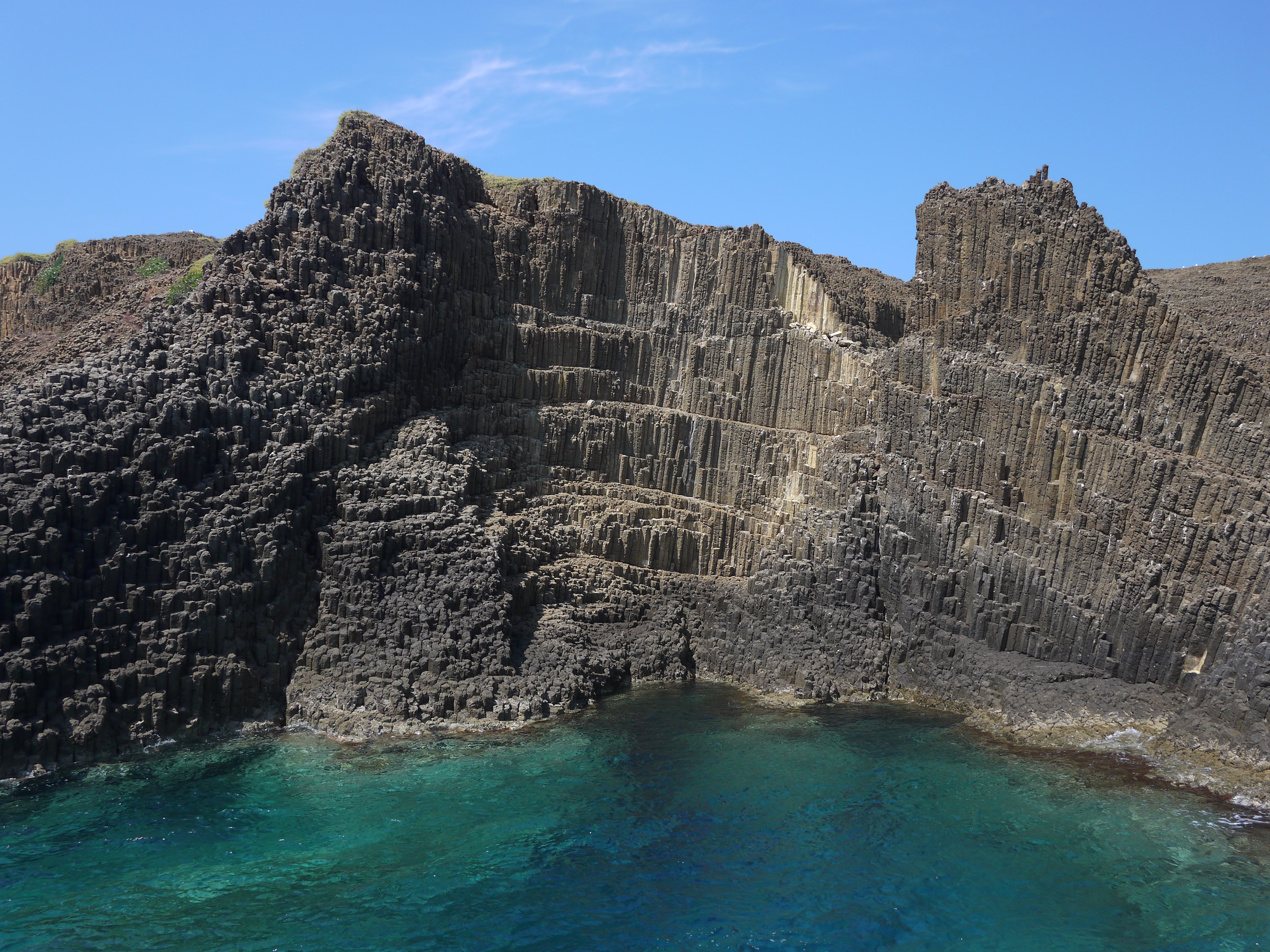}
 $\qquad$
 \includegraphics[scale=.55]{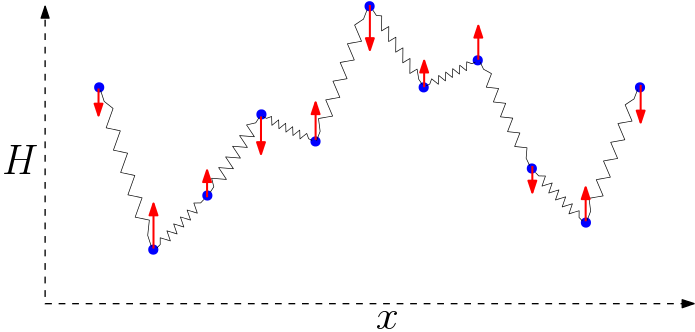}
\end{center}
We are interested in modeling the evolution of such interfaces.

To drastically simplify the situation, let us assume our interface is a one-dimensional function $H(t,x)$ for $x\in \R$. The simplest scenario is that the upward growth and downward drop of $H$ occurs {\it equally} likely, though randomly. In this case, it turns out that
the interface behaves similar to the one-dimensional version of the polymer in Section~\ref{sec:linear}
whose beads are kicked by {\it isotropic} random force (see the right image in the above figure).  Thus, the large scale fluctuation should be given by the linear SPDE \eqref{e:LSHE}.

In the asymmetric scenario where the interface is more likely to grow up than to drop down, one expects to see nontrivial fluctuation described by \eqref{e:LSHE} perturbed by a nonlinearity.
In particular, the asymmetry should not be too strong or it will overwhelm the local smoothing (the $\partial_x^2 H$ term) and randomness (the $\xi$ term) and not too weak or it will not change the limiting equation.
This critical tuning is called ``weakly asymmetric'' and results (under the same sort of scaling as in the symmetric case) in the following SPDE description for fluctuation of $H(t,x)$:
\[
  \partial_t H = \partial_x^2 H + (\partial_x H)^2  + \xi.
  \tag{KPZ Equation}
\]
Due to the asymmetry, the interface establishes an overall height shift. So, for the above limit, we must recenter into this moving frame. 
The KPZ equation was first proposed by Kardar, Parisi and Zhang  \cite{KPZ86},
see the nice review \cite{krug1991kinetic} for more background.
We will return to this equation later in Section~\ref{sec:universal} and describe why it arises from various models of interface growth.

\subsection{Some other important nonlinear SPDEs}

Besides the $\Phi^4$ and KPZ equations discussed above, there are a number of physically important equations---some of which we briefly review now.
The reader is warned that it is a formidable challenge to define the meaning of a solution to  nonlinear SPDEs driven by very singular noises. We postpone this important issue until Section \ref{sec:Challenge}.
\begin{itemize}[leftmargin=.3in]
   \setlength\itemsep{1ex}
\item {\bf Stochastic Navier-Stokes equation} (with  spatial dimensions $d=2,3$ of particular physical interest)
\begin{equ}[e:SNS]
\partial_t \vec u+\vec u\cdot \nabla \vec u =  \Delta \vec u-\nabla p+\vec \zeta ,
\quad \mbox{div } \vec u = 0
\end{equ}
where $p$ is the pressure, $\vec \zeta$ is a $d$-vector valued noise. 
When the noise $\vec \zeta$ is taken to be singular, for instance each component of $\vec \zeta$ is an independent space-time white noise, it models motion of fluid with randomness arising from microscopic scales, and in this case we refer to \cite{MR1941997} and \cite{Zhu2014NS} for well-posedness results.

We remark that while this article focuses on singular noises, when modeling {\it large scale} random stirring of the fluid,
the noise
$\vec \zeta$ is often  assumed to be smooth (called ``colored noise'' in contrast with white noise), 
and in fact the most important case is that the equation is driven by only a few number of random Fourier modes.
In these situations the {\it long-time behavior} is of primary interest, 
and various dynamical system questions such as ergodicity and mixing are studied.
There is a vast literature on this topic, and we only refer to 
the book \cite{MR3443633} and the survey articles \cite{MR2050597,flandoli2008introduction,kupiainen2010ergodicity}.
\item
{\bf Stochastic heat equation with multiplicative noise} in one spatial dimension:
 \begin{equ}[e:NSHE]
\partial_t u = \Delta u  + f(u)\xi
  \end{equ}
 where $f$ is some  continuous function. The It\^o solution theory successful for stochastic ordinary differential equations can extend to this stochastic PDE,
 see for instance the lecture notes \cite{MR876085}.
The specialization $f(u)=u$, i.e.
 \begin{equ}[e:SHE]
\partial_t u = \Delta u  + u\xi
  \end{equ}
  has a significant connection to the KPZ equation: one can formally check that if $H$ solves KPZ, then the Hopf-Cole transform $u:=e^H$ solves \eqref{e:SHE}.
 Other choices of $f$ are
 \[
 f(u)=\alpha\sqrt{u} \qquad
 \mbox{and}
 \qquad f(u)=\alpha\sqrt{u(1-u)}
 \qquad (\alpha\in \R)
 \]
 which (along with $f(u)=u$) arise in modeling population dynamics and genetics, see for instance  \cite{MR0321178,MR0378124}.
\item
{\bf Nonlinear parabolic Anderson model}, in  spatial dimensions $d=2,3$:
 \begin{equ}[e:gPAM]
\partial_t u = \Delta u  + f(u)\zeta
  \end{equ}
  where $f$ is a continuous function and $\zeta$ is a noise which typically is assumed to be spatially independent (i.e. white), but constant in time. This models the motion of mass through a random media. The assumption of constant in time noise is consistent with the regime where the mass is assumed to move much faster
   than the time scale in which the media changes. We refer to \cite{MR3406823,MR3779690} and references therein for well-posedness results.
   
 The parabolic Anderson model especially the linear case  ($f(u)=u$)  is a simple model which exhibits \emph{intermittency} over long time; for the study of long time behavior, one often considers the spatial-discrete equation with $\zeta$ being independent noises on lattice sites, see
 for instance the reviews \cite{MR1185878} and~\cite{MR3526112} for further discussion and references regarding long time behaviors of the parabolic Anderson model.
 \item
   {\bf Dynamical sine-Gordon  equation}
\begin{equ}[e:sine-G]
  \partial_t u = \Delta u  + \sin(\beta u) +\xi   \qquad
  (t,x)\in \R_+\times  \T^2.
     \end{equ}
This equation describes the natural dynamics of a  class of  two dimensional  systems that exhibits the Berezinskii-Kosterlitz-Thouless (BKT) phase transition \cite{berezinsky1970destruction,kosterlitz1973ordering,jos201340},
such as two-dimensional Coulomb gas and
certain condensed matter materials\footnote{These include thin disordered super-conducting granular films. The phase transition is from bound vortex-antivortex pairs at low temperatures to unpaired vortices and anti-vortices at some critical temperature.}. 
 See \cite{MR0441179,MR634447} for earlier studies of the model in equilibrium.
Here $\beta^2$ represents the inverse temperature,
and  $\beta^2=8\pi$ is the BKT critical point. See \cite{MR3452276,sineGordonwholeregime} for the construction of local solutions of this dynamic.
\item
{\bf Random motion of a curve} in an $N$ dimensional manifold  $M$ driven by $m$ independent space-time white noises (see  \cite{hairer2016motion,hairer2019geo,MR3787728}):
\begin{equ}  \label{e:Hairer-String}
\partial_t h^\alpha =  \partial_x^2 h^\alpha
	+  \sum_{\beta,\gamma=1}^N
		\Gamma_{\beta\gamma}^\alpha (h)  \partial_x h^\beta \partial_x h^\gamma
	+ \sum_{i=1}^m \sigma^\alpha_i (h) \xi_i
\qquad \alpha\in\{1,2,\cdots, N\} \;,
\end{equ}
where $h$ is a map from an interval to $M$,
$\Gamma_{\beta\gamma}^\alpha $ are  the Christoffel symbols for the Levi-Civita connection, and  $\{\sigma_i \}_{i=1}^m$ is a collection of vector fields on the manifold. This is a non-Euclidean generalization of \eqref{e:LSHE}.
\item
{\bf Stochastic Yang-Mills flow} in spatial dimensions $d\le 4$
\begin{equ}[e:YM]
 \partial_t A = - d^*_{A} F_A + \xi
 \end{equ}
 where the deterministic part  (without $\xi$) is the Yang-Mills  gradient flow
 introduced in \cite{MR1079726} which is extensively studied in geometry (see the monograph \cite{feehan2014global}).
Here,  in the setting of differential geometry, one fixes a Lie group, $A$ is a connection (or Lie algebra valued 1-form), $F_A$ is the curvature of $A$, $d_{A}$ is the covariant derivative operator and $d^*_{A}$ is its adjoint.
The noise  $\xi=\xi_1dx_1 + \cdots + \xi_d dx_d$ is a 1-form
with each component $\xi_i$ being an (independent copy of)
Lie-algebra valued space-time white noise.
See \cite{Shen2018Abelian} for some initial progress in $d=2$ in the case that the Lie group is Abelian. 
Note that \eqref{e:YM} is not a parabolic equation, and one usually adds 
an additional term $- d_A d^*_{A} A$ on the right hand side to obtain a parabolic equation
\begin{equ}[e:YM-heat]
 \partial_t A = - d^*_{A} F_A - d_A d^*_{A} A + \xi
 \end{equ}
 which is gauge equivalent with the original equation (``Donaldson-De Turck trick'').

The study of geometric equations with randomness such as \eqref{e:Hairer-String} and \eqref{e:YM} is of general interest.  Eq. \eqref{e:YM} is motivated by  the  problem of quantization of the Yang-Mills field theory; see also the next item.

\item
{\bf Stochastic quantization}.
This refers to a large class of singular SPDEs
arising from Euclidean quantum field theories defined via Hamiltonians (or actions, energy etc.). They were introduced by Parisi and Wu in \cite{ParisiWu}. Given a Hamiltonian $\mathcal H(\Phi)$ which is a functional of $\Phi$,
one considers a gradient flow of $\mathcal H(\Phi)$  perturbed by space-time white noise $\xi$:
\begin{equ}[e:SQE]
\partial_t \Phi = - \frac{\delta \mathcal H(\Phi)}{\delta \Phi} + \xi \;.
\end{equ}
Here $\frac{\delta \mathcal H(\Phi)}{\delta \Phi}$ is  the variational derivative of the functional $\mathcal H(\Phi)$;
for instance, when $\mathcal H(\Phi)=\frac12\int (\nabla \Phi)^2 dx$
is the Dirichlet form, $\frac{\delta \mathcal H(\Phi)}{\delta \Phi} =-\Delta \Phi$
and \eqref{e:SQE} boils down to the stochastic heat equation \eqref{e:LSHEd}.
Note that 
 $\Phi$ can be also multi-component fields, with $\xi$ being likewise multi-component. The aforementioned $\Phi^4$ equation, sine-Gordon
equation and stochastic Yang-Mills flow all belong to this class of stochastic quantization equations, each corresponding to a Hamiltonian $\mathcal H$.

The significance of these ``stochastic quantization equations'' \eqref{e:SQE} is that
given a Hamiltonian $\mathcal H(\Phi)$,
the formal measure
\begin{equ}[e:QFT]
\tfrac{1}{Z}e^{-\mathcal H(\Phi)} D\Phi
\end{equ}
 is formally an invariant measure\footnote{Being invariant means that if the initial condition of \eqref{e:SQE} is random with ``probability law'' given by \eqref{e:QFT}, then the solution at any $t>0$ will likewise be distributed according to this same ``probability law''.
 For readers familiar with stochastic ordinary differential equations, one simple example is given by the Ornstein-Uhlenbeck process $dX_t=-\frac12 X_tdt+dB_t$ where $B_t$ is the Brownian motion, and its invariant measure is the (one-dimensional) Gaussian measure $\frac{1}{\sqrt{2\pi}}e^{-\frac{X^2}{2}}dX$.}
 for Eq. \eqref{e:SQE}.
 Here  $ D\Phi$ is the formal Lebesgue measure and
 $Z$ is a ``normalization constant''.
We emphasize that \eqref{e:QFT} are only formal measures because, among several other reasons, there is no ``Lebesgue measure'' $D\Phi$ on an  infinite dimensional space and it is a priori not clear at all if the measure can be normalized.
These measures arise from {\it Euclidean quantum field theories}. In their path integral formulations quantities of physical interest are defined by expectations with respect to these measures.
The task of {\it constructive quantum field theory} is to give precise meaning or constructions to these formal measures, see the book \cite{MR1773042}. 

Given the very recent progress of SPDE,
a new  approach to construct 
the measure of the form \eqref{e:QFT} is to construct the {\it long-time} solution to the stochastic PDE \eqref{e:SQE} and average the distribution of the solution over time.
This approach has shown to be successful for the $\Phi^4$ model in $d\le 3$ 
in a series of very recent works, which starts with
 \cite{MR3719541} on the torus $\T^3$ where a priori estimates were obtained to rule out the possibility of finite time blow-up.
Then in \cite{gubinelli2018global,gubinelli2018pde} established a priori estimates for solutions on the full space $\R^3$ yielding the construction of $\Phi^4$ quantum field theory 
on the whole $\R^3$, as well as verification of some key properties that this invariant measure must satisfy as desired by physicists such as reflection positivity. See also \cite{albeverio2017invariant}. Similar uniform a priori estimates are obtained by \cite{moinat2018space} using  maximum principle.

 \item
{\bf Random (nonlinear) Schr\"odinger equation}:
\begin{equ}[e:RSchrE]
i \frac{du}{dt} = \Delta u +\lambda |u|^2 u + u \xi
\end{equ}
where $u$ is complex valued,
$\xi$ is a real valued {\it spatial} white noise.
The linear case $\lambda=0$ is a model for Anderson localization (a complex version of \eqref{e:gPAM}, see the recent works \cite{AllezChouk,MR3762097}).
In the nonlinear case, it describes the evolution of nonlinear dispersive waves in a totally disordered medium,
 with  $\lambda>0$ corresponding to the ``focusing case'' and $\lambda < 0$ to the ``defocusing case'', see \cite{conti2012solitonization,ghofraniha2012shock} for its physical background, and \cite{MR3785398,DebusscheMartin,gubinelli2018semilinear} for recent mathematical results.
 \item
{\bf (Nonlinear) stochastic wave equation}:
 \begin{equ}[e:NSWE]
\partial_t^2 u - \Delta u + F(u)  
 = \xi
\end{equ}
with given initial data $(u,\partial_t u)|_{t=0}$.
The linear case ($F=0$) in $d=1$ spatial dimension,
as Walsh  explained in  \cite{MR876085},  describes ``a guitar left outdoors during a sandstorm. The grains of sand hit the strings continually but irregularly.''
If $\xi(dt,dx)$ is the random measure of the number of grains hitting in $(dt,dx)$ (centered by subtracting the mean),
then $\xi$ should be space-time white noise since the numbers hitting over different time intervals or string portions
 will be essentially independent. The position $u(t,x)$ of the string should satisfy \eqref{e:NSWE} with $F=0$.
 
 Eq.~\eqref{e:NSWE} with non-zero $F$ are investigated in earlier works by
\cite{MR1396758,MR1640497} in spatial dimensions $d= 2,3$
and they proved that with just a  function $F$ (satisfying some ``nice'' properties)
the solution to \eqref{e:NSWE} is trivial, namely the same with the solution for $F=0$;
the reason of this triviality will be clear in the next subsection.

More recently \cite{MR3841850}
obtained nontrivial solutions 
with $F$  given  (formally! see the next subsection) by  $F(u)=\lambda u^k$ in $d=2$,
and $F(u)=u +u^2$ in $d=3$ in \cite{gubinelli2018paracontrolled}.
\cite{gubinelli2018semilinear} then studied a stochastic wave equation 
with $F(u)=u^3$ and multiplicative noise ($u\xi$ on the right hand side) in $d=2,3$.
\end{itemize}

\begin{rmk}\label{rmk:scaling}
Note that these nonlinear SPDEs are generally not scaling invariant,
 unlike the linear stochastic heat equation \eqref{e:LSHEd} which is scaling invariant in any dimension $d$ (recall property (P5) in the end of Section~\ref{sec:linear}).
 For instance, for the KPZ equation,
 $\tilde H(t,x):= \lambda^{\frac{d-2}{2}} H(\lambda^2 t, \lambda x)$ 
 will satisfy $  \partial_t \tilde H = \partial_x^2 \tilde H + \lambda^{\frac{2-d}{2}}(\partial_x \tilde H)^2  + \tilde \xi$
 where $ \tilde \xi$ is a new  white noise, thus not invariant unless $d=2$.
 (Indeed, any choice of three scaling components for $t,x,H$ cannot make four terms invariant.)
 This will turn out to be important for defining solutions to these  equations,
 see Remark~\ref{rmk:subcritical}.
 Also, as we will see in Section~\ref{sec:Renormalize-physical} and Section~\ref{sec:universal},
it would not be possible to derive these equations (that are not scaling invariant)
as limits of scaling certain physical models,
unless the physical models have a weak asymmetry, a long interaction range, or a weak intensity of noise, which sets an ``additional scale''.
\end{rmk}

%
%

\subsection{Challenge of solution theory}
\label{sec:Challenge}

Solution theories for SPDEs have been developed since 1970s.
The earlier progress was recorded by the books written in the '80s such as
\cite{MR683274,MR876085},
and more recent books such as \cite{MR1661761,MR2329435,MR1500166,DPZ,MR3410409}.
However, many very important equations including some of those listed above remained poorly understood---that is, until very recently. 

The difficulty in building  solution theories to nonlinear SPDEs
is that often these equations are {\it too singular}, namely the solution (if it exists) would have low enough regularity so that certain parts of the equation do not a priori make sense.
Indeed, recall that for the linear equation \eqref{e:LSHEd}
in $d$ spatial dimensions,
the solution is almost surely an element of
$\mathcal C^{\alpha}$ for $\alpha <\frac{2-d}{2}$,
which is continuous when $d=1$ and is a distribution when $d\ge 2$.
Since with nonlinear terms the solutions are not expected
 to be more regular, the ``$\Phi^3$" term in the $\Phi^4 $ equation
when $d\ge 2$
 is a priori meaningless, because distributions in general can not be multiplied.
 Similarly, for the KPZ equation,
if $d\ge 1$, $\partial_x H$ is distribution valued and thus ``$(\partial_x H)^2$" does not have a clear meaning.
For this type of singular SPDE, it is challenging to even {\it interpret what one means} by a solution.
\footnote{
The same issue is present for the dispersive equations; for instance the solution to the the stochastic wave equation \eqref{e:NSWE} in spatial dimensions $d\ge 2$ is distributional, 
and this is exactly the reason that the triviality result  by \cite{MR1396758,MR1640497} for the nonlinear problem should be expected. Indeed, their proofs are based on a ``Colombeau distribution'' machinery.}

Starting in the 80s, the idea of ``renormalization"  entered the study of SPDE \cite{MR815192,AlbRock91,MR2016604}. Recently, this idea has received far-reaching generalization
in work of Hairer \cite{Hairer14}, Gubinelli-Imkeller-Perkowski \cite{MR3406823} and many subsequent works.
The idea is to subtract terms with infinite constants
from the nonlinearities. Taking the $\Phi^4$ equation as an example with spatial dimension $d\le 3$,
one needs to consider the following ``renormalized'' $\Phi^4$ equation
\begin{equ}[e:Phi4-inf]
 \partial_t \Phi = \Delta \Phi - (\Phi^3 -\infty \Phi)  + \xi \;.
\end{equ}

This precisely means the following.
Since the origin of the problem is the singularity of the driving noise $\xi$,
one starts by  regularizing $\xi$. For instance, we take a space-time convolution of $\xi$ with a mollifier $\varphi_\eps$ that is a  smooth  function of space and time with support of size $\eps$
so that  $\varphi_\eps \to \delta$  (the Dirac delta distribution) as $\eps \to 0$. Now we consider the $\Phi^4$ equation driven by $\xi_\e$:
 \begin{equation} \label{e:Phi4eps}
 \partial_t \Phi_\eps = \Delta \Phi_\eps - \Phi_\eps^3  + \xi_\eps
 \qquad
 \mbox{where}
 \qquad
  \xi_\eps = \xi * \varphi_\eps
\end{equation}
For any $\eps>0$, due to the smoothness of the noise we can solve the above PDE in the classical sense, and
 $ \xi_\eps  \in \mathcal C^\infty$ implies $\Phi_\eps\in \mathcal C^\infty$.
 As $\eps\to 0$, $\xi_\eps \to \xi$, but $\Phi_\eps$ do not converge to any nontrivial limit!

 The idea is that
before the limit, one should insert renormalization terms (also often called ``counter-terms" in the context of
quantum field theory\footnote{In fact the corresponding quantum field theory requires a renormalization $\int \frac12 |\nabla\Phi|^2+\frac14 \Phi^4 - \frac{C_\e}{2} \Phi^2 \,d x$ for dimensions $2$ and $3$ which is well-known in physics.})
 \begin{equation} \label{e:RPhi4eps}
 \partial_t \Phi_\eps = \Delta \Phi_\eps
 - ( \Phi_\eps^3  - C_\eps \Phi_\eps) + \xi_\eps
\end{equation}
where   $C_\eps$ {\it diverges} as $\eps\to 0$ at a suitable rate.
If the sequence of constants $C_\eps$ is suitably chosen,
the sequence of smooth solutions $\Phi_\eps$ of \eqref{e:RPhi4eps} will converge to a nontrivial limit as $\eps \to 0$:
\[
\Phi = \lim_{\eps\to 0} \Phi_\eps
\]
This is what we mean by a solution to the (renormalized) $\Phi^4$ equation. Note that we do not attempt to make sense of a limiting equation \eqref{e:Phi4-inf}, but we construct $\Phi$ via by a limit procedure,
the limit  of solutions to a sequence of regularized and renormalized equations \eqref{e:RPhi4eps}.

The same renormalization procedure applies to KPZ equation (in one spatial dimension) and many  of the other singular SPDEs listed above.

This discussion prompts  several questions:

\begin{itemize}
\item Why does $\Phi_\eps$ converge and how does one choose suitable constants $C_\eps$ to make this convergence happen? Is the resulting limit unique or does it depend on the mollification?
This is essentially the question of ``well-posedness'' which will be addressed in Section~\ref{sec:Well-posedness}.
\item Why are we allowed to ``change" the equation by inserting new terms that are not negligible - in fact infinite? We address this ``renormalization'' question in Section~\ref{sec:Renormalize-physical}, in which we will see that
the SPDEs such as $\Phi^4$  arise as ``scaling limits'' of physical systems and the ``renormalization'' will turn out to have physical meanings in these systems.
\item How robust are these singular SPDEs  under different approximation schemes? This is a ``universality'' question, meaning that one singular SPDE should be able to  serve as the continuum large scale description of {\it a class of} systems which may have different small scale details.
We discuss this in Section~\ref{sec:universal}.
\end{itemize}

These questions are of course entangled in many ways. In terms of ``approximations and convergence'', the  procedure described as in \eqref{e:Phi4eps} is the simplest way of approximation and approaching a limit, but ``scaling limits''  of physical models are essentially also ways of obtaining  the limits.
In terms of ``uniqueness'', one expects to get the same SPDE limit
 not only for different choices of mollifications in \eqref{e:Phi4eps} but also via scaling limits of perhaps apparently very different models, which is what ``universality'' means.
Section~\ref{sec:Renormalize-physical} below will focus on the meaning of renormalization in physical models, and Section~\ref{sec:universal} will
provide more detailed discussions on deriving an SPDE from these physical models, and of course in all these endeavors one needs to first understand the meaning of a solution as discussed in Section~\ref{sec:Well-posedness}.

\section{Well-posedness of singular SPDEs}
\label{sec:Well-posedness}

We discuss how to choose suitable renormalization constants so that one can obtain a nontrivial limit for solutions to renormalized equations. Our exposition consists of two parts.
\begin{itemize}
\item
Starting from the '90s, solutions to renormalized singular SPDEs have been constructed, see for instance \cite{AlbRock91,MR2016604}.
Here we present an elegant  argument due to \cite{MR2016604} which illustrates a simple example of renormalization, plus a standard Picard iteration  (fixed point) PDE argument.
This argument, despite of its simplicity, yields solutions to several singular (but ``not too singular'') equations, such as the $\Phi^4$ equation in two spatial dimensions.
\item
The above argument fails for more singular SPDEs, such as $\Phi^{4}$ in three spatial dimensions and the KPZ equation in one spatial dimension. This motivates us to turn to a more robust approach, the theory of regularity structures introduced in \cite{Hairer14}.
We will also mention some alternative theories or methods such as paracontrolled distributions \cite{MR3406823} or renormalization group \cite{Antti}.
\end{itemize}

To focus our discussion in this section, we will work with the $\Phi^4$ equation.

\subsection{A PDE argument and renormalization}
\label{sec:PDE-renorm}

Consider the $\Phi^4$ equation, where the spatial variable takes values in the two-dimensional torus.
As explained above, we take a sequence of mollified noises $\xi_\eps$ and consider the mollified equation \eqref{e:Phi4eps}.

Write $u_\eps := (\partial_t - \Delta)^{-1} \xi_\eps$ for the ``stationary solution''\footnote{This corresponds to dropping the term involving initial data in \eqref{e:solve-SHE} and integrating time in
\eqref{e:heat-kernel-conv} from $-\infty$ instead of $0$. Stationarity means that the distribution of $u_{\eps}$ does not depend on $t$. This assumption will be convenient when performing moment calculations such as \eqref{e:Eu2}. Namely, the moments won't depend on space-time points.}
to  the mollified linear stochastic heat equation \eqref{e:LSHEd} $\partial_{t} u_{\eps} = \Delta u_{\eps} + \xi_{\eps}$. The key observation is that the most singular part of $\Phi$ is $u$, so if we write
\begin{equ}[e:DD1]
\Phi_\eps = u_\eps + v_\eps,
\end{equ}
we can expect the remainder $v_\eps$ to converge in a space of better regularity.
Subtracting this linear equation from \eqref{e:Phi4eps} gives
\begin{equ} [e:v-equ]
\partial_{t} v_{\eps}
=
\Delta v_{\eps}-(v_{\eps} + u_{\eps})^{3}
=
\Delta v_{\eps}-v_{\eps}^{3}
-3u_{\eps}v_{\eps}^{2}
-3u_{\eps}^{2}
v_{\eps}
-u_{\eps}^{3} \;.
\end{equ}
This equation looks more promising since the rough driving noise $\xi_{\eps}$ has dropped out.
This manipulation has not solved the problem of multiplying distributions, since the limit of $u_\eps$ is still a distribution valued in two spatial dimensions (as we discussed earlier, see the fact (P2) in the end of Section~\ref{sec:linear}).
However $u_{\eps}$ is a rather concrete object since it is {\it Gaussian} distributed ((P3) in the end of Section~\ref{sec:linear}). This makes it possible to study the behavior of $u_{\eps}^{2}$ and $u_{\eps}^{3} $ via probabilistic methods.

As an illustration, consider the expectation
\begin{equ}[e:Eu2]
\mathbf E [u_{\eps}(t,x)^{2}]
= \int P(t-s,x-y) P(t-r,x-z) \varphi_\eps^{(2)}(s-r,y-z) \,dsdydrdz
\end{equ}
where $\varphi_\eps^{(2)}$  is the convolution of $\varphi_\eps$ introduced in \eqref{e:Phi4eps} with itself and $P$ is the heat kernel introduced in \eqref{e:heat-kernel-conv}.
Due to singularity of the heat kernel $P$ at origin, this integral diverges like $O(\log \eps)$ as $\eps\to 0$ in two spatial  dimensions. Denoting $C_\eps:= \mathbf E [u_{\eps}(t,x)^{2}]$ (which does not depend on $(t,x)$ by stationarity of $u_\eps$),
this calculation indicates that in \eqref{e:v-equ} we should subtract $C_\eps$ from $u_{\eps}^{2}$, and subtract\footnote{The factor $3$ arises from three ways of choosing two powers of $u_\e$ from the cubic term $u_\e^3$.
A ``Wick theorem'' allows one to compute expectation of a product of arbitrarily many Gaussian variables, and in two dimensions the Wick renormalized power $u^{:n:}:= \lim_{\e\to 0} C_\e^{\frac{n}{2}} H_n(u_\e/\sqrt{C_\e})$.} $3C_\eps u_\eps$ from $u_{\eps}^{3}$.

This amounts to considering the renormalized $\Phi^4$ equation
 \begin{equ} [e:renormalizedPhi4]
 \partial_t \Phi_\eps = \Delta \Phi_\eps
 - ( \Phi_\eps^3 - 3C_\eps \Phi_\eps) + \xi_\eps \;.
\end{equ}
These ``renormalized powers'' of $u_\eps$ do converge to nontrivial limits.
In fact, thanks to Gaussianity of $u_\eps$, given a smooth test function
$f$ one can explicitly compute any probabilistic moment of
$ (u_{\eps}^{2}-C_\eps) (f)$ and prove its convergence.
By choosing $f$ from a  suitable set of wavelets or Fourier basis,
one can apply Kolmogorov theorem to  prove that $u_{\eps}^{2}-C_\eps $ and $u_{\eps}^{3}-C_\eps u_\eps $ converge in $\mathcal C^\alpha$ for any $\alpha<0$.
We denote these limits $u^{:2:}$ and $u^{:3:}$; they are elements of $\mathcal C^\alpha$ for any $\alpha<0$.

To summarize, we have found that the renormalization constants can be found through {\it explicit} moment calculations/expectations.

Passing \eqref{e:v-equ} to the limit, we get
\begin{equ} [e:v-equ-R]
\partial_{t} v
=\Delta v
- \Big(v^{3}+3u v^{2}
+3u^{:2:}v + u^{:3:} \Big)\;.
\end{equ}
We can prove local well-posedness of this equation as a {\it classical} PDE, by a standard fixed point argument. For this, we use a classical result in harmonic analysis  called ``Young's theorem'' which states that if $f \in \mathcal C^{\alpha}$,
$g \in \mathcal C^{\beta}$, and $\alpha+\beta>0$, then $f\cdot g \in \mathcal C^{\min(\alpha, \beta)}$. Thus if we assume that $v \in \mathcal C^{\beta}$ for, say, $\beta=1$, then the worst term in the parenthesis in  \eqref{e:v-equ-R} has regularity $\alpha$.
By the classical Schauder estimate which states that the heat kernel improves regularity by $2$, the following fixed point map
\begin{equ} [e:fix-pt-Phi42]
v \mapsto
(\partial_t - \Delta)^{-1}
\Big(v^{3}+3u v^{2}
+3u^{:2:}v + u^{:3:} \Big)
\end{equ}
is well-defined. Namely, it maps a generic element $v \in \mathcal C^{\beta}$  to a new element which is {\it again} in $\mathcal C^{\beta}$ for $\beta=1$. This is since for $\alpha<0$ sufficiently close enough to $0$ one has $\mathcal C^{\alpha+2}\subset \mathcal C^{\beta}$.
With a bit of extra effort, one can show that over short time interval the fixed point map is contractive and thus has a fixed point in $\mathcal C^{\beta}$, and this fixed point $v$ is the solution.
(The sharp result is $\mathcal C^{\beta}$ for any $\beta<2$.)
To conclude, one has $\Phi = u+ v$, which is the local solution to the renormalized $\Phi^4$ equation in two spatial dimensions.

The above argument was first used by Da Prato and Debussche in \cite{MR2016604} and it applies to other equations, for instance the stochastic Navier-Stokes equation \eqref{e:SNS} with space-time white noise on two dimensional torus \cite{MR1941997}.
Let us mention another, somewhat surprising application, that is the dynamical sine-Gordon equation \eqref{e:sine-G} in two spatial dimensions in the regime $\beta^2\in (0,4\pi)$.
The renormalized equation reads
\[
\partial_{t}u_\eps =\tfrac{1}{2}\Delta u_\eps + C_\eps \sin(\beta u_\eps )+\xi_\eps \;,
\qquad \beta\in \mathbf R\;,
\]
where $C_\eps$ is a renormalization constant which diverges like $O(\eps^{-\frac{\beta^2}{4\pi}})$.
By writing $u_\eps = \phi_\eps + v_\eps$ with $\phi_\eps := (\partial_t - \Delta)^{-1} \xi_\eps$, one finds that
\[
\partial_{t} v_\eps
=\tfrac{1}{2}\Delta v_\eps  + C_\eps \sin(\beta \phi_\eps) \cos(\beta v_\eps) + C_\eps \cos(\beta \phi_\eps) \sin(\beta v_\eps) \;.
\]
 \cite{MR3452276}  proved that  $C_\eps e^{\pm i \beta \phi_\eps} =C_\eps ( \cos(\beta \phi_\eps) \pm i \sin(\beta \phi_\eps))$ converges  to a nontrivial limit in $\mathcal C^{-\frac{\beta^2}{4\pi}}$,
 so  $\beta^2\in (0,4\pi)$ is precisely\footnote{Recall, for the fixed point argument to work, we must have that the regularity plus 2 is at least 1.} the regime where the above classical PDE argument applies.
 Note that the constant $ C_\eps$ can be again found by calculating the expectation of $e^{\pm i \beta \phi_\eps}$, i.e. the ``characteristic function'' of the Gaussian random variable $\phi_\eps$.

The same idea (but with a slightly different transformation than \eqref{e:DD1}) applies to the linear parabolic Anderson model in \cite{MR3358965}:
 \begin{equ}
\partial_t u_\e = \Delta u_\e  + u_\e \left(\zeta_\e  -C_\e\right)
\qquad (t,x)\in \R_+\times \T^2
  \end{equ}
 where $\zeta_\e$ is regularized {\it spatial} white noise on $\T^2$.
 With a transformation $v_\e=u_\e e^{Y_\e}$
 where $\Delta Y_\e=\zeta_\e$, one can simply check
 \[
\partial_t v_\e = \Delta v_\e
+ v_\e \left( |\nabla Y_\e|^2 -C_\e\right) - 2 \nabla Y_\e \cdot \nabla v_\e \;.
 \]
 Again, $Y_\e$ is a Gaussian process, and with $C_\e:= \E|Y_\e|^2 = -\frac{1}{2\pi} \log\e +O(1)$, $ |\nabla Y_\e|^2 -C_\e$ converges to a nontrivial limit, and the equation for $v$ is shown to be locally well-posed by standard PDE methods as above.
 (In fact \cite{MR3358965} constructed global solution on $\R_+\times \R^2$, making use of linearity of the equation.) \cite{MR3785398} studies the stochastic Schr\"odinger equation \eqref{e:RSchrE} on $\T^2$, in which a similar transformation as in \cite{MR3358965} can be applied.

This type of strategy has also been applied to stochastic hyperbolic equations. Consider the stochastic nonlinear wave equations
\[
\partial_t^2 u - \Delta u \pm u^k = \xi
\qquad
(x,t)\in \R_+\times \T^2
\]
with given initial data $(u,\partial_t u)|_{t=0}$, where $\xi$ is the space-time white noise on $\R_+\times \T^2$, and $k\ge 2$.
\cite{MR3841850} adopts the above Da Prato-Debussche trick to write $u$ as a linear part plus a remainder. Such an idea previously appears in the context of deterministic dispersive PDEs with random initial data in earlier work of McKean \cite{MR1328250} and Bourgain \cite{MR1374420}.
The proof in \cite{MR3841850} is based on a fixed point argument for the remainder equation (as above), but with the Schauder estimates replaced by Strichartz estimates for the wave equations. The key point is to use function spaces where the wave equation allows for a gain in regularity.
This gain is sufficient to prove that the remainder has better regularity than the linear solution and gives a well-defined nonlinearity for which suitable local-in-time estimates can be established.


\subsection{Regularity structures and paracontrolled distributions}
\label{sec:RS}

The above argument fails for the $\Phi^4$ equation in three spatial dimensions. As dimension increases, the space-time white noise (and thus the solution) becomes more singular.
To see how this problem rears its head, consider  the term $u^{:2:}v$ in \eqref{e:v-equ-R}.
One can show\footnote{In three spatial dimensions $\xi\in \mathcal C^{-\frac52 -}$, therefore $u\in \mathcal C^{-\frac12 -}$. From this one can show that $u^{:2:}\in \mathcal C^{-1-}$ (the rigorous proof of this fact is done by moment analysis).}
 that  in three spatial dimensions, as a space-time distribution, $u^{:2:}\in \mathcal C^\alpha$ for any $\alpha<-1$. Thus, to multiply $u^{:2:}$ with $v$, we would have to formulate the fixed point map \eqref{e:fix-pt-Phi42} for $v\in \mathcal C^\beta$ with $\beta>1$.
The product $u^{:2:}v$  would then lie in $\mathcal C^\alpha$ for any $\alpha<-1$. Unfortunately,  $(\partial_t - \Delta)^{-1}$ only provides two more degrees of regularities and thus the fixed point map \eqref{e:fix-pt-Phi42} will not bring an element in $\mathcal C^\beta$  back to the same space.

A natural idea is to go one step further in the expansion \eqref{e:DD1}. 
In view of the equation~\eqref{e:v-equ},
we define a ``second order perturbative term'' $w_\e:= (\partial_t - \Delta)^{-1} u_\e^{:3:}$, and rewrite the expansion as
\begin{equ}[e:DD2]
\Phi_\eps = u_\eps - w_\e+R_\eps \;.
\end{equ}
It turns out that one can prove that $w_\e$ converges in $\mathcal C^\alpha$ for $\alpha<1$  to a limit $w$. It remains to see whether $R_\e$ converges to a limit $R$ with even better regularity.
Using  \eqref{e:renormalizedPhi4} it is straightforward to derive an equation for $R_\eps$:
\begin{equ} [e:Picard-R]
R_\e
=
 (\partial_t - \Delta)^{-1}
\Big(
-3(u_\e^{2}-C_\e) R_\e - 3 (u_\e^{2}-C_\e) w_\e+\cdots
\Big)\;.
\end{equ}
There should be eight terms in the parenthesis, but we have only written down the two terms that are important for our discussion;
the other terms (in ``$\cdots$'') can be treated by the standard PDE argument as in the two dimensional case above.
It turns out that even after this higher order expansion \eqref{e:DD2}, the above PDE fixed point argument still does not work because of the two terms written on the right hand side of \eqref{e:Picard-R}.

For the second term, $ (u_\e^{2}-C_\e) w_\e$, we have
$u_\e^{2}-C_\e\to u^{:2:} \in \mathcal C^{-1-}$ and $w_\e\to w\in \mathcal C^{\frac12-}$. This is below the borderline of applicability of Young's theorem. It is not hard to overcome this difficulty.
 In fact, in three spatial dimensions, the term $ (u_\e^{2}-C_\e) w_\e$ requires further renormalization to converge to a nontrivial limit.
 Since this term is nothing but a convolution of several heat kernels and Gaussian noises, one can again carry out a moment analysis to find a suitable renormalization constant $\bar C_\e$ which turns out to diverge logarithmically such that
 \begin{equ}[e:u2w]
  (u_\e^{2}-C_\e) w_\e - \bar C_\e u_\e  \qquad \mbox{converges
 in $\mathcal C^{-\frac12-}$}\;.
\end{equ}
 This amounts to  renormalizing the $\Phi^4$ equation in three spatial dimensions
in the following way
\begin{equ}[e:2nd-ren-phi43]
 \partial_t \Phi_\eps = \Delta \Phi_\eps
 - ( \Phi_\eps^3 - (3C_\eps +\bar C_\e )\Phi_\eps ) + \xi_\eps \;,
\end{equ}
where $C_\eps \sim \e^{-1}$ and $\bar C_\e\sim |\log\e|$. See also Remark~\ref{rem:BPHZ}.

The first term $-3(u_\e^{2}-C_\e) R_\e$ is of the same nature as $u^{:2:}v$ in \eqref{e:fix-pt-Phi42}, so we suffer exactly the same vicious circle  of difficulty as in the discussion for $u^{:2:}v$; namely, the fixed point argument does not close.
In fact,  higher order expansions  beyond \eqref{e:DD2}
 will  always end up with such a term so the same problem will remain.
This is the real obstacle. The idea of regularity structures (which overcomes this obstacle) is that the solutions to the following two equations
\begin{equ} [e:RRt]
R=
 (\partial_t - \Delta)^{-1}
\big(  -3 u^{:2:} R
\big)\;,
\qquad
\tilde R
= (\partial_t - \Delta)^{-1}
\big(-3 u^{:2:}\big)
\end{equ}
should have the same ``small scale behavior'', because $u^{:2:}$ is  more singular than $R$  and it is the factor $u^{:2:}$ that dominates the small scale roughness.
(Here we have ignored all the other terms in \eqref{e:Picard-R},
which have better regularities than that of $u^{:2:} R $,
in order to focus on the main issue of the problem.)
This a priori knowledge that $R$  should locally look like $\tilde R$
can be formulated 
as that when the space-time points
 $z$ and $z_0$ are close, one should expect that\footnote{Eq.~\eqref{e:local-inc-R} is reminiscent to
	a Taylor expansion
	$F(z)=F(z_0)+F'(z_0)(z-z_0)+...$ where one approximates a differentiable function by Taylor polynomials. Here we approximate $R$ by $\tilde R$ which is also an object that is simply a convolution of heat kernels with white noises. Taylor polynomials are special examples of regularity structures and the theory of regularity structures is a generalization of Taylor expansion.}
 \begin{equ}[e:local-inc-R]
 R(z)-R(z_0) = R(z_0) (\tilde R(z)-\tilde R(z_0)) + O(|z-z_0|) \;.
 \end{equ}
Namely, the  local increment of $R$ is approximately the same as the  local increment of $\tilde R$ -- up to multiplying a factor $R(z_0)$ which depends on the base point $z_0$; the reason that this multiplicative factor should be $R(z_0)$  is clear from the structure of Eq.~\eqref{e:RRt}.

Since $\tilde R$ is again a {\it concrete} object which is simply convolutions of heat kernels with powers of Gaussians,
it is easy to prove $\tilde R\in \mathcal C^{1-}$ by analyzing its moments as before.
Thus if $R$ satisfies \eqref{e:local-inc-R}, that is, locally looks like $\tilde R$,  one  has $R\in \mathcal C^{1-}$ as well.
The converse is not true; the set of $R$ satisfying  \eqref{e:local-inc-R}  is a strictly smaller set than $\mathcal C^{1-}$.
The key is to formulate a fixed point problem in the space of all functions $R$ that have ``prescribed local expansion'' \eqref{e:local-inc-R} (rather than in standard function spaces such as $\mathcal C^\alpha$).
The aforementioned vicious term $u^{:2:} R$ which could not be defined 
for arbitrary $R\in \mathcal C^{1-}$, can now be defined if $R$ locally looks like $\tilde R$,
because $u^{:2:} \tilde R$ is again simply a concrete combination of Gaussian processes and heat kernels!
It turns out that the fixed point argument closes in the space of functions  having such ``prescribed local expansions'', and the fixed point $R$ together with \eqref{e:DD2} yields the solution to the $\Phi^4$ equation in three dimensions.

The above idea of solving stochastic equations  in a space of functions or distributions that have prescribed local approximations by certain canonical objects to some extent had precursor in the simpler setting of stochastic ordinary differential equations, which is called rough path theory, see \cite{MR1654527} or the book \cite{MR3289027}, in particular a formulation by \cite{MR2091358}.
Constructing the solution to the $\Phi^4$ equation on three dimensional torus was
the first example of the theory built in \cite{Hairer14}. The review articles 
\cite{MR3336866,hairer2015regularity} have more detailed pedagogical explanations on the theory and the application to this equation.


\begin{rmk}\label{rem:BPHZ}
We have found the renormalization
of $u_{\eps}^{2}$, $u_{\eps}^{3}$ and $ u_\e^{2} w_\e$
and proved convergence of the renormalized objects
by moment analysis. Analyzing moments of these random objects
are the only probabilistic component of  \cite{Hairer14}.
As the equation in question becomes more singular,
the number of such random objects to be studied increases,
and it is tedious or even impossible to analyze each of them by hand.
\cite{chandra2016analytic} develops a ``blackbox'' that provides systematic and automatic treatment for renormalization and moment analysis for these perturbative objects arising from general singular SPDEs.
Moreover, there are also algebraic aspects for the renormalization procedure (so called ``renormalization groups''),
which has been systematically treated in \cite{bruned2016algebraic}.
Finally, there is a question regarding how the renormalized equation (e.g. \eqref{e:2nd-ren-phi43})
will look like after  renormalizing these random objects,
and this is answered in \cite{bruned2017renormalising}.
\end{rmk}

Hairer's theory has been applied to provide solutions to other very singular SPDEs, for instance, a  generalized parabolic Anderson model (a generalization of \eqref{e:gPAM})
 \[
 \partial_t u = \Delta u + \sum\nolimits_{ij} f_{ij}(u) \partial_i u \partial_j u+ g(u)\xi
 \]
 where $f$ and $g$ are sufficiently regular functions.
The well-posedness of the  KPZ equation in one spatial dimension 
was solved in \cite{Hairer13} using 
 the theory of ``controlled rough paths'' \cite{MR2091358}, which can now be viewed as 
 a special case of regularity structures; see the book \cite{MR3289027} for rough paths, regularity structures and applications to KPZ.
 
Other applications include (but are not limited to) the stochastic Navier-Stokes equation \eqref{e:SNS} with white noise on the three-dimensional torus \cite{Zhu2014NS}, 
the stochastic heat equation with multiplicative noise
\eqref{e:NSHE} \cite{MR3417505,MR3779690},
the dynamical sine-Gordon equation \eqref{e:sine-G} on two-dimensional torus for arbitrary $\beta^2 < 8\pi$ \cite{MR3452276,sineGordonwholeregime},
the  stochastic quantization of abelian gauge theory / stochastic gauged Ginzburg-Landau equation by \cite{Shen2018Abelian},
and the random motion of string in manifold \eqref{e:Hairer-String} \cite{hairer2019geo}.

Besides Hairer's theory \cite{Hairer14}, some alternative methods have also been introduced.
The {\it paracontrolled distribution} method of Gubinelli-Imkeller-Perkowski \cite{MR3406823} is based on a similar idea of controlling the local behavior of solutions, but is implemented in a different way by using the Littlewood-Paley theory and paraproducts \cite{MR2768550}.
See \cite{MR3828162} for a review on paracontrolled distribution.  
The paracontrolled distribution method has been also successfully applied to, for instance, the KPZ equation \cite{GubPerk,MR3769661} in $d=1$ (and more recently the construction of solution on the entire real line instead of a circle \cite{perkowski2018line}), 
a multi-component coupled KPZ equation \cite{MR3653951},
the $\Phi^4$ equation \cite{MR3846835} in $d=3$ (and more interestingly, its global solution by \cite{MR3719541}).

The paracontrolled distribution method has not only allowed to prove well-posedness results for stochastic PDE, but also resulted in the
 construction of other singular objects which could not be made sense before. 
For instance, \cite{AllezChouk} constructed the
Anderson Hamiltonian (i.e. Schr\"odinger operator) on the two dimensional torus, formally defined as
$\mathscr H=-\Delta+\xi$ where $\xi$ is a singular potential such as white noise.
As another example, \cite{MR3785598} proved existence and uniqueness of solution for stochastic 
{\it ordinary} differential equations $dX_t =V(t,X_t)dt +dB_t$ with {\it distributional} valued drift $V$ where  $B$ is a $d$-dimensional Brownian motion, and this is achieved via the study of the
generator of the above stochastic 
 ordinary differential equations  given by
$\partial_t+\frac12 \Delta+V\cdot \nabla$.
\cite{MR3785598}  also managed to make sense of a
a singular ``polymer measure''  on the space of continuous functions  formally given by $\mathbb Q_T(d\omega)=Z_0^{-1}\exp\left(\int_0^T\xi(\omega_s)ds\right)\mathbb W_T(d\omega)$
where $\mathbb W$ is the Wiener measure (i.e. the Gaussian measure for Brownian motions) on $C([0, T ], \R^d)$ for $d = 2, 3$, $\xi$ is a spatial white noise on the d-dimensional torus  independent of $\mathbb W$, and $Z_0$ is an (infinite) renormalization constant.

We will discuss another application of the paracontrolled distribution method  on scaling limit problem 
with a bit more details in Section~\ref{phy-PAM}.

In the line of this paracontrolled distribution approach,
\cite{MR3475460} provided a semigroup approach,
and has been applied to the generalized parabolic Anderson model \eqref{e:gPAM}
on a potentially unbounded 2-dimensional Riemannian manifold.

Another method based on renormalization group flow was introduced by Kupiainen \cite{Antti},
which for instance has been applied to prove local well-posedness for a generalized KPZ equation \cite{MR3607594} introduced by H. Spohn \cite {MR3176405} in the context of stochastic hydrodynamics.

With all these alternative methods, the theory of regularity structures is by far the most systematic and  general approach; for instance 
 it has developed the ``blackbox theorems'' as mentioned in Remark~\ref{rem:BPHZ}
 which makes the implementation of this theory very ``automatic'',
and 
it can deal with equations which are extremely singular (that is, very close or even arbitrarily close to ``criticality'', see Remark~\ref{rmk:subcritical}) such as  the random string in a manifold \eqref{e:Hairer-String} or the dynamical sine-Gordon  \eqref{e:sine-G}  for arbitrary $\beta^2 < 8\pi$.

\begin{rmk}\label{rmk:subcritical}
(Subcriticality of stochastic PDE.)  The methods developed in  \cite{Hairer14}, \cite{MR2768550} and \cite{Antti} are all for  ``subcritical'' semilinear stochastic PDEs. For stochastic PDEs with white noise,
the equation being subcritical means that the nonlinear term has better regularity than the linear terms;
namely, small scale roughness is dominated by the linear solution.
For instance, for the $\Phi^4$ equation in three spatial dimensions, the term $\Phi^3$ has regularity $-\frac32-$ while $\Delta\Phi$ and $\xi$ have regularities $-\frac52-$. 
Subcriticality often depends on spatial dimensions: KPZ, \eqref{e:NSHE} and \eqref{e:Hairer-String}
are subcritical in $d<2$, while $\Phi^4$, parabolic Anderson model \eqref{e:gPAM}, Navier-Stokes equation \eqref{e:SNS} with space-time white noise and the stochastic  Yang-Mills heat flow \eqref{e:YM-heat} are subcritical in $d<4$.
The dynamical sine-Gordon equation \eqref{e:sine-G} however, is subcritical for $\beta^2 < 8\pi$.

The stochastic PDEs being discussed here in supercritical regimes (i.e. above the aforementioned criticalities)
are not expected to have nontrivial meanings of solutions. We only expect to get Gaussian limit,
although the Gaussian variances may be nontrivial; the reader is referred to \cite[Theorem~1.1]{MR3790153} for a flavor of such a result for the KPZ equation in $d\ge 3$.

Critical dimensions are much more subtle. We refer to \cite{ChatterjeeDunlapKPZ,MR3719953,Caravenna2Dentire,Gu2018KPZ}
for the very new progress on the KPZ equation in $d=2$.
\end{rmk}

\begin{rmk}
We remark that although we have focused on semilinear equations in our expositions, 
the methods developed in \cite{Hairer14}, \cite{MR2768550} have also extended to quasilinear equations, see \cite{otto2018quasilinear,MR3916943,MR3916262,GerencserQuasilinear}.
\end{rmk}

\subsection{A brief discussion on weak solutions}
\label{sec:weak-sol}

The solutions to SPDEs that we have discussed so far are called ``strong solutions'', as opposed to the ``weak solutions'' that we will now briefly discuss\footnote{Not all equations which admit weak solutions admit strong solutions. A famous stochastic differential equation
example is called Tanaka's equation, see \cite[Example~3.5]{MR1121940}.}.
Let us immediately point out that the ``weak solutions'' in the stochastic context have nothing to do with the weak solutions in deterministic PDE theory; one sometimes adopts the terminology ``probabilistically weak solutions''.

For a strong solution, one starts with a probability space on which the noise is defined and then builds a mapping from that probability space and the initial data space to a space of functions (or distributions) which satisfies the prescribed equation\footnote{As we saw, making sense of what
it means to satisfying the equation often takes significant work an involves regularizations and renormalizations. There are also some measurability assumptions which should be imposed on strong solutions so future noise cannot effect the evolution before its time.}
with probability one (i.e., for almost every point in the probability space).
Though subtle, it is important to understand that a strong solution to an SPDE need not be function valued (as we saw, in some instances it is distribution valued, living in some spaces of negative regularity).

For a standard PDE, a weak solution requires that the equation holds when tested against a suitable class of functions.
For SPDEs, the analog of this involves treating solutions statistically as probability measures on the solution space, rather than as random variables supported on the probability space on which the noise is defined.
Roughly, a weak solution means that we can define some noise (with the right distribution and measurability assumptions satisfied) so that the canonical process\footnote{The canonical process is the random variable whose probability space is defined as the solution space equip with the probability measure
of the proposed solution.} on the solution space, along with the noise satisfying the desired equation\footnote{Here is a, hopefully, more intuitive explanation for this different notions of strong versus weak. Imagine that human life were governed by an SPDE.
Then, a strong solution would tell us how each individual's life would unfold given knowledge of all of the randomness which befalls them, in addition to the world around them. A weak solution is statistical---it tells us that people with certain characteristics,
have certain probabilities of having their life unfold in various ways. Given such a prescription of probabilities, how can be verify that this is, indeed, a weak solution to the ``SPDE of life''? Well, we need to demonstrate that there exists randomness which would, in fact, result in the
aforementioned probabilities. Then, we would need to verify that the randomness is distributed in the way that the SPDE of life claims (e.g. space-time white noise). While this is all a bit tongue and cheek, hopefully it helps explain the difference.}.
As we will see, ``martingale problems'' provide a very convenient way to demonstrate that a weak solution solves an SPDE (instead of demonstrating the existence of a suitable noise as above).

Let us illustrate these ideas in the simplest setting of stochastic differential equations (SDE). Let $\xi(t)$ denote temporal white noise. Like its space-time counterpart, this can be defined in various ways (e.g. as a series with random coefficients). Consider the SDE $dX_t = \xi(t) dt$, which in
integrated form reads $X_t = X_0 + \int_0^t \xi(ds)$ (let us assume that $X_0=0$ for simplicity). Once integration with respect to white noise is understood, this defines a solution map (and hence a strong solution) from $(\xi,X_0)$ to the full trajectory of $X_t$ for $t\geq 0$.
One checks that $t\mapsto X_t$ is continuous and that its marginal distributions are Gaussian with covariance of $X_s$ and $X_t$ given by the minimum of $s$ and $t$.
This, in fact, implies that the distribution of the function $t\to X_t$ is Wiener measure---that is, the distribution of Brownian motion.
If instead of $X_t$ we had another Brownian motion $\tilde{X}_t$ (for instance, we could have $\tilde{X}_t = -X_t$ or just an independent Brownian motion), then $\tilde{X}_t$ would be a weak solution, but not a strong solution. This is because $X_t$ and $\tilde{X}_t$ have the same distribution, even
if they are not ``driven'' by the same noise.

The ``martingale problem'' provides an alternative characterization to the Gaussian description above for Brownian motion.
The L\'{e}vy characterization theorem says that $X_t$ is distributed as a Brownian motion if it is almost surely continuous and both $X_t$ and $X_t^2-t$ are martingales\footnote{In fact, ``local'' martingales.}
A measure on $X_t$ that satisfies this is said to satisfy the ``martingale problem'' characterizing Brownian motion.

What does it mean that $X_t$ (or $X_t^2-t$) is a ``martingale''?  Roughly speaking, this means that given the history of $X_t$ up to time $t$, the expected value of its future location is exactly $X_t$.
This is like a fair gambling system in which your future expected profit is always zero. Martingales are essentially a particular class of centered noise.

Martingale problems exist for general classes of SDEs and are often very useful for proving convergence results.
For instance, to show that a discrete time Markov chain converges to an SDE (e.g. a random walk converges to Brownian motion), one can demonstrate that the discrete chain satisfies a discrete version of the SDE's martingale problem.
Then, provide one can demonstrate compactness of the measures (on the evolution of the Markov chain), all limit points must satisfying the limiting SDE's martingale problem. This generally proves uniqueness of the limit points, and hence convergence.

Weak solutions to linear SPDEs can also be characterized  in terms of martingale problems. Let us describe how this works for the multiplicative noise stochastic heat equation \eqref{e:SHE}, recalled here\footnote{This equation also admits a strong solution which can be written as a ``chaos series'' of
multiple stochastic integrals against space-time white noise $\xi$.}:
 \begin{equ}[e:SHE2]
\partial_t u(t,x) - \Delta u(t,x)  = u(t,x)\xi (t,x)
\qquad
(t,x)\in \R_+ \times \R
\;.
  \end{equ}
Let us write $u_t(x) = u(t,x)$ and think of $u$ as a measure on $C(\R_+,C(\R))$ (continuous maps from $t\in \R_+$ to continuous spatial functions). For any test function $\varphi\in C^\infty(\R)$, write $u_t(\phi)=\int u_t(x)\phi(x)dx$.
With this notation, define the processes
\begin{equ}[e:N-and-Q]
	N_t(\varphi )
	:=
	u_t(\varphi )- u_0(\varphi)
	-  \int_0^t u_s(\Delta \varphi)\,ds
\qquad \mbox{and}\qquad
	Q_t(\varphi)
	 :=
	N_t(\varphi)^2 - \int_0^t (u_s^2,\varphi^2 )\,ds.
\end{equ}
We say that $u$ satisfies the martingale problem for the multiplicative noise stochastic heat equation if both $N_t$ and $Q_t$ are (local) martingales for all test functions $\phi$. Any $u$ that satisfies this is a weak solution, see for instance \cite[Definition~4.10]{Bertini1997}. Just as martingale problems are
useful in proving convergence of Markov chains to SDEs, so too can they be used in SPDE convergence proofs---see Section~\ref{sec:ASEP} and Section~\ref{sec:6V} for some examples where this type of martingale problem has been used for such a purpose.

It is generally hard to formulate a martingale problem characterization for weak solutions to singular nonlinear SPDEs. For the KPZ equation (in one spatial dimension) $H$, one can use the Hopf-Cole transform and define $H$ as a ``Hopf-Cole solution''
to the KPZ equation if $u=e^H$ is a solution to the multiplicative noise stochastic heat equation. This notion of solution agrees with those discussed earlier in this text. However, such linearizing transformations are uncommon and this should be
thought of as a rather useful trick, but not a general theory.

Remarkably, for the stochastic Burgers equation (which is formally the equation satisfied by the spatial derivative of the KPZ equation $\partial_x H$)
  \begin{equation}\label{eq:SBE}
     \partial_t u = \tfrac{ 1}{2}\partial_x^2 u - \tfrac{1}{2} \partial_x (u^2) +\partial_x \xi,
  \end{equation}
\cite{GJ14a} found a way to formulate a martingale problem characterization and \cite{GP2015a} (with a slightly improved formulation) matches the solution to this martingale problem  to the Hopf-Cole solution (see \eqref{e:SHE}) whereby showing uniqueness
of  the solution to this martingale problem; they call it an {\it energy solution}. 
There were some limitations of this notion of solution, namely it only works for particular types of initial data;
very recently,  however,
\cite{MR3828193} has generalized the notion of energy solution to system configurations with finite entropy with respect  to stationarity; and
 \cite{Yang} has extended this method to include more general initial data such as flat initial data.
It has proved to be quite useful in demonstrating convergence results, as we explain later in Section~\ref{sec:energy}.
Finally, let's mention that very recently \cite{gubinelli2018infinitesimal} developed
a martingale approach for a class of singular stochastic PDEs of Burgers type, including fractional and multi-component Burgers equations.

Let us end this discussion by mentioning (without any explanation) another powerful approach to defining weak solutions of SPDEs---Dirichlet forms. For example, for the $\Phi^4$ equation in two spatial dimensions, before Da Prato and Debussche constructed
  their strong solution in \cite{MR2016604}, the paper \cite{AlbRock91} constructed a weak solution via Dirichlet forms (which involves significant functional analysis). For comprehensive discussion on this topic we refer to the book \cite{Fukushima2011} and the references therein.

\section{Renormalization in physical models}
\label{sec:Renormalize-physical}

Let us take stock of what we have learned so far. In Section \ref{sec:linear} we observed that (at least in the linear case) SPDEs arise as scaling limits for microscopic models of physical systems. In Section \ref{sec:NSPDE} we introduced a number of nonlinear SPDEs and claimed that they model
various interesting physical systems. However, before trying to justify that claim, we had to confront the challenge of well-posedness. Namely, how to make sense of what a ``solution'' means to these equations. Section \ref{sec:Well-posedness} surveyed the main techniques for doing this.

In defining a ``solution'' we mollified (or smoothed) the noise and defined the solution through a limit transition.
From that perspective, it is reasonable to hope that the same methods can be applied to other types of regularizations of the noise, or equation---for instance to show that discrete systems converge to the SPDEs.
We will address this further in Section \ref{sec:universal}.

In Section \ref{sec:Well-posedness} we found that besides regularizing the noise, we also had to introduce certain ``renormalizations'' to our equations for them to admit limits. At first glance this tweaking of the equation seems a bit crooked.
In this section we will explain how these renormalizations have concrete physical meaning, thus justifying our definitions.
For instance, a diverging renormalization constant may relate to a tuning for a microscopic system of the overall scale, reference frame, temperature, or other physically meaningful parameters.
We will focus our discussion on two systems: the dynamic $\Phi^4$ equation and parabolic Anderson model. For the KPZ equation, we save this discussion until Section \ref{sec:universal} where we also highlight the notion of universality. 
%

\subsection{$\Phi^4$ equation}
\label{sec:phy-Phi4}
Let us consider the example of a magnet near its critical temperature in Section~\ref{sec:NSPDE}. Many mathematical models have been proposed to describe various behaviors of magnetic systems. Here we investigate one particular example called the ``Kac-Ising model''\footnote{The Ising model was introduced in 1920 by Lenz and named after his student Ising who showed that in one spatial dimension it did not admit any phase transition. That original model involves nearest neighbor pair interactions of $-1$ and $+1$ spins. There are many generalizations of this model besides the long-range
Kac-Ising model we consider here. For instance, one can consider the ``higher spin'' versions with $q$ different types of spins and interactions which award equal spins along edges. This is known as the $q$-state Potts model.}.

The model in two dimensions is defined as follows. Denote by $\Lambda_{N}=\Z^d/(2N+1)\Z^d$  a large two-dimensional discrete torus of ``radius'' $N=\e^{-1}$ (we introduce $\e$ here to later take to zero), which represents the space in which our magnetic material lives.
Each site $k\in \Lambda_{N} $ is decorated with a spin $\sigma(k)$ which for simplicity is assumed to take values in  $\{ -1, +1 \}$.
 Denote by $\Sigma_{N} = \{ -1, +1 \}^{\Lambda_{N}}$ the set of spin configurations on $\Lambda_{N}$. For a spin configuration $\sigma = (\sigma(k), \, k \in \Lambda_N )\in \Sigma_{N}$ we define the ``Hamiltonian'' as
\begin{equ}\label{e:Hamiltonian}
H(\sigma) := -  \frac12 \sum_{k, j \in \Lambda_{N}}  \kappa_{\gamma}(k-j) \sigma(j)   \sigma(k)  \;,
\end{equ}
where, for $\gamma\in (0,1$), $\kappa_{\gamma}(k)$ is a non-negative function\footnote{A concrete choice of the interaction kernel
$\kappa_\gamma$ is to set $\kappa_\gamma(k) = c_\gamma \, \gamma^d \, \mathfrak R(\gamma k) $
where $\mathfrak R\colon \R^d \to \R$ is a smooth, nonnegative function with compact support and $ c_\gamma$ is chosen to ensure that $\sum_{k \in \Lambda_N} \kappa_\gamma(k)=1$.} supported on $|k|\le \gamma^{-1}$ which integrates to $1$.
%
%
Then for any ``inverse temperature'' $\beta>0$ we can define the ``Gibbs measure'' $\lambda_{\gamma}$ on $\sigma\in \Sigma_{N}$ as
\begin{equation*}
\lambda_{\gamma} (\sigma) := \frac{1}{Z(\beta)}\exp\Big( - \beta H(\sigma) \Big)\; ,
\end{equation*}
where $Z(\beta) := \sum_{\sigma \in \Sigma_{N}} \exp\big( - \beta H(\sigma) \big) $  makes $\lambda_{\gamma}$ into a probability measure.

The measure $\lambda_{\gamma}$ is known as a equilibrium measure since the probability of finding a configuration $\sigma$ is proportional to the exponential of the energy of that configuration. It is also known as equilibrium because it arises as the equilibrium (or stationary, or invariant) measure for various simple, local stochastic dynamics on the configuration space. We will consider one such example known as ``Glauber dynamics'' \cite{doi:10.1063/1.1703954}. For $j\in\Lambda_N$, let $\sigma^j \in \Sigma_N$ denote the spin configuration that coincides with $\sigma$ except the spin at position $j$ is flipped:  $\sigma^j(j) = -\sigma(j)$.
The Glauber dynamic is the following continuous time Markov process: For each $j\in \Lambda_N$, the configuration $\sigma$ is updated to $\sigma^j$ at rate\footnote{For those not used to continuous time Markov processes, let us explain this more precisely. Starting from some configuration $\sigma$, to each $j\in \Lambda_N$ we associate independent exponentially distributed random variables of rate $c_{\gamma}(\sigma,j)$ (i.e., with mean inverse to the rate). One then compares all of these random variables and for the $j$ whose random variable is minimal, the configuration updates to $\sigma^j$. The time at which this occurs is the value of the associated random variable. From that point on, one repeats the whole story, choosing new exponential random variables with rates given by the updated $\sigma^j$. Due to the ``memoryless'' property of exponential random variables (i.e., for $X$ exponential of rate 1, conditional on $X>x$, the law of $X-x$ is that of a rate 1 exponential random variable), this constructs a continuous time Markov process.} $c_\gamma(\sigma, j)$ where
\begin{equ}[e:def-c-gamma]
c_\gamma(\sigma, j) := \frac{\lambda_{\gamma}(\sigma^j)}{\lambda_{\gamma}(\sigma) + \lambda_{\gamma}(\sigma^j)}\;.
\end{equ}
Once an update occurs for some $j$, all rates are recalculated relative to the new configuration. It is standard to show that the measure $\lambda_{\gamma}$ is the unique invariant measure for the Glauber dynamic, meaning that for any starting measuring, eventually the measure will converge in distribution to $\lambda_{\gamma}$. Likewise, if started according to $\lambda_{\gamma}$, then the distribution at any later time will still be distributed according to that measure.

The following picture illustrate the dynamics where, for each fixed time $t$, one has a spin configuration $\sigma(t)$. We would like to observe a scaling limit of the system from
a large distance scale of $\e^{-1}$ and long time scale
of $\alpha^{-1}$. At larger scales, the  $\sigma$ oscillating between $\pm 1$ would yield a field in a very weak topology;
it will be more convenient to consider an averaged field\footnote{Using the interaction kernel $\kappa_\gamma$ to average out the field $\sigma$ is merely a matter of convenience, and it will lead to a clean form of \eqref{e:Ising-Phi4}. Also, convergence of $\sigma$  follows a fortiori.}
\begin{equ}[e:def-h-gamma]
h_\gamma (\sigma(t),k):= \sum_{\ell \in \Lambda_N} \kappa_\gamma(k - \ell) \sigma(t,\ell).
\end{equ}
We may also abuse notation and write $h_\gamma(t,k) = h_\gamma (\sigma(t),k)$, suppressing the explicit dependence on $\sigma(t)$.
\begin{center}
\includegraphics[scale=0.3]{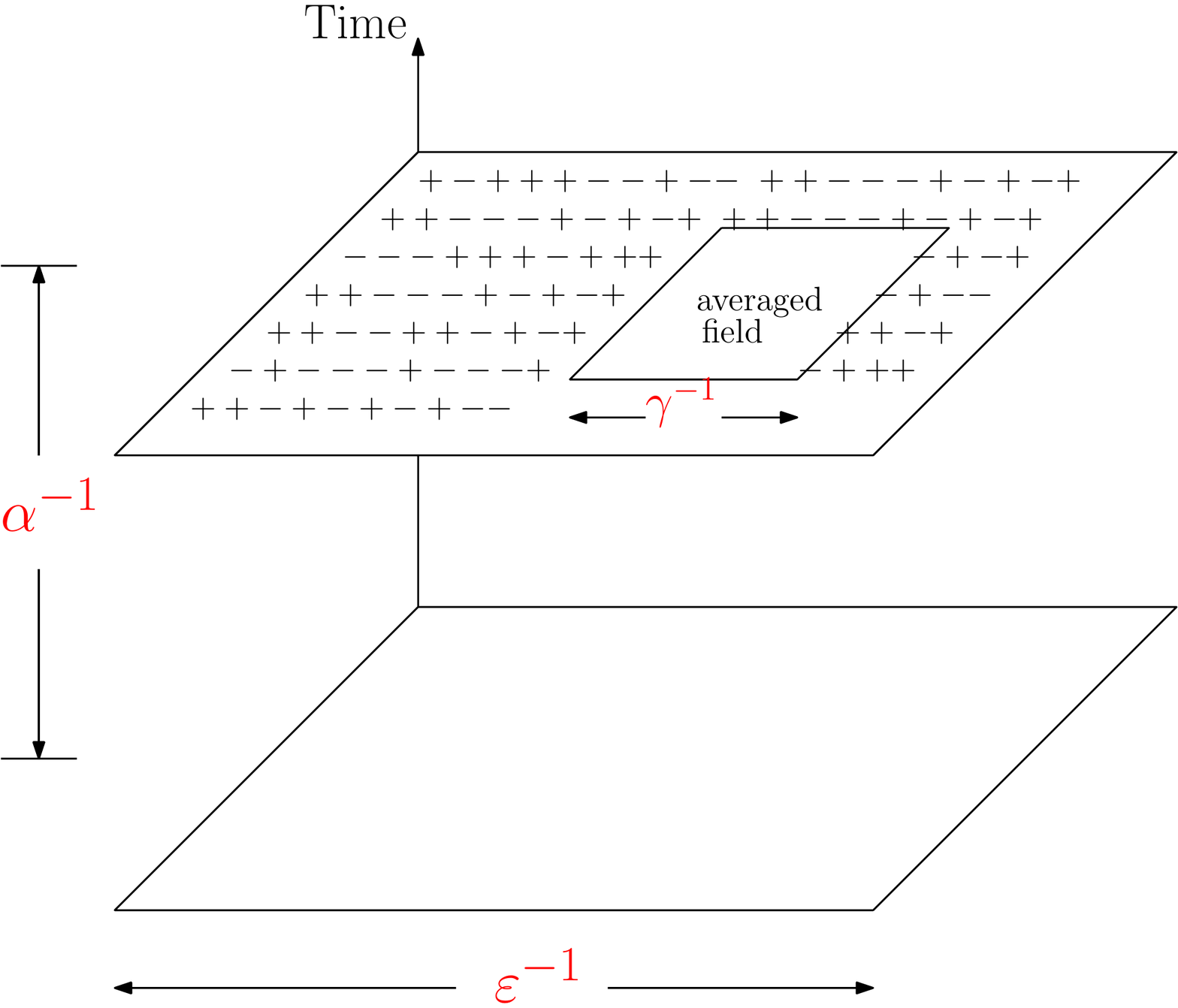}
\end{center}

\begin{rmk}\label{rem:Markov process}
In order to prove an SPDE limit  we first write down a discretized SPDE. Generally, this involves a coupled system of stochastic differential equations driven by martingales (recall the brief discussion from Section \ref{sec:weak-sol}). Without delving deep into details of theory of Markov processes, let us illustrate this with a simple example. Our applications of this general idea will be more involved, though we will avoid going into details there. Consider a continuous time random walk $x(t)$ on $\Z$ where a left jump (by 1) occurs at rate $\ell$ and a right jump (by 1) occurs at rate $r$. For any function $f:\Z\to \R$, the expected value $u(t,y)= \E_y\big[f\big(x(t)\big)\big]$ (where $\E_y$ means the expected value assuming initial data $x(0)=y$) satisfies the system of ODEs $\tfrac{\partial}{\partial t} u(t,y) = Lu(t,y)$ where the operator $L$ acts on functions $g(y)$ as $\big(Lg\big)(y) := \ell \big(g(y-1)-g(y)\big) + r \big(g(y+1)-g(y)\big)$. Without taking the expectations, $f\big(x(t)\big)$ will satisfy the evolution equation
$$
d f\big(x(t)\big) = \big(Lf \big)\big(x(t)\big)dt + dm(t)
$$
where the first term is a ``drift'' and the extra term $dm(t)$ is the time derivative of a martingale. The martingale $m(t)$ can be explicitly described in terms of what are called compensated Poisson processes; or equivalently in terms of its ``quadratic variation''. Under the diffusive scaling which takes $x(t)$ to a Brownian motion, the martingale converges to the time derivative of Brownian motion (what we called 1-dimensional white noise earlier in Section \ref{sec:weak-sol}).  This can be shown, for example, via the martingale problem for Brownian motion. For more complicated Markov processes, the story is analogous, albeit the analogs of $L$ (known as the instantaneous generator) and $m(t)$ are more complicated.
\end{rmk}

There are three steps to find an SPDE limit for the Glauber dynamic on the Kac-Ising model:

\noindent \underline{\it Microscopic evolution}.
In the spirit of Remark~\ref{rem:Markov process} we may write down the evolution the Markov process $h_\gamma$ \eqref{e:def-h-gamma} in terms of a drift part (that we make explicit) and a martingale part $m$ (that we do not precisely describe below):
\begin{equ}[e:dyn-h-gamma]
dh_\gamma (\sigma(t),x) =
\sum_{j \in \Lambda_N} c_\gamma(\sigma(t),j) \big(h_\gamma (\sigma^j(t),x) -h_\gamma (\sigma(t),x) \big)\,dt + dm (\sigma(t),x)\;.
\end{equ}
Recall that $c_\gamma$ is given in \eqref{e:def-c-gamma}. We may  look for scaling so that
 $X_\gamma(t,x): = \frac{1}{\delta} h_\gamma ( \sigma(\alpha^{-1}t), \e^{-1} x)$ converges to a limiting (nonlinear) SPDE.
After Taylor expanding the nonlinear dynamic \eqref{e:dyn-h-gamma}
into polynomials in $h_\gamma$ and passing to $X_\gamma$
we get the  following discrete equation
\begin{equ}[e:Ising-Phi4]
dX_\gamma
=\Big(\frac{\varepsilon^{2}}{\gamma^{2}}\frac{1}{\alpha}\Delta_\gamma X_\gamma
+\frac{\beta-1}{\alpha} X_\gamma
-\frac{\delta^{2}}{\alpha} X^{3}+...\Big)\,dt
+dM\;.
\end{equ}
Here $\Delta_\gamma$ is a difference operator (based on the kernel $\kappa_\gamma$) which is  approximately the  Laplacian. The martingale $M$ is a rescaled version of $m$, whose quadratic variation can be also explicitly calculated. 
We omit this, though note that  one should think of the noise term $dM$ as being of order $O(\frac{\eps}{\delta \sqrt{\alpha}})$, in the sense that the quadratic variation of $M$ is of order $O(\frac{\eps^2}{\delta^2 \alpha})$.

\noindent \underline{\it Scaling}.
We may now seek suitable choices for $\e,\alpha,\gamma,\delta$ such that the dynamics converge to a limit (the $\Phi^4$ equation). In particular, we choose  $\e \approx \gamma^{2}$, $\alpha\approx \gamma^{2}$ and $\delta \approx \gamma$ so that the Laplacian, cubic and martingale all terms balance:
\[
dX_\gamma
=\Big(\overbrace{\frac{\varepsilon^{2}}{\gamma^{2}}\frac{1}{\alpha}}^{O(1)} \Delta_\gamma X_\gamma
+\!\!\!\!\! \overbrace{\frac{\beta-1}{\alpha}}^{\mbox{\scriptsize
	tune $\beta\approx \beta_c =1$}} \!\!\!\!\!\!X_\gamma
- \overbrace{\frac{\delta^{2}}{\alpha}}^{O(1)} X_\gamma^{3}+...\Big)\,dt
+\overbrace{dM}^{O(1)} \;.
\]
The ``critical''  inverse temperature (at which the magnet loses magnetization) $\beta_c =1$ is precisely the value such that sending $\beta\to \beta_c$ at suitable rate will suppress the linear term from blowing-up.
Tuning $\beta=1+ o(\alpha)$, at  large scales (i.e., $\gamma\searrow 0$ with $\e,\alpha,\delta$ as above) one would expect that $X_{\gamma}$ converges to the solution to the equation
\footnote{One might be puzzled  why we can obtain this limiting equation that is not scaling invariant (recall Remark~\ref{rmk:scaling})
via a scaling limit procedure. The reason is that the interaction range $\gamma^{-1}$ in the model \eqref{e:Hamiltonian} sets an additional scale.}
\[
\partial_t \Phi = \Delta \Phi - \Phi^3 + \xi \;.
\]
This, however, is not the case. $X_{\gamma}^{3}$ only admits a nontrivial limit when suitable renormalized. 

\noindent \underline{\it Renormalization}.
Recalling the discussion in Section~\ref{sec:Well-posedness} (in particular Eq.~\eqref{e:renormalizedPhi4}), the correct way of taking the limit is to subtract a renormalization term $C_\gamma X_\gamma$ after the cubic term where $C_\gamma$ is the suitably divergent constant. We then add this term back in the linear term so that the equation remains unchanged
 \[
dX_\gamma
=\Big(\frac{\varepsilon^{2}}{\gamma^{2}}\frac{1}{\alpha}\Delta X_\gamma
+ \Big(\frac{\beta-1}{\alpha}  -\frac{\delta^2 C_\gamma}{\alpha} \Big)X_\gamma
-\frac{\delta^{2}}{\alpha} \left(
	X_\gamma^{3} - C_\gamma X \right)+...\Big)\,dt
+dM \;.
\]
In two spatial dimensions $C_{\gamma}$ diverges logarithmically according to the calculation in Eq~\eqref{e:Eu2}. Scaling $\e,\alpha,\gamma,\delta$ as above, but now  tuning the inverse temperature in the {\it correct}  way $\beta=1+\delta^2 C_\gamma+ o(\alpha)$,
so that $\beta\to \beta_c=1$ but with a slightly modified rate, $X_\gamma$ converges to the solution to the  renormalized $ \Phi^4$ equation. This is proved in \cite{MR3628883}\footnote{A similar result in one spatial dimension was shown in the nineties in \cite{MR894407,MR1358083}. In one or two spatial dimensions the $ \Phi^4$ equation is not the only possible limit. \cite{MR3820327} considered a generalized ferromagnetic model (called the Kac-Blume-Capel model) and proves that the Glauber dynamic converge to either the $\Phi^4$ or $\Phi^6$ equations in various regimes. In three spatial dimensions the $ \Phi^4$ equation is believed to be the only nontrivial SPDE limit one can obtain from ferromagnetic models, though the proof is still open.}.

We emphasize two important points: First, the model \eqref{e:Hamiltonian}  is an interpolation between the standard nearest-neighbor Ising model and a mean-field Ising model (also called the Curie-Weiss model) where all sites interact equally with each other. These two extreme cases have rather different behaviors than the Kac-Ising model. For instance, limits of the nearest-neighbor Ising model lead to conformal invariant objects; see, for example, \cite{MR3296821} and references therein. On the other hand, the mean-field Ising model with interaction length of order $N$ has Gaussian fluctuations for the magnetization. The $\Phi^4$ equation arises in an intermediate or ``cross-over'' scale where it is possible to balance all desired terms in \eqref{e:Ising-Phi4}.

The second point is that the renormalization constant in this model represents the delicate rate at which the temperature approaches criticality. Instead, if we did not tune $\beta$ to $\beta_c$ properly, the (averaged) magnetic field would either become deterministic (when $\beta >\beta_c$) or completely disordered (like smoothed white noise (when $\beta<\beta_c$). It is only in this critical scaling window that we see the balance between these two ordered and disordered phases.

\subsection{Parabolic Anderson equation}
\label{phy-PAM}

We turn to study the parabolic Anderson model (PAM), which describes population dynamics. 

\noindent \underline{\it Microscopic evolution}. Consider the discrete PAM for $v:\R_+\times \Z^d \to \R$:
\begin{equ}\label{e:dPAM}
   \partial_t v (t, x) = \Delta_{\rm d} v (t, x)
   +  v (t, x) \eta (x)\;.
\end{equ}
Here  $\Delta_{\rm d}$ can be taken as the generator of a general symmetric random walk, but
for simplicity we will assume it is the discrete Laplacian on $\Z^d$ (i.e., the generator of a simple symmetric random walk). We also assume for simplicity that $\{\eta (x) : i \in \Z^d\}$ is family of independent and identically distributed mean-zero Gaussian random variables (though in general
this randomness does not have to be Gaussian distributed).

The PAM \eqref{e:dPAM} models random walks which can {\it branch or die} in a given random environment $\eta$. These are particles on the lattice $\Z^d$ which all independently follow the dynamics of $\Delta_{\rm d}$-generator  random walk, and
which at each lattice point $x\in \Z^d$ get killed with rate $\eta(x)^-= \max(-\eta(x), 0)$ and branch into two particles with rate $\eta(x)^+= \max(\eta(x), 0)$; after the branching the two particles follow the same dynamics. All particles evolve, branch and die independently of each other. The function $v(t,x)$ is the expected number of particles at time $t$ in location $x$, conditioned on the random environment $\eta$. This model is used to study, for example, the population evolution for microorganism which flourish in regions with high concentration of nutrition (i.e. $\eta$ large) and perish in regions with little concentration of nutrition (i.e. $\eta$ small).

\begin{figure}[h]
\begin{center}
 \includegraphics[scale=.4]{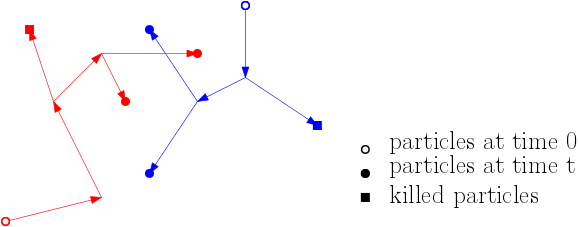}
\end{center}
\caption{An illustration of a situation starting from 2 particles,
after branching and dying, end up with 4 particles at time $t$.}
\end{figure}

In what follows we will assume that the dimension $d\le 3$. 
%
%


\noindent\underline{\it Scaling}.
Intuitively, $v$ should have high peaks in regions where the environment is most favorable for the particles ($\eta$ large), and it should have deep valleys elsewhere.
After long time, the bulk of the mass will be concentrated in different scales in these isolated islands, a phenomena known as intermittency\footnote{We refer to the surveys~\cite{MR1185878} and~\cite{MR3526112} for  the intermittency properties of the discrete PAM.}.
Being intermittent, we cannot expect to see nontrivial behavior on large spatial scales, since peaks and valleys are microscopic in their nature. 

On the other hand, the heat equation (when we set $\eta\equiv 0$) smooths particle density. This suggests that in order to obtain a meaningful SPDE scaling limit of the discrete PAM, we  ought to tune the strength of the potential $\eta$ so that its influence  (and the associated intermittency) balances with the smoothing effect of the heat equation. 
This should be thought as  analogous to the necessary tuning of the magnetic interaction length scale $\gamma^{-1}$ in Section~\ref{sec:phy-Phi4},
or,
the ``weakly asymmetric scaling'' that arises in KPZ equation convergence in Section \ref{sec:universal} below. This kind of tunings of the physical models are generally necessary
in order to obtain limiting SPDEs that are not scaling invariant,
 see  Remark~\ref{rmk:scaling}.

Indeed of \eqref{e:dPAM}, let us consider the discrete PAM with weakened noise: 
\begin{equ}\label{e:weak-potential}
   \partial_t v (t, x) = \Delta v (t, x) + \e^{2- \frac{d}{2}} v (t, x) \eta (x),\qquad (t,x)\in\R_+\times \Z^d.
\end{equ}
It is reasonable (though false) to believe that the diffusively rescaled solution $v(\varepsilon^{-2}t, \varepsilon^{-1}x)$ of~\eqref{e:weak-potential} should converge to the solution $w$  of  the continuous PAM:
\begin{equ}  \label{e:pam-no-renorm}
   \partial_t w (t, x) = \Delta w (t, x) + \sigma w (t, x) \xi (x), \qquad
    (t,x) \in \R_+ \times \R^d,
\end{equ}
where $\sigma^2 $ is the variance of $\eta(0)$ and $\xi$ is a spatial white noise (i.e., the Gaussian random field with mean zero and covariance $\E[\xi(x) \xi(y)] = \delta(x-y)$ described in Section \ref{sec:linear}).
Note that the choice of weak noise strength $\e^{2-\frac{d}{2}}$,  $\e^{2}$ reflects the  scaling dimension of the Laplacian operator, and $ \e^{-\frac{d}{2}}$ reflects the scaling dimension of the noise $\xi$ (since the scaling dimension of $\delta$ is $-d$).

\noindent\underline{\it Renormalization}.
It turns out that the above naive derivation of~\eqref{e:pam-no-renorm} from \eqref{e:weak-potential} is not correct (nor the result true) due to the singular nature of Eq. \eqref{e:pam-no-renorm}. In fact, even to make sense of  Eq. \eqref{e:pam-no-renorm} in the continuum setting, one has to introduce some renormalization.
In dimensions $d\ge 2$,  the total number of particles grows exponentially (even with the weak noise) and thus we have to renormalize by this expected growth rate in order to see a non-trivial behavior. More precisely, for $t > 0$ the expected number of particles at time $\e^{-2} t$ will be of order $e^{t C_\e}$ with a specific constant $C_\e =O(\log \e)>0$ in $d=2$ and $C_\e =O(\e^{-1})>0$ in $d=3$. So, one should instead consider 
\[
 u_\e(t,x) := e^{- t C_\e} v(\e^{-2}t, \e^{-1}x)
\]
which solves the modified and scaled discrete PAM:
\begin{equation}\label{eq:scaled lattice pam renorm}
   \partial_t u_\e(t,x) = \Delta^\e u_\e(t,x) + u_\e(t,x) (\eta_\e(x) - C_\e).
\end{equation}
Here $\Delta^\e$ is a scaling of the discrete Lapacian which  approximates the continuous Laplacian as $\e\to 0$, and $\eta_\e$ is a scaled version of $\eta$ that converges to $\xi$.
This is precisely the form of the renormalized
parabolic Anderson equation given by
 regularity structures in~\cite{Hairer14} and paracontrolled distributions in~ \cite{MR3406823}  if $\eta_\e$ is a mollification of the white noise.

It was rigorously proved in \cite{MR3736653} that, when $d=2$,  the solution $u_\e$ to~\eqref{eq:scaled lattice pam renorm} converges in law to the solution of the renormalized~\eqref{e:pam-no-renorm}, where
the potential $\eta$  is assumed to be a generally distributed (under certain very weak assumptions), and the discrete Laplacian can be  generator of any symmetric random walk whose increments have sufficiently many moments. In $d=3$ the same result is expected to hold under such general assumptions on the random walks and random environment. The proof of  \cite{MR3736653} is based on paracontrolled distributions as we introduced in Section~\ref{sec:RS}. The result   \cite{MR3736653} is further generalized by \cite{martin2017paracontrolled} which proves such a convergence result for the 
{\it nonlinear} parabolic Anderson model \eqref{e:gPAM} where the factor $f(u)$
models some interaction between the individual particles.
The proofs require showing convergence of the perturbative objects  discussed in Section~\ref{sec:RS}
in the discrete settings, which is one of the main technical challenge -- the argument of  \eqref{e:gPAM} relies on some general tools developed by \cite{MR3584558}.

\section{Singular SPDEs as universal objects}
\label{sec:universal}

The singular SPDEs we have studied here are universal objects which arise under various different approximation schemes. It can be shown that continuous mollification in position space as in \cite{Hairer14},  regularization in Fourier space as in \cite{MR3406823}, and lattice approximations (for instance \cite{GubPerk,CanMat} for KPZ, \cite{HaiMat,MR3758734} for $\Phi^4$, and references therein)
will all lead  to {\it the same} limiting  solution for the SPDEs  we discussed in the previous sections.  The choices of renormalization constants generally depend on the specific way of approximation in order to obtain the same limit.

Being universal objects means more:  each of these singular SPDEs governs the large scale fluctuation of a large class of physical models that have apparently very different microscopic interactions and details. In this section we demonstrate this universality for the KPZ equation by reviewing several recent scaling limit results. We choose to focus on the KPZ equation because, on one hand, there has been quite a lot of progress on KPZ in the last decade, and, on the other hand, several different approaches to solution theories of the KPZ equation have been found. These approaches all yield equivalent notions of solutions to this equation. However, when proving that a convergence results, some notions are better adapted to certain circumstances than others. Our examples will illustrate the application of some of these solution theories.

We start by discussing a continuous formulation of such a ``universality'' result proved for the KPZ equation, and then move to more discrete models. 

\subsection{KPZ equation with general nonlinearity}
\label{sec:CLT}

Consider the KPZ equation with quadratic nonlinear strength $\lambda$:
\begin{equ} [e:KPZ-CLTsec]
\partial_t H =\partial^2_x H + \lambda(\partial_x H)^2 +\xi
\end{equ}
and Gaussian space-time white noise $\xi$. This is a widely accepted model for a growing interface subject to three types of local forces: the term $\partial^2_x H $ models a  smoothing mechanism; the term $\lambda(\partial_x H)^2 $ models lateral growth (the interface tends to grow in the direction normal to the local slope), and $\xi$  models space-time randomness arising on a microscopic scale. 
In the seminal work of Kardar, Parisi and Zhang \cite{KPZ86}, they justified their eponymous equation by saying that ``{\it  the noise has a Gaussian distribution ... although the actual form of the distribution is  irrelevant}''. They continued to argue that (see the following figure) ``{\it growth occurs in a direction locally normal to the interface ... the increment projected along the $h$ axis is $\delta h = [(v\delta t)^2 + (v \delta t \partial_x h)^2]^{\frac12} \approx [v + \frac{v}{2} (\partial_x h)^2+\cdots]\delta t$ ... the original equation is regained. Such a nonlinear term is clearly expected \underline{in all situations} where lateral growth is allowed}''.

\begin{figure}[h]
  \begin{center}
    \includegraphics[width=0.4\textwidth]{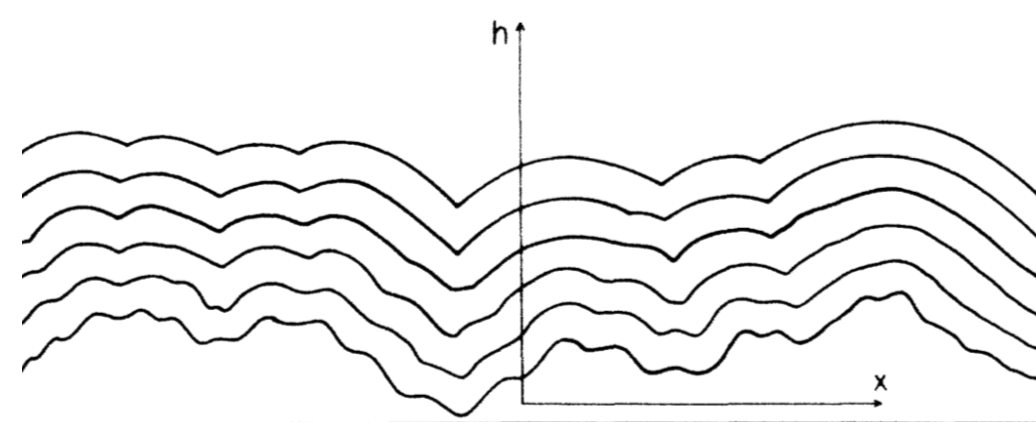}
    \qquad
     \includegraphics[width=0.15\textwidth]{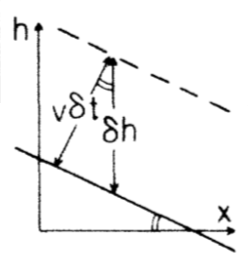}
  \end{center}
\caption{The original pictures in the seminal paper \cite{KPZ86} that explain why the growth term is chosen as $\lambda(\partial_x H)^2$.}
 \end{figure}

This heuristic derivation of the KPZ equation represents a claim of universality that can be formulated and proved mathematically. 
Consider the following class of continuous interface growth models where the {\it microscopic} growth equation is given by
\begin{equ} \label{e:JMmodel}
\partial_t h =\partial^2_x h +\sqrt\eps F(\partial_x h) +\eta
\end{equ}
where $F$ is a continuous function modeling the microscopic lateral growth (that could be rather complicated), and $\eta$ is a  continuous random field modeling the microscopic randomness (that is generally distributed,  not necessarily Gaussian).
The coefficient $\sqrt\eps$ here corresponds to the ``weakly asymmetric regime" meaning that the small scale interactions are tuned to be small (otherwise they would blow up upon scaling).
The challenge is to show that by scaling suitably and  applying the correct renormalization, one always obtains the standard KPZ equation \eqref{e:KPZ-CLTsec}.

%
%

Much progress has been achieved thanks to the recent developments of singular SPDE solution theory.
The following two results are achieved by the theory of regularity structures.

Hairer and Quastel \cite{KPZJeremy}  considered the above model \eqref{e:JMmodel} assuming that $\eta$ is Gaussian, but $F$ is an arbitrary even polynomial. They proved that for the rescaled height function defined by 
$\tilde h_\eps(x,t)= \eps^{\frac12}h(\eps^{-1} x,\eps^{-2}t)$ there exist constants $v_\eps$ and $\lambda$ such that $\tilde h_\eps(x, t)-v_\eps t$ converges as $\eps\to 0$ to the solution to the KPZ equation \eqref{e:KPZ-CLTsec}  with nonlinearity $\lambda (\partial_x H)^2$.
 This might sound not very surprising  because  under the above scaling, $\tilde h_\eps$ will satisfy a KPZ equation with ``error terms'' of the form $ \e^{k}(\partial_x \tilde h_\e)^{2k+2}$ with $k\ge 1$. However, a very  nontrivial fact is that the mean interface growth velocity $v_\eps$ (which is a ``renormalization constant" in this context) and the limiting interaction strength $\lambda$ depend explicitly on {\it all coefficients} of the polynomial $F$. The essence behind this nontrivial fact is again renormalization: via a similar  discussion as in Section~\ref{sec:PDE-renorm}, an ``error term'' such as $\e (\partial_x \tilde h_\e)^4$ needs to be renormalized as
 \begin{equ}[e:JM-ren-term]
 \e [(\partial_x \tilde h_\e)^4 -  C_\e (\partial_x \tilde h_\e)^2 -\bar C_\e]
 \end{equ}
  in order to converge; it turns out that $C_\e \sim O(\e^{-1})$, so that $\e C_\e $ is a finite constant contributing to the limiting coefficient $\lambda$. The constant $\e \bar C_\e$ is still divergent and makes a nontrivial contribution to the velocity $v_\eps$.
 Note that one  {\it does not} change the model \eqref{e:JMmodel} in this renormalization procedure; in fact one just re-shuffles the polynomial $F$ into linear combinations of renormalized terms like \eqref{e:JM-ren-term}.

In parallel,   \cite{HaiShen} focused on generally distributed $\eta$ in \eqref{e:JMmodel} with $F$ assumed to be quadratic $F(\partial_x h)=\lambda (\partial_x h)^2$. Under a very weak ``mixing'' assumption on the non-Gaussian field $\eta$ (so that $\eta_\eps (t,x)=\eps^{-\frac32} \eta(\frac{t}{\eps^2},\frac{x}{\eps})$  converges to the Gaussian white noise $\xi$ by classical central limit theorem), they demonstrated that with the correct choices of velocities $v^x$ and $v^y_{\eps}$,
\begin{equ}  \label{e:ref-frame}
\tilde h_\eps(x-v^x t, t)-v^y_{\eps} t
\end{equ}
 converges in law to the same solution  to KPZ  \eqref{e:KPZ-CLTsec} driven by Gaussian white noise $\xi$. It is important to note that  the convergence of noise $\eta_\e \to \xi$ does not imply convergence of solutions. In fact, this is the whole point of building the perturbative random objects in Section~\ref{sec:PDE-renorm} and  \ref{sec:RS}.
 Interestingly, the velocity  $v^x$ 
 in the Galilean transformation 
 shows up since  the covariance of $\eta$ 
is  generally not symmetric under spatial reflection $x\mapsto -x$;
 and the mean velocity  $v^y_{(\eps)} \sim \eps^{-1}C_0 + \eps^{-\frac12}C_1 + C_2$ where $C_1$ and $C_2$ depend on {\it the third and fourth cumulants\footnote{Gaussian random variables have zero cumulants of order three or higher. Therefore the third and fourth cumulants of $\eta$ represent some of the ``non-Gaussian bits'' of $\eta$.}} of $\eta$ respectively.
The velocities $v^x$ and $v^y_{\eps}$ are renormalization constants which,
as in our discussion in Section~\ref{sec:Well-posedness}, have explicit expressions in terms of the moments or cumulants of $\eta$ and heat kernels.


These two results \cite{KPZJeremy} and \cite{HaiShen} to a large extent justify the heuristic derivation performed  in \cite{KPZ86} of the KPZ equation.
In addition, they showed that the microscopic details -- such as higher order polynomial interactions and higher order cumulants in the microscopic randomness --
can contribute to the limiting coefficients or the reference frame in which the limit is observed, even when these details seem to just vanish by scaling.

Proving universality theorems in this continuum setting remains an active direction. More recently, \cite{HairerXu2018} generalized the result of \cite{KPZJeremy} by allowing more general function $F$ (not necessarily polynomial); they essentially only need to assume certain decay of the local distributional norm of the Fourier transform of $F$. In particular, it includes the case $F(u)= \sqrt{1+u^2}$ as originally considered in \cite{KPZ86}.  Under a strong assumption that  the system is put into the equilibrium state, \cite{Gubinelli2016hairer} provided a very simple proof using the notion of energy solutions (which is a weak solution as mentioned in Section~\ref{sec:weak-sol}).

\begin{rmk}
One can ask the same universality questions for the so called ``phase coexistence models'' (or ``reaction-diffusion models'') in three spatial dimensions which have the dynamical $\Phi^4$ equation as universal limit.
Essentially this is a continuous version of the model discussed in \eqref{e:Ising-Phi4},
just like the problem \eqref{e:JMmodel} can be viewed as continuous version of the models 
that we will discuss in the following subsections.
We refer to \cite{MR3772400} where polynomial nonlinearities are treated,
and   
\cite{furlan2017weak} where assumptions on  nonlinearities  are generalized to $C^9$ functions. In the case of polynomial nonlinearities the noise in the phase coexistence models can be non-Gaussian \cite{ShenXu}, but the case of general function nonlinearity with non-Gaussian noise is still open. Note that the dynamical $\Phi^4$ equation in two spatial dimensions does not have this universal property, because in two spatial dimensions
the equation with any polynomial nonlinearity is subcritical in the sense of Remark~\ref{rmk:subcritical} while in three spatial dimensions
only the cubic nonlinearity is subcritical.
\end{rmk}

\subsection{The solid on solid growth model}
\label{sec:ASEP}

In the previous subsection we saw how the KPZ equation arises as a universal limit of a large class of continuous interface growth models.
Another direction which has received a lot of interest is in studying how the KPZ equation arises as a scaling limit of discrete growth models. A particularly simple and well-studied one-dimensional interface growth model is the ``solid on solid interface growth model'', which is also called the ``corner growth model'' or ``asymmetric simple exclusion process'' height process. The   ``interface'' here is an integer-valued function $h(t,x)$ which depends on continuous time variable $t\in \R_+$ and discrete spatial variable $x\in \Z$, subject to the restriction that $|h(x+1)-h(x)|=1$ for all $x$.
\begin{figure}[h]
\begin{center}
 \includegraphics[scale=0.6]{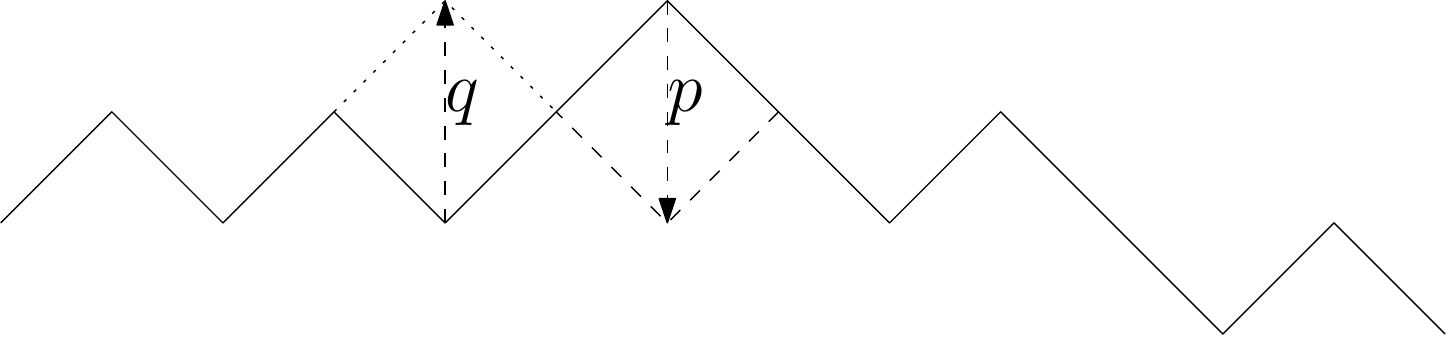}
\end{center}
\caption{A typical plot of a portion of $h(t,\cdot)$ at a given time $t$.}
\end{figure}
%
%
%
%
%

The height function $h$ evolves according to a simple Markov process.
As shown in the figure, a $\wedge$ shaped corner is flipped down by $2$ into a $\vee$ shaped corner at rate $p$;
and a $\vee$ shaped corner is randomly flipped up by $2$  into a $\wedge$ shaped corner at rate $q$. These flips occur independently and according to exponential waiting times (just as described in Section \ref{sec:phy-Phi4} for the Kac-Ising model Glauber dynamics).
Analogously to the discussion in Remark~\ref{rem:Markov process}, we may derive a microscopic evolution equation for $h$:
\begin{equ}[e:discreteKPZ1]
d h (t,x)
= 2 \big(
q\cdot \mathbf 1_{\Delta h (t,x)>0}
-
p\cdot  \mathbf 1_{\Delta h (t,x)<0}
\big) \,dt + dm(t,x) \;,
\end{equ}
where $m$ is an explicit martingale.
Denoting  $\nabla f(x):=f(x+1)-f(x)$ and using our assumption $|\nabla h(x)|=1$,
\eqref{e:discreteKPZ1} can be written into the form of a discrete KPZ type equation:
\begin{equ}[e:discreteKPZ2]
d h (t,x)
=\frac{p+q}{2}\Delta h \,dt+ \frac{q-p}{2} (1- \nabla h(t,x-1)\nabla h(t,x))\,dt+ dm(t,x)\;.
\end{equ}
In principle, one could implement the theories discussed in the previous sections such as regularity structures to prove that for a suitable choice of $v_\eps$,  $\varepsilon^{\frac{1}{2}}h(\varepsilon^{-1}x,\varepsilon^{-2}t) - v_\eps t$
converges to the solution of the KPZ equation. The main technical challenge in proving this convergence will be showing convergence of the perturbative objects (as described in the case of $\Phi^4$ equation in Section~\ref{sec:RS}),
which are objects built from the discrete heat kernel and the martingale noise $dm$. The martingale noise $dm$ is difficult to deal with since it  depends nontrivially on $h$, and there are a fairly large number of such objects that need to be handled.
Works such as \cite{MR3628883,MR3820327,Matetski2018martingale} have made progress in studying convergence of approximate SPDEs driven by martingales, but a complete treatment for the KPZ is still under progress.

Here we present a short-cut approach due to \cite{Bertini1997}
which applies well to the KPZ equation and certain approximation thereof. It is not, however applicable to general SPDEs.
Recall from \eqref{e:SHE} that if $u$ solves the stochastic heat equation (SHE)
 \begin{equ}[e:SHE-walsh]
\partial_t u =\frac12 \Delta u  + u\xi
  \end{equ}
then $H=\log u$ solves KPZ. This defines the ``Hopf-Cole solution'' and agrees with the other solution theories defined much later.
The Hopf-Cole transform makes rigorous sense as pointed out by  \cite{Bertini1995} since the SHE  is well-defined in a classical {\it It\^o sense} \cite{MR876085} and since \cite{Mueller91} proved {\it strict positivity} of $u$ (given rather general initial data). 
If we can find a version of Hopf-Cole transform at microscopic level, i.e., for the solid on solid model, and derive an approximate SHE for the exponentiated height function, we can just work at the level of the SHE. 

\cite{GJ88,Bertini1997} introduced\footnote{With slight tweak of notation adapted to our convention} the G\"artner  (or microscopic Hopf-Cole) transform:
 \begin{equ}[e:Gartner]
 Z(t,x):=(p/q)^{\frac12 h(t,x)} e^{t(1-2\sqrt{pq})} \;.
 \end{equ}
As in Section~\ref{sec:Renormalize-physical} the evolution of $Z$
can be decomposed into a drift part and a martingale.
The drift part turns out to precisely match a {\it discrete Laplacian},
namely
 \begin{equ}[e:dSHE]
 dZ (t,x) = \sqrt{pq} \Delta Z(t,x) dt+ dM(t,x)\;,
 \end{equ}
with $M$ a martingale. This is rather surprising since the change $h(t,x)\mapsto h(t,x)\pm 2$ during time $dt$ would seem to effect $Z$ in a nonlinear way. 
The fact that hold crucially relies on the fact that $\nabla h $ only takes two possible values $\{\pm 1\}$. In fact, for more general discrete growth processes, this miracle\footnote{There is, in fact, a broader class of discrete models which enjoy a version of this exact microscopic transform. These are models which enjoy  (at least one particle) Markov duality with respect to exponential functions of the height. See \cite{Corwin2016} for a general class of systems which enjoy such a relation, and \cite{Amir11,CT15,CSL16,Labbe16b,corwin2016open,G17,Shalin,Promit,CTinhomo} for examples of applications of this Markov duality relation to KPZ equation convergence.} is generally lost and the approach through the SHE is stymied (see, however, the work of \cite{Dembo2016} discussed further below).

\cite{Bertini1997} proved that with
$p/q=e^{-\sqrt\eps}$ (i.e., weakly asymmetry),
$Z(\eps^{-2}t,\eps^{-1}x)$ converges to the solution to the SHE.
Under this scaling
$Z(\eps^{-2}t,\eps^{-1}x) \approx e^{-\frac12 \sqrt{\eps}h(\eps^{-2}t,\eps^{-1}x) -O(\eps^{-1})t}$
so that the exponent proportional to $t$ in \eqref{e:Gartner}
is a ``vertical shift'' being subtracted as a renormalization.

 The proof in \cite{Bertini1997} relies on the notion of weak solution (or martingale problems) discussed in
Section~\ref{sec:weak-sol} (see, in particular, the discussion around Eq.~\eqref{e:N-and-Q}).
In the present context, one only needs to prove the discrete analogues of the processes $N_t$ and $Q_t$ in Eq.~\eqref{e:N-and-Q} are approximately martingales.
Indeed the fact that the discrete analogue of the process $N_t$ is martingale can be immediately read off from \eqref{e:dSHE}.
Showing the discrete analogue of  $Q_t$ is approximately martingale amounts to arguing that the quadratic variation of $dM$ in Section~\ref{e:dSHE} is approximately the quadratic variation of the limiting term $u\xi$, that is, $u^2$.
An explicit calculation shows that the quadratic variation involves a term which is approximately $u^2$ plus an ``error term'' which does not vanish point-wise when passing to the limit. The most challenging part in the proof is to show that this  ``error term''  vanishes in a suitable {\it averaged} sense; we refer to  \cite[Eq.~(3.16) and Sect.~4.2 ``Key estimate'']{Bertini1997} or \cite[Appendix~A]{Promit} for further details.

\begin{rmk}\label{rmk:ASEP}
The solid on solid interface growth model is equivalent to the ``asymmetric simple exclusion process'' which is a paradigmatic model for an ``interacting particle system''.
Particles occupy sites indexed by $\Z$ with the restriction of at most one particle per site (indicated by an occupation function at $x\in \Z$ of $\eta(x)=1$ for a particle or $\eta(x)=-1$ for a hole). 
Each particle attempts to jump left or right according to exponential rates $q$ and $p$, though only takes the jump if the destination site is not occupied at that time.
The occupation process $\eta(t,x)$ coincides with the discrete derivative of the solid so that  $\eta(t,x) = \nabla h(t,x)$ (where $\nabla$ acts in the $x$ variable).
Due to this, the result of \cite{Bertini1997} also shows that the fluctuation of $\eta$  converges (in a suitable sense) to the solution to the stochastic Burgers equation \eqref{eq:SBE} (which arises as the spatial derivative of the KPZ equation).
\end{rmk}


The work \cite{Dembo2016} generalized the result of \cite{Bertini1997}
by allowing growth or recession of {\it a section} (of length at most $4$) of the interface, see the left figure below. In terms of an interacting particle system, this corresponds with allowing jumps left and right further than distance 1. 
The challenge  is that the exact match with discrete Laplacian as in \eqref{e:dSHE} is no longer available and instead ``hydrodynamic limit'' techniques are used to control the defect in this matching.

\begin{figure}[h]
\includegraphics[scale=0.4]{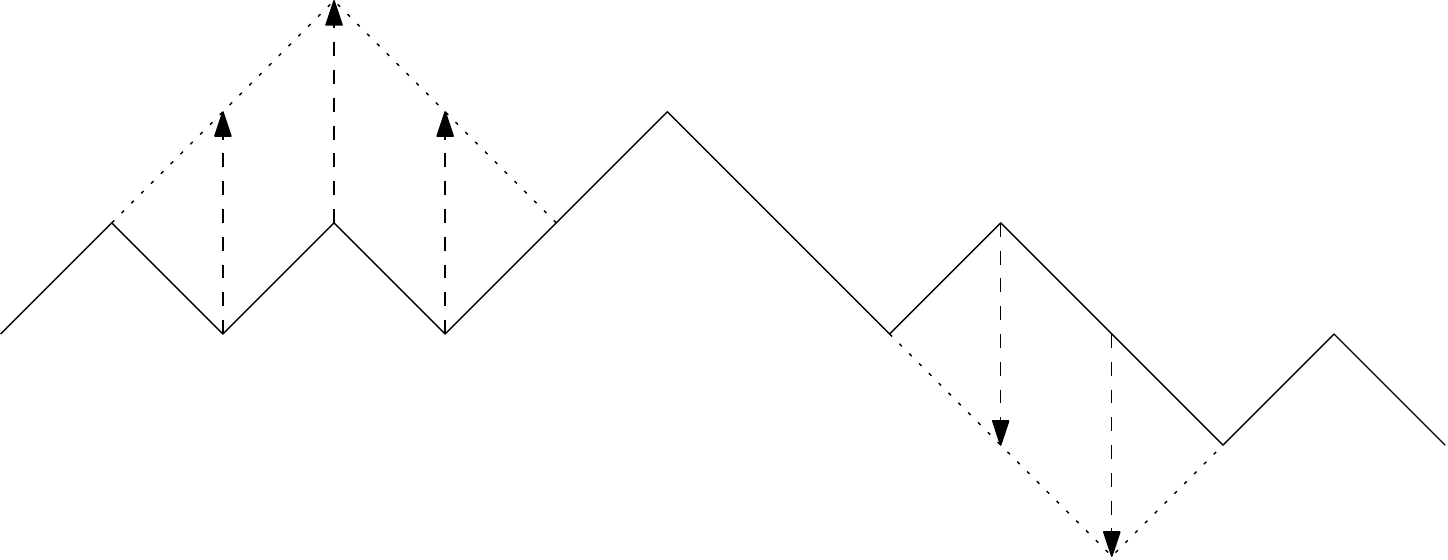}
$\qquad$
\includegraphics[scale=0.7]{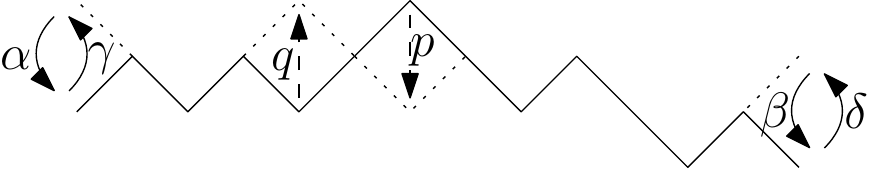}
\end{figure}

Another direction in which \cite{Bertini1997} has been generalized in  \cite{corwin2016open,Shalin} is to growth models on finite intervals with ``open boundaries'' (see the right figure).
This means that within the interval, the growth rule is as usual, but at the two ends, the height function randomly flips up or down with rates $\alpha,\beta,\gamma,\delta$.
With a fine tuning of these parameters\footnote{As in Remark~\ref{rmk:ASEP} there is an underlying interacting particle system called the ``OpenASEP'' which is equivalent to the discrete derivative of this growth model. The OpenASEP exhibits three phases: high density, low density and maximal current depending on the choices of rate parameters  $\alpha,\beta,\gamma,\delta$. The limiting SPDE arises when the  parameters  are tuned to approach a ``triple critical point'' separating these three phases.}
 it is  proved that $\varepsilon^{\frac{1}{2}}h(\varepsilon^{-1}x,\varepsilon^{-2}t) - C_\eps t$ (where  $C_\eps t$ is the same renormalization as the infinite interface case) converges to $H$ which is the solution of the KPZ equation on the spatial interval $[0,1]$, with inhomogeneous Neumann boundary condition $\partial_x H(t,x) |_{x=0} = A$  and $\partial_x H(t,x) |_{x=1} = B$. These two boundary conditions are only formal because $\partial_x H$ is distribution and some care is needed to properly define this solution\footnote{\cite{corwin2016open} defines it via the Hopf-Cole transform of the stochastic heat equation (SHE) with Robin boundary condition which can be defined in the style discussed in Section~\ref{sec:weak-sol}. \cite{GerencserHairer} studied this type of equation with boundary condition from the perspective or regularity structures, and \cite{GonPerMar} via energy solutions.}.

%
%

\subsection{Six vertex model}
\label{sec:6V}


The six vertex model (6V), originally introduced by \cite{LPauling35}, is generally  defined on a finite box in the 2D square lattice. Each lattice site is occupied by a vertex of one of the six types, with the restriction that vertices join up in a coherent manner (as shown in the following picture) and respect given boundary conditions.
Each vertex type has a weight parameter (so in general the model has six parameters), and the probability of a configuration is proportional to the product of all the vertex weights. The height function associated with a configuration is the $\Z$-valued function, denoted by $N(x,y)$, that changes by $1$ when crossing a line, as shown by the numbers in the figure.

  \begin{center}
    \includegraphics[width=0.28\textwidth]{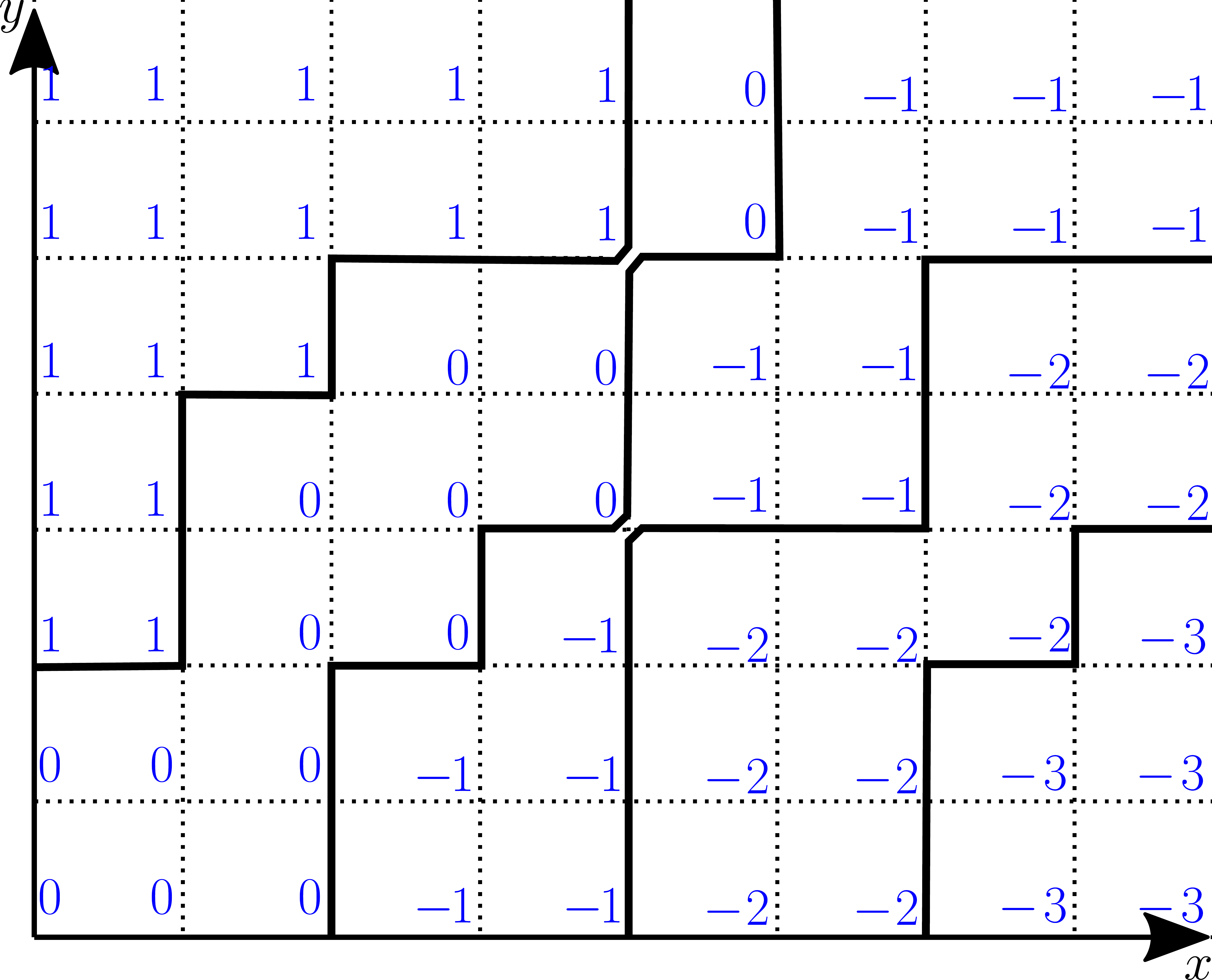}
  \end{center}

Two special situations are of particular interests: the {\it stochastic 6V} model and the {\it symmetric 6V} model. Their vertex weights are given by:
 \[
       \begin{tabular}{|c|c|c|c|c|c|c|}
           \hline
\footnotesize{Vertex type}
&
 \begin{tikzpicture}[scale=0.3]
 \draw[thick,white] (0,-1.4) -- (0,1.4);
 \draw[thick,white] (-1.4,0) -- (1.4,0);
            \draw[dotted] (-1,0) -- (0,0);
            \draw[dotted] (0,0) -- (1,0);
            \draw[dotted] (0,-1) -- (0,0);
            \draw[dotted] (0,0) -- (0,1);
            \end{tikzpicture}
&
            \begin{tikzpicture}[scale=0.3]
 \draw[thick,white] (0,-1.4) -- (0,1.4);
 \draw[thick,white] (-1.4,0) -- (1.4,0);
            \draw[ultra thick] (-1,0) -- (0,0);
            \draw[ultra thick] (0,0) -- (1,0);
            \draw[ultra thick] (0,-1) -- (0,0);
            \draw[ultra thick] (0,0) -- (0,1);
            \draw[thick][white] (-.1,-.1) -- (.1,.1);
            \end{tikzpicture}
&
            \begin{tikzpicture}[scale=0.3]
  \draw[thick,white] (0,-1.4) -- (0,1.4);
 \draw[thick,white] (-1.4,0) -- (1.4,0);
            \draw[dotted] (-1,0) -- (0,0);
            \draw[dotted] (0,0) -- (1,0);
            \draw[ultra thick] (0,-1) -- (0,0);
            \draw[ultra thick] (0,0) -- (0,1);
            \end{tikzpicture}
&
            \begin{tikzpicture}[scale=0.3]
 \draw[thick,white] (0,-1.4) -- (0,1.4);
 \draw[thick,white] (-1.4,0) -- (1.4,0);
            \draw[ultra thick] (-1,0) -- (0,0);
            \draw[ultra thick] (0,0) -- (1,0);
            \draw[dotted] (0,-1) -- (0,0);
            \draw[dotted] (0,0) -- (0,1);
            \end{tikzpicture}
&
            \begin{tikzpicture}[scale=0.3]
 \draw[thick,white] (0,-1.4) -- (0,1.4);
 \draw[thick,white] (-1.4,0) -- (1.4,0);
            \draw[dotted] (-1,0) -- (0,0);
            \draw[dotted] (0,0) -- (0,1);
            \draw[ultra thick] (0,-1) -- (0,0);
            \draw[ultra thick] (0,0) -- (1,0);
            \draw[thick][white] (-.1,-.1) -- (.1,.1);
            \end{tikzpicture}
&
            \begin{tikzpicture}[scale=0.3]
 \draw[thick,white] (0,-1.4) -- (0,1.4);
 \draw[thick,white] (-1.4,0) -- (1.4,0);
            \draw[ultra thick] (-1,0) -- (0,0);
            \draw[ultra thick] (0,0) -- (0,1);
            \draw[dotted] (0,-1) -- (0,0);
            \draw[dotted] (0,0) -- (1,0);
            \draw[thick][white] (-.1,-.1) -- (.1,.1);
            \end{tikzpicture}
  \\
           \hline
\footnotesize{Stochastic 6V weights} &        $1$ & $1$ & $b_1$ & $b_2$ & $1-b_1$ & $1-b_2$ \\
            \hline
\footnotesize{Symmetric 6V weights} &
        $a$ & $a$ & $b$ & $b$ & $c$ & $c$ \\
            \hline
        \end{tabular}
    \]
The {\it stochastic 6V}, proposed by Gwa and Spohn \cite{GS92} in 1992, depends on two parameters $b_1,b_2 \in (0,1)$, as shown in the table.
This choice is special since if we treat the bottom and left edges coming into vertices as ``inputs'' and the top and right edges as ``outputs'', then the sum of weights over all outputs always equals 1. This enables us to define the model on the entire first quadrant through a Markovian update. Essentially, the lines can be thought of as paths taken by particles, and the vertex weights determine the probabilities associated to the local moves made by each particle. It is, in fact, also possible to define this particle system on the full upper-half-plane, in which case it is natural to think of the $x$-axis as space and the $y$-axis as time (so we use $t$ instead of $y$ to label this direction). 

In \cite{Promit} it was proved that under a suitable ``weak asymmetry'', where $b_2\to b_1$, upon scaling and renormalizing, the fluctuation of the height function $N(t,x)$ converges to the solution of the KPZ equation. More precisely,
fixing any density $\rho \in(0,1)$ of lines entering from the horizontal axis, $b_1\in(0,1)$, and $b_2=b_1 e^{-\sqrt\e}$,
\begin{equ}
	\label{eq:S6VtoKPZ}
	\sqrt\e \Big( N \big(\e^{-2}t,\e^{-1}x+\mu \e^{-2} t\big)
	 - \rho(\e^{-1}x+\mu \e^{-2} t) \Big) - \e^{-2}t\log\lambda
	\Longrightarrow
	H(t,x)
\end{equ}
in the topology $C(\R_+, C(\R))$, for some constants $\mu$ and $\lambda$, where $H$ is the solution to the KPZ equation with more general coefficients
\begin{align}
	\label{eq:KPZmkD}
	\partial_t  H
	 = \tfrac{\nu}{2} \partial_x^2  H
	-   \tfrac{\kappa}{2} \big( \partial_x  H\big)^2
	+ \sqrt{D} \xi\;.
\end{align}
Here the equation coefficients $\nu=\frac{2b_1}{1-b_1} $, $\kappa=\frac{2b_1}{1-b_1}$ and $D = \frac{2b_1\rho(1-\rho)}{1-b_1}$. The $\lambda$ and $\mu$ required in the convergence \eqref{eq:S6VtoKPZ} also depend on $b_1,b_2,\rho$ via explicit formulae, which behave like $\lambda=	1 - \rho \sqrt\e + \mathcal{O}(\e)$ and $\mu = 	1 + \tfrac{b_1-2b_1\rho}{b_1-1} \sqrt{\e} + \mathcal{O}(\e)$.

 \eqref{eq:S6VtoKPZ} shows that in order to observe the KPZ equation fluctuation, one needs to follow a characteristic direction $\mu$, then tilt the height function by $\e^{-1}\rho x$ (since by definition the function $N$ was tilted), and finally center the height function by  subtracting an overall growth speed multiplied by $t$ (i.e. the terms proportional to $t$ in  \eqref{eq:S6VtoKPZ}).

The proof of this result relies upon the Hopf-Cole solution, as in Section~\ref{sec:ASEP}, namely, for the KPZ equation as in \eqref{eq:KPZmkD}, $H =-\frac{\nu}{\kappa} \log u$ where $u$ solves the stochastic heat equation with multiplicative noise as introduced in \eqref{e:SHE}
\begin{equ}[e:6VSHE]
\partial_t u
	= \tfrac{\nu}{2} \partial_x^2 u
	+ \tfrac{\kappa\sqrt{D}}{\nu}\xi u \;.
\end{equ}
 Remarkably (and owing to the one-particle duality from \cite{Corwin2016}), just like for the solid-on-solid model in Section~\ref{sec:ASEP}, the stochastic 6V model admits an exact microscopic analog of this Hopf-Cole transform, transforming it into a discrete  version of \eqref{e:6VSHE}.
The proof  that this discrete version of \eqref{e:6VSHE} converges to the limiting continuous one requires control over the martingale quadratic variation. The method of \cite{Bertini1997} discussed in Section \ref{sec:ASEP} does not seem to apply here and \cite{Promit} had to introduce a new method relying upon the one and two-particle duality enjoyed by this model.

%
%

Turning very briefly to the older and very well-studied {\it symmetric 6V} model as originally proposed by \cite{LPauling35} in 1935, it turns out that the KPZ equation (or rather stochastic Burgers equation) also arises here. Among the (widely believed, though not generally proved or constructed) two-parameter family of infinite-volume, translation invariant ergodic Gibbs states, \cite{A16aa} (see also earlier observations of \cite{Bukman1995} in physics) constructed a one-parameter subfamily of ``stochastic Gibbs states'' using the stochastic 6V model. This construction only applies in the ``ferroelectric'' regime for parameters $a,b$ and $c$. Zooming into the ``ferroelectric-disorder'' phase transition point, and scaling along ``characteristic directions'' in the manner described above, these stochastic Gibbs states are shown in \cite{Promit} to converge to stationary solutions to the stochastic Burgers equation.

Let us mention that there are some other works studying SPDE limits of this or closely related models, see, for example, \cite{CT15,BO16} or \cite{BorGorforth, ST18} (the later involves a limiting linear hyperbolic SPDE termed the ``stochastic telegraph equation''.

%

\subsection{Interacting Brownian motions}
\label{sec:energy}

We recall one more convergence result for interacting Brownian motions, as well as the method of energy solutions used.

Consider a collection of independent standard Brownian motions $\{W_t(i):i\in \Z\}$, and define an ``interface profile'' $\phi_t:\Z\to \R$ via
\begin{equ}\label{e:InteractingBM}
   \mathrm{d}\phi_t (i) = \left\{ pV' (\phi_t (i + 1) - \phi_t (i)) - qV' (\phi_t (i) - \phi_t (i-1)) \right\} \mathrm{d}t + \mathrm{d}W_t (i),
\end{equ}
for $(t,i) \in \R_+ \times \Z$ and $p,q \ge 0$ with $p+q=1$.
Here, $V$ is a potential function. When $V$ is a constant potential, $\{\phi_t(i)\}_{i\in \Z}$
is simply a collection of independent Brownian motions.

$\phi_t$ may be thought of as a one-dimensional interface separating two phases and we are interested in the random dynamics of this interface. The case $p=q=1/2$ describes a type of balance between the two phases and the interface dynamics have no net growth; in this case the system is known as the (one-dimensional) \emph{Ginzburg-Landau $\nabla \phi$ interface model} which has been intensely studied, see  \cite{MR894407,MR954674,MR826856,MR1044937,MR1162798}.

\cite{MR3663617} proved  that under ``weak asymmetry''
\begin{equ}[e:BMweak-asym]
p=(1+\sqrt{\varepsilon})/2 \;,
\qquad
q=(1-\sqrt{\varepsilon})/2 \;,
\end{equ} 
the properly scaled fluctuation of $\phi$ converges to the solution to the KPZ equation. This result relies on the use of energy solutions as mentioned in Section~\ref{sec:weak-sol}, and has only been proved for stationary initial data. 
In fact, the ``slope function'' (or discrete derivative in $i$) 
\[
u_t(i) := \phi_t (i) - \phi_t (i - 1)
\]
admits
 a {\it one parameter family} of stationary measures parametrized by $\lambda \in \mathbb{R}$:
\begin{equ}[e:1-par-stat-mea]
   {\mu}_{\lambda} (d u) := \prod_{j = - \infty}^{\infty} p_\lambda\big(u(j)\big) d u(j)\;,
   \qquad
  \mbox{where}\;\; p_\lambda(u) :=Z_{\lambda}^{-1}e^{\lambda u - V (u)}
\end{equ}
and $Z_{\lambda} =\int_{\R} e^{\lambda u - V (u)}du$ is a normalization constant.
More precisely, \cite{MR3663617} proved the following result.
Assuming $u_0(i) = \phi_0 (i) - \phi_0 (i - 1)$ has the probability distribution $\mu_\lambda$ for a fixed $\lambda \in\mathbb{R}$,
and let $ \rho'(\lambda) := \int_\R u \, p_\lambda(u)  d u$ be
   the mean of the coordinates $u(j)$ under $\mu_\lambda$. Then, under weak asymmetry \eqref{e:BMweak-asym}, the scaled and renormalized profile process 
   \begin{equ}[e:scaled and renormalized phi]
    \e^{\frac12} (\phi (\e^{-2} t, \lfloor \e^{-1} x  - c_\e t \rfloor) ) - C_\e
    \end{equ}
      converge to the solution to
\[
  \partial_t  H = \frac{ 1}{2 \sigma_{\lambda}^{2}}
     \partial_x^2  H
      - \frac{m_{3,\lambda}}{2 \sigma_{\lambda}^6} (\partial_x  H)^2 + \xi,
  \]
  where $ c_\e := \e^{-1 / 2} \sigma_{\lambda}^{-2}$, the diverging renormalization constants $C_\e$ is explicit, and where
$ m_{k,\lambda}$  are the moments of the measure $p_\lambda$
\[
 m_{k,\lambda} := \int_\R (u - \rho'(\lambda))^k p_\lambda(u)d u \;,
 \qquad
 \sigma^2_\lambda := m_{2,\lambda} \;.
\]
Since the energy solution is defined for the stochastic Burgers equation, and the solution to  stochastic Burgers equals the spatial derivative of that of the KPZ equation  (Remark~\ref{rmk:ASEP}),
the above result is proved via showing that
$ \e^{-1 / 2} (u (\e^{-2} t, \lfloor \e^{-1} x  - c_\e t \rfloor) - \rho'(\lambda))$ converges to the energy solution to the stochastic Burgers equation.

The energy solution method for the KPZ / stochastic Burgers equation convergence was initiated in the work of the Jara and Gon\c calves \cite{GJ10} (cf. \cite{Assing}). Initially this approach only provided ``tightness'', in other words, existence of limit, and it was not known whether energy solutions were unique. Uniqueness (and hence the identification with the Hopf--Cole solution of KPZ as discussed in the previous sections) was proved in \cite{GP2015a}. This approach has been applied to prove that a wide variety of particle systems converge to the \ac{KPZ} equation, see \cite{GJ14a,GJS15,TGS16,GJ13,GJ16,GonPerMar}. These results all require that the system has an invariant measure (i.e., stationary, or equilibrium measure) -- like the measures  \eqref{e:1-par-stat-mea} for the model \eqref{e:InteractingBM} above --
and that the initial condition  is thusly distributed.
The method of proof relies heavily on having well-developed hydrodynamic theory estimates available. Quite recently,  however,  \cite{MR3828193} and \cite{Yang} have extended this method to include more general initial data as we mentioned in Section~\ref{sec:weak-sol}.

\begin{rmk}\label{rmk:recap}
Let us recapitulate this section with several remarks.

First, it would be interesting to compare the expressions 
\eqref{e:ref-frame} in the continuous interface model,
\eqref{eq:S6VtoKPZ} in the six vertex model, and 
\eqref{e:scaled and renormalized phi} for the interface driven by Brownian motions.
The common feature is that in all of these cases, the KPZ equation arises via suitable adjustments of  reference frames  dictated by the renormalizations in these discrete models.

Also, we have seen that the proof strategies of these convergence results depend on which notion of solution 
is being chosen in the context. The Hopf-Cole solution to the KPZ equation has the obvious advantage in turning the problem to a linear equation \eqref{e:SHE-walsh} and \eqref{e:6VSHE}, as we demonstrated in the solid on solid growth model and the six vertex model.
Although it has been shown for a number of models one can implement this Hopf-Cole transform at the microscopic level, this relation certainly does not exactly hold in general;
it is also a special feature for  KPZ  that is absent for other nonlinear SPDEs discussed in Section~\ref{sec:NSPDE}. The energy solution has been applied to proving KPZ / stochastic Burgers equation convergence results of systems which have an equilibrium invariant measure and start from this equilibrium.  The recent works \cite{MR3828193,Yang} extended the type of initial data and \cite{gubinelli2018infinitesimal} studied more general class of equations; further extending the scope of applicability of energy solution to physical models seems to be an interesting direction. The theory of regularity structures and paracontrolled distributions method provide robust solution theories for very large classes of equations, and yielded universality results in the continuum setting of Section~\ref{sec:CLT}. Proving such convergence results for discrete physical systems driven by martingale noises in a more {\it systematic} way would be another interesting direction; besides the results cited above let us mention \cite{HaiMat,erhard2017discretisation,Matetski2018martingale} for some initial progress on discrete regularity structures.
\end{rmk}

%
%
%
%
%
%
%
%
%
%

%
%
%
\bibliographystyle{alphaabbr}
\bibliography{Reference}

\newcommand{\etalchar}[1]{$^{#1}$}
\begin{thebibliography}{DKM{\etalchar{+}}09}

\bibitem[AC15]{AllezChouk}
R.~Allez and K.~Chouk.
\newblock {The continuous Anderson Hamiltonian in dimension two}.
\newblock {\em arXiv:1511.02718}, 2015.

\bibitem[ACQ11]{Amir11}
G.~Amir, I.~Corwin, and J.~Quastel.
\newblock Probability distribution of the free energy of the continuum directed
  random polymer in {$1+1$} dimensions.
\newblock {\em Commun. Pure Appl. Math.}, 64(4):466--537, 2011.

\bibitem[Agg16]{A16aa}
A.~Aggarwal.
\newblock Current fluctuations of the stationary {ASEP} and six-vertex model.
\newblock {\em Duke Math. J.}, 167:269--384, 2016.

\bibitem[AHR96]{MR1396758}
S.~Albeverio, Z.~Haba, and F.~Russo.
\newblock Trivial solutions for a non-linear two-space-dimensional wave
  equation perturbed by space-time white noise.
\newblock {\em Stoch. and Stoch. Rep.}, 56(1-2):127--160, 1996.

\bibitem[AK17]{albeverio2017invariant}
S.~Albeverio and S.~Kusuoka.
\newblock The invariant measure and the flow associated to the $\phi^4_3
  $-quantum field model.
\newblock {\em arXiv:1711.07108}, 2017.

\bibitem[AR91]{AlbRock91}
S.~Albeverio and M.~R{\"o}ckner.
\newblock Stochastic differential equations in infinite dimensions: solutions
  via {D}irichlet forms.
\newblock {\em Probab. Theory Related Fields}, 89(3):347--386, 1991.

\bibitem[Ass13]{Assing}
S.~Assing.
\newblock A rigorous equation for the {C}ole--{H}opf solution of the
  conservative {KPZ} equation.
\newblock {\em Stoch. PDEs: Analy. Comp.}, 1(2):365--388, 2013.

\bibitem[BB16]{MR3475460}
I.~Bailleul and F.~Bernicot.
\newblock Heat semigroup and singular {PDE}s.
\newblock {\em J. Funct. Anal.}, 270(9):3344--3452, 2016.
\newblock With an appendix by F. Bernicot and D. Frey.

\bibitem[BC95]{Bertini1995}
L.~Bertini and N.~Cancrini.
\newblock The stochastic heat equation: {F}eynman--{K}ac formula and
  intermittence.
\newblock {\em J. Statist. Phys.}, 78(5-6):1377--1401, 1995.

\bibitem[BCCH17]{bruned2017renormalising}
Y.~Bruned, A.~Chandra, I.~Chevyrev, and M.~Hairer.
\newblock Renormalising spdes in regularity structures.
\newblock {\em arXiv:1711.10239}, 2017.

\bibitem[BCD11]{MR2768550}
H.~Bahouri, J.-Y. Chemin, and R.~Danchin.
\newblock {\em Fourier analysis and nonlinear partial differential equations},
  volume 343 of {\em Grundlehren der Mathematischen Wissenschaften [Fundamental
  Principles of Mathematical Sciences]}.
\newblock Springer, Heidelberg, 2011.

\bibitem[BDH19]{MR3916262}
I.~Bailleul, A.~Debussche, and M.~Hofmanov\'{a}.
\newblock Quasilinear generalized parabolic {A}nderson model equation.
\newblock {\em Stoch. Partial Differ. Equ. Anal. Comput.}, 7(1):40--63, 2019.

\bibitem[Ber70]{berezinsky1970destruction}
V.~Berezinsky.
\newblock Destruction of long range order in one-dimensional and
  two-dimensional systems having a continuous symmetry group. {I}. classical
  systems.
\newblock {\em Zh. Eksp. Teor. Fiz.}, 32:493--500, 1970.

\bibitem[BG97]{Bertini1997}
L.~Bertini and G.~Giacomin.
\newblock Stochastic {B}urgers and {KPZ} equations from particle systems.
\newblock {\em Commun. Math. Phys.}, 183(3):571--607, 1997.

\bibitem[BG18]{BorGorforth}
A.~Borodin and V.~Gorin.
\newblock A stochastic telegraph equation from the six-vertex model.
\newblock {\em arXiv:1803.09137}, 2018.

\bibitem[BGHZ19]{hairer2019geo}
Y.~Bruned, F.~Gabriel, M.~Hairer, and L.~Zambotti.
\newblock Geometric stochastic heat equations.
\newblock {\em arXiv:1902.02884}, 2019.

\bibitem[BHZ16]{bruned2016algebraic}
Y.~Bruned, M.~Hairer, and L.~Zambotti.
\newblock Algebraic renormalisation of regularity structures.
\newblock {\em Inventiones mathematicae}, pages 1--118, 2016.

\bibitem[BO17]{BO16}
A.~Borodin and G.~Olshanski.
\newblock The {ASEP} and determinantal point processes.
\newblock {\em Commun. Math. Phys.}, 353:853--903, 2017.

\bibitem[Bou96]{MR1374420}
J.~Bourgain.
\newblock Invariant measures for the {$2$}{D}-defocusing nonlinear
  {S}chr\"{o}dinger equation.
\newblock {\em Comm. Math. Phys.}, 176(2):421--445, 1996.

\bibitem[BP15]{BP16b}
A.~Borodin and L.~Petrov.
\newblock Integrable probability: stochastic vertex models and symmetric
  functions.
\newblock {\em Stochastic Processes and Random Matrices: Lecture Notes of the
  Les Houches Summer School}, 104, 2015.

\bibitem[Bru68]{Brush}
S.~G. Brush.
\newblock {A History of Random Processes: I. Brownian Movement from Brown to
  Perrin}.
\newblock {\em Archive for History of Exact Sciences}, 5(1):1--36, 1968.

\bibitem[BS95]{Bukman1995}
D.~J. Bukman and J.~D. Shore.
\newblock The conical point in the ferroelectric six-vertex model.
\newblock {\em J. Stat. Phys.}, 78(5):1277--1309, 1995.

\bibitem[CC18a]{MR3785598}
G.~Cannizzaro and K.~Chouk.
\newblock Multidimensional {SDE}s with singular drift and universal
  construction of the polymer measure with white noise potential.
\newblock {\em Ann. Probab.}, 46(3):1710--1763, 2018.

\bibitem[CC18b]{MR3846835}
R.~Catellier and K.~Chouk.
\newblock Paracontrolled distributions and the 3-dimensional stochastic
  quantization equation.
\newblock {\em Ann. Probab.}, 46(5):2621--2679, 2018.

\bibitem[CD18]{ChatterjeeDunlapKPZ}
S.~Chatterjee and A.~Dunlap.
\newblock Constructing a solution of the $(2+ 1) $-dimensional {KPZ} equation.
\newblock {\em arXiv:1809.00803}, 2018.

\bibitem[CGP17]{MR3736653}
K.~Chouk, J.~Gairing, and N.~Perkowski.
\newblock An invariance principle for the two-dimensional parabolic {A}nderson
  model with small potential.
\newblock {\em Stoch. Partial Differ. Equ. Anal. Comput.}, 5(4):520--558, 2017.

\bibitem[CGST18]{Promit}
I.~Corwin, P.~Ghosal, H.~Shen, and L.-C. Tsai.
\newblock Stochastic {PDE} limit of the {Six Vertex Model}.
\newblock {\em arXiv:1803.08120}, 2018.

\bibitem[CH16]{chandra2016analytic}
A.~Chandra and M.~Hairer.
\newblock An analytic {BPHZ} theorem for regularity structures.
\newblock {\em arXiv:1612.08138}, 2016.

\bibitem[CHI15]{MR3296821}
D.~Chelkak, C.~Hongler, and K.~Izyurov.
\newblock Conformal invariance of spin correlations in the planar {I}sing
  model.
\newblock {\em Ann. of Math. (2)}, 181(3):1087--1138, 2015.

\bibitem[CHS18]{sineGordonwholeregime}
A.~Chandra, M.~Hairer, and H.~Shen.
\newblock The dynamical sine-{G}ordon model in the full subcritical regime.
\newblock {\em arXiv:1808.02594}, 2018.

\bibitem[CM94]{MR1185878}
R.~A. Carmona and S.~A. Molchanov.
\newblock Parabolic {A}nderson problem and intermittency.
\newblock {\em Mem. Amer. Math. Soc.}, 108(518):viii+125, 1994.

\bibitem[CM16]{CanMat}
G.~Cannizzaro and K.~Matetski.
\newblock Space-time discrete {KPZ} equation.
\newblock {\em Commun. Math. Phys.}, 358(2):521--588, 2016.

\bibitem[Con12]{conti2012solitonization}
C.~Conti.
\newblock Solitonization of the {A}nderson localization.
\newblock {\em Physical Review A}, 86(6):061801, 2012.

\bibitem[Cor14]{CorwinICM2014}
I.~Corwin.
\newblock {M}acdonald processes, quantum integrable systems and the
  {K}ardar-{P}arisi-{Z}hang universality class.
\newblock {\em Proceedings of the 2014 ICM, arXiv:1403.6877}, 2014.

\bibitem[CP16]{Corwin2016}
I.~Corwin and L.~Petrov.
\newblock Stochastic higher spin vertex models on the line.
\newblock {\em Commun. Math. Phys.}, 343(2):651--700, 2016.

\bibitem[CR99]{MR1661761}
R.~A. Carmona and B.~Rozovskii, editors.
\newblock {\em Stochastic partial differential equations: six perspectives},
  volume~64 of {\em Mathematical Surveys and Monographs}.
\newblock American Mathematical Society, Providence, RI, 1999.

\bibitem[CS18]{corwin2016open}
I.~Corwin and H.~Shen.
\newblock Open {ASEP} in the weakly asymmetric regime.
\newblock {\em Comm. Pure Appl. Math.}, 71(10):2065--2128, 2018.

\bibitem[CST18]{CSL16}
I.~Corwin, H.~Shen, and L.-C. Tsai.
\newblock {${\rm ASEP}(q,j)$} converges to the {KPZ} equation.
\newblock {\em Ann. Inst. Henri Poincar\'{e} Probab. Stat.}, 54(2):995--1012,
  2018.

\bibitem[CSZ17a]{MR3584558}
F.~Caravenna, R.~Sun, and N.~Zygouras.
\newblock Polynomial chaos and scaling limits of disordered systems.
\newblock {\em J. Eur. Math. Soc. (JEMS)}, 19(1):1--65, 2017.

\bibitem[CSZ17b]{MR3719953}
F.~Caravenna, R.~Sun, and N.~Zygouras.
\newblock Universality in marginally relevant disordered systems.
\newblock {\em Ann. Appl. Probab.}, 27(5):3050--3112, 2017.

\bibitem[CSZ18]{Caravenna2Dentire}
F.~Caravenna, R.~Sun, and N.~Zygouras.
\newblock The two-dimensional {KPZ} equation in the entire subcritical regime.
\newblock {\em arXiv:1812.03911}, 2018.

\bibitem[CT17]{CT15}
I.~Corwin and L.-C. Tsai.
\newblock {KPZ} equation limit of higher-spin exclusion processes.
\newblock {\em Ann. Probab.}, 45:1771--1798, 2017.

\bibitem[CT18]{CTinhomo}
I.~Corwin and L.-C. Tsai.
\newblock {SPDE Limit of Weakly Inhomogeneous ASEP}.
\newblock {\em arXiv:1806.09682}, 2018.

\bibitem[CW17]{ChandraWeber}
A.~Chandra and H.~Weber.
\newblock Stochastic {PDE}s, regularity structures, and interacting particle
  systems.
\newblock {\em Ann. Fac. Sci. Toulouse Math. (6)}, 26(4):847--909, 2017.

\bibitem[CY92]{MR1162798}
C.~C. Chang and H.-T. Yau.
\newblock Fluctuations of one-dimensional {G}inzburg-{L}andau models in
  nonequilibrium.
\newblock {\em Comm. Math. Phys.}, 145(2):209--234, 1992.

\bibitem[Daw72]{MR0321178}
D.~A. Dawson.
\newblock Stochastic evolution equations.
\newblock {\em Math. Biosci.}, 15:287--316, 1972.

\bibitem[DGP17]{MR3663617}
J.~Diehl, M.~Gubinelli, and N.~Perkowski.
\newblock The {K}ardar-{P}arisi-{Z}hang equation as scaling limit of weakly
  asymmetric interacting {B}rownian motions.
\newblock {\em Comm. Math. Phys.}, 354(2):549--589, 2017.

\bibitem[DK90]{MR1079726}
S.~K. Donaldson and P.~B. Kronheimer.
\newblock {\em The geometry of four-manifolds}.
\newblock Oxford Mathematical Monographs. The Clarendon Press, Oxford
  University Press, New York, 1990.
\newblock Oxford Science Publications.

\bibitem[DKM{\etalchar{+}}09]{MR1500166}
R.~Dalang, D.~Khoshnevisan, C.~Mueller, D.~Nualart, and Y.~Xiao.
\newblock {\em A minicourse on stochastic partial differential equations},
  volume 1962 of {\em Lecture Notes in Mathematics}.
\newblock Springer-Verlag, Berlin, 2009.
\newblock Held at the University of Utah, Salt Lake City, UT, May 8--19, 2006,
  Edited by Khoshnevisan and Firas Rassoul-Agha.

\bibitem[DM17]{DebusscheMartin}
A.~Debussche and J.~Martin.
\newblock Solution to the stochastic {S}chrodinger equation on the full space.
\newblock {\em arXiv:1707.06431}, 2017.

\bibitem[DPD02]{MR1941997}
G.~Da~Prato and A.~Debussche.
\newblock Two-dimensional {N}avier-{S}tokes equations driven by a space-time
  white noise.
\newblock {\em J. Funct. Anal.}, 196(1):180--210, 2002.

\bibitem[DPD03]{MR2016604}
G.~Da~Prato and A.~Debussche.
\newblock Strong solutions to the stochastic quantization equations.
\newblock {\em Ann. Probab.}, 31(4):1900--1916, 2003.

\bibitem[DPZ14]{DPZ}
G.~Da~Prato and J.~Zabczyk.
\newblock {\em Stochastic equations in infinite dimensions}, volume 152 of {\em
  Encyclopedia of Mathematics and its Applications}.
\newblock Cambridge University Press, Cambridge, second edition, 2014.

\bibitem[DT16]{Dembo2016}
A.~Dembo and L.-C. Tsai.
\newblock Weakly asymmetric non-simple exclusion process and the
  {K}ardar--{P}arisi--{Z}hang equation.
\newblock {\em Commun. Math. Phys.}, 341(1):219--261, 2016.

\bibitem[DW18]{MR3785398}
A.~Debussche and H.~Weber.
\newblock The {S}chr\"{o}dinger equation with spatial white noise potential.
\newblock {\em Electron. J. Probab.}, 23:Paper No. 28, 16, 2018.

\bibitem[EH17]{erhard2017discretisation}
D.~Erhard and M.~Hairer.
\newblock Discretisation of regularity structures.
\newblock {\em arXiv:1705.02836}, 2017.

\bibitem[Fee14]{feehan2014global}
P.~M. Feehan.
\newblock Global existence and convergence of smooth solutions to yang-mills
  gradient flow over compact four-manifolds.
\newblock {\em arXiv:1409.1525}, 2014.

\bibitem[FG17]{furlan2017weak}
M.~Furlan and M.~Gubinelli.
\newblock Weak universality for a class of 3d stochastic reaction--diffusion
  models.
\newblock {\em Probability Theory and Related Fields}, pages 1--66, 2017.

\bibitem[FG19]{MR3916943}
M.~Furlan and M.~Gubinelli.
\newblock Paracontrolled quasilinear {SPDE}s.
\newblock {\em Ann. Probab.}, 47(2):1096--1135, 2019.

\bibitem[FGS16]{TGS16}
T.~Franco, P.~Gon\c{c}alves, and M.~Simon.
\newblock Crossover to the stochastic {B}urgers equation for the {WASEP} with a
  slow bond.
\newblock {\em Commun. Math. Phys.}, 346(3):801--838, 2016.

\bibitem[FH14]{MR3289027}
P.~K. Friz and M.~Hairer.
\newblock {\em A course on rough paths: with an introduction to regularity
  structures}.
\newblock Universitext. Springer, Cham, 2014.

\bibitem[FH17]{MR3653951}
T.~Funaki and M.~Hoshino.
\newblock A coupled {KPZ} equation, its two types of approximations and
  existence of global solutions.
\newblock {\em J. Funct. Anal.}, 273(3):1165--1204, 2017.

\bibitem[Fla08]{flandoli2008introduction}
F.~Flandoli.
\newblock An introduction to 3d stochastic fluid dynamics.
\newblock In {\em SPDE in hydrodynamic: recent progress and prospects}, pages
  51--150. Springer, 2008.

\bibitem[Fle75]{MR0378124}
W.~H. Fleming.
\newblock Distributed parameter stochastic systems in population biology.
\newblock pages 179--191. Lecture Notes in Econom. and Math. Systems, Vol. 107,
  1975.

\bibitem[FOT11]{Fukushima2011}
M.~Fukushima, Y.~Oshima, and M.~Takeda.
\newblock {\em Dirichlet forms and symmetric {M}arkov processes}, volume~19 of
  {\em De Gruyter Studies in Mathematics}.
\newblock Walter de Gruyter \& Co., Berlin, extended edition, 2011.

\bibitem[FR95]{MR1358083}
J.~Fritz and B.~R\"{u}diger.
\newblock Time dependent critical fluctuations of a one-dimensional local mean
  field model.
\newblock {\em Probab. Theory Related Fields}, 103(3):381--407, 1995.

\bibitem[Fri87]{MR894407}
J.~Fritz.
\newblock On the hydrodynamic limit of a one-dimensional {G}inzburg-{L}andau
  lattice model. {T}he a priori bounds.
\newblock {\em J. Statist. Phys.}, 47(3-4):551--572, 1987.

\bibitem[FS81]{MR634447}
J.~Fr\"{o}hlich and T.~Spencer.
\newblock The {K}osterlitz-{T}houless transition in two-dimensional abelian
  spin systems and the {C}oulomb gas.
\newblock {\em Comm. Math. Phys.}, 81(4):527--602, 1981.

\bibitem[Fun83]{MR692348}
T.~Funaki.
\newblock Random motion of strings and related stochastic evolution equations.
\newblock {\em Nagoya Math. J.}, 89:129--193, 1983.

\bibitem[G{\"a}r88]{GJ88}
J.~G{\"a}rtner.
\newblock Convergence towards {B}urgers' equation and propagation of chaos for
  weakly asymmetric exclusion processes.
\newblock {\em Stoch. Proc. Appl.}, 27(2):233--260, 1988.

\bibitem[GGF{\etalchar{+}}12]{ghofraniha2012shock}
N.~Ghofraniha, S.~Gentilini, V.~Folli, E.~DelRe, and C.~Conti.
\newblock Shock waves in disordered media.
\newblock {\em Physical review letters}, 109(24):243902, 2012.

\bibitem[GH17a]{GerencserHairer}
M.~Gerencs{\'e}r and M.~Hairer.
\newblock Singular {SPDEs} in domains with boundaries.
\newblock {\em Probability Theory and Related Fields}, pages 1--62, 2017.

\bibitem[GH17b]{GerencserQuasilinear}
M.~Gerencs{\'e}r and M.~Hairer.
\newblock A solution theory for quasilinear singular {SPDEs}.
\newblock {\em Communications on Pure and Applied Mathematics}, 2017.

\bibitem[GH18a]{gubinelli2018pde}
M.~Gubinelli and M.~Hofmanov{\'a}.
\newblock {A PDE construction of the Euclidean $\Phi^4_3 $ quantum field
  theory}.
\newblock {\em arXiv:1810.01700}, 2018.

\bibitem[GH18b]{gubinelli2018global}
M.~Gubinelli and M.~Hofmanov{\'a}.
\newblock Global solutions to elliptic and parabolic $\phi^4$ models in
  {E}uclidean space.
\newblock {\em arXiv:1804.11253}, 2018.

\bibitem[Gho17]{G17}
P.~Ghosal.
\newblock {Hall-Littlewood-PushTASEP and its KPZ limit}.
\newblock {\em arXiv:1701.07308}, 2017.

\bibitem[GIP15]{MR3406823}
M.~Gubinelli, P.~Imkeller, and N.~Perkowski.
\newblock Paracontrolled distributions and singular {PDE}s.
\newblock {\em Forum Math. Pi}, 3:e6, 75, 2015.

\bibitem[GJ10]{GJ10}
P.~{Goncalves} and M.~{Jara}.
\newblock {Universality of KPZ} equation.
\newblock {\em arXiv:1003.4478}, 2010.

\bibitem[GJ13]{GJ13}
M.~Gubinelli and M.~Jara.
\newblock Regularization by noise and stochastic {B}urgers equations.
\newblock {\em Stoch. Partial Differ. Equ. Anal. Comput.}, 1(2):325--350, 2013.

\bibitem[GJ14]{GJ14a}
P.~Gon\c{c}alves and M.~Jara.
\newblock Nonlinear fluctuations of weakly asymmetric interacting particle
  systems.
\newblock {\em Arch. Ration. Mech. Anal.}, 212(2):597--644, 2014.

\bibitem[GJ17]{GJ16}
P.~Gon\c{c}alves and M.~Jara.
\newblock Stochastic {B}urgers equation from long range exclusion interactions.
\newblock {\em Stochastic Process. Appl.}, 127(12):4029--4052, 2017.

\bibitem[GJS15]{GJS15}
P.~Gon{\c{c}}alves, M.~Jara, and S.~Sethuraman.
\newblock A stochastic {B}urgers equation from a class of microscopic
  interactions.
\newblock {\em Ann. Probab.}, 43(1):286--338, 2015.

\bibitem[GKO18a]{gubinelli2018paracontrolled}
M.~Gubinelli, H.~Koch, and T.~Oh.
\newblock Paracontrolled approach to the three-dimensional stochastic nonlinear
  wave equation with quadratic nonlinearity.
\newblock {\em arXiv preprint arXiv:1811.07808}, 2018.

\bibitem[GKO18b]{MR3841850}
M.~Gubinelli, H.~Koch, and T.~Oh.
\newblock Renormalization of the two-dimensional stochastic nonlinear wave
  equations.
\newblock {\em Trans. Amer. Math. Soc.}, 370(10):7335--7359, 2018.

\bibitem[GKR18]{MR3762097}
Y.~Gu, T.~Komorowski, and L.~Ryzhik.
\newblock The {S}chr\"{o}dinger equation with spatial white noise: the average
  wave function.
\newblock {\em J. Funct. Anal.}, 274(7):2113--2138, 2018.

\bibitem[Gla63]{doi:10.1063/1.1703954}
R.~J. Glauber.
\newblock Time dependent statistics of the ising model.
\newblock {\em Journal of Mathematical Physics}, 4(2):294--307, 1963.

\bibitem[GP16]{Gubinelli2016hairer}
M.~Gubinelli and N.~Perkowski.
\newblock The {Hairer--Quastel} universality result in equilibrium.
\newblock {\em arXiv:1602.02428}, 2016.

\bibitem[GP17a]{GP2015a}
M.~Gubinelli and N.~Perkowski.
\newblock {E}nergy solutions of {KPZ} are unique.
\newblock {\em J. Amer. Math. Soc.}, 31:427--471, 2017.

\bibitem[GP17b]{GubPerk}
M.~Gubinelli and N.~Perkowski.
\newblock {KPZ} reloaded.
\newblock {\em Commun. Math. Phys.}, 349:165--269, 2017.

\bibitem[GP18a]{gubinelli2018infinitesimal}
M.~Gubinelli and N.~Perkowski.
\newblock The infinitesimal generator of the stochastic {B}urgers equation.
\newblock {\em arXiv:1810.12014}, 2018.

\bibitem[GP18b]{MR3828162}
M.~Gubinelli and N.~Perkowski.
\newblock An introduction to singular {SPDE}s.
\newblock In {\em Stochastic partial differential equations and related
  fields}, volume 229 of {\em Springer Proc. Math. Stat.}, pages 69--99.
  Springer, Cham, 2018.

\bibitem[GP18c]{MR3828193}
M.~Gubinelli and N.~Perkowski.
\newblock Probabilistic approach to the stochastic {B}urgers equation.
\newblock In {\em Stochastic partial differential equations and related
  fields}, volume 229 of {\em Springer Proc. Math. Stat.}, pages 515--527.
  Springer, Cham, 2018.

\bibitem[GPS17]{GonPerMar}
P.~Goncalves, N.~Perkowski, and M.~Simon.
\newblock Derivation of the stochastic {B}urgers equation with {D}irichlet
  boundary conditions from the {WASEP}.
\newblock {\em arXiv:1710.11011}, 2017.

\bibitem[GPV88]{MR954674}
M.~Z. Guo, G.~C. Papanicolaou, and S.~R.~S. Varadhan.
\newblock Nonlinear diffusion limit for a system with nearest neighbor
  interactions.
\newblock {\em Comm. Math. Phys.}, 118(1):31--59, 1988.

\bibitem[GS92]{GS92}
L.-H. Gwa and H.~Spohn.
\newblock Six-vertex model, roughened surfaces, and an asymmetric spin
  {H}amiltonian.
\newblock {\em Phys. Rev. Lett.}, 68:725--728, 1992.

\bibitem[Gu18]{Gu2018KPZ}
Y.~Gu.
\newblock Gaussian fluctuations of the {2D KPZ} equation.
\newblock {\em arXiv:1812.07467}, 2018.

\bibitem[Gub04]{MR2091358}
M.~Gubinelli.
\newblock Controlling rough paths.
\newblock {\em J. Funct. Anal.}, 216(1):86--140, 2004.

\bibitem[Gub18]{Gubinelli2018panorama}
M.~Gubinelli.
\newblock A panorama of singular {SPDEs}.
\newblock In {\em Proc. Int. Cong. of Math}, volume~2, pages 2277--2304, 2018.

\bibitem[GUZ18]{gubinelli2018semilinear}
M.~Gubinelli, B.~E. Ugurcan, and I.~Zachhuber.
\newblock Semilinear evolution equations for the {Anderson Hamiltonian} in two
  and three dimensions.
\newblock {\em arXiv preprint arXiv:1807.06825}, 2018.

\bibitem[Hai13]{Hairer13}
M.~Hairer.
\newblock Solving the {KPZ} equation.
\newblock {\em Ann. Math.}, 178(2):559--664, 2013.

\bibitem[Hai14a]{HairerICM2014}
M.~Hairer.
\newblock Singular stochastic {PDEs}.
\newblock {\em Proceedings of the 2014 ICM, arXiv:1403.6353}, 2014.

\bibitem[Hai14b]{Hairer14}
M.~Hairer.
\newblock A theory of regularity structures.
\newblock {\em Invent. Math.}, 198(2):269--504, 2014.

\bibitem[Hai15a]{MR3336866}
M.~Hairer.
\newblock Introduction to regularity structures.
\newblock {\em Braz. J. Probab. Stat.}, 29(2):175--210, 2015.

\bibitem[Hai15b]{hairer2015regularity}
M.~Hairer.
\newblock Regularity structures and the dynamical $\phi^4_3$ model.
\newblock {\em arXiv:1508.05261}, 2015.

\bibitem[Hai16]{hairer2016motion}
M.~Hairer.
\newblock The motion of a random string.
\newblock {\em arXiv:1605.02192}, 2016.

\bibitem[HL15]{MR3358965}
M.~Hairer and C.~Labb\'{e}.
\newblock A simple construction of the continuum parabolic {A}nderson model on
  {${\bf R}^2$}.
\newblock {\em Electron. Commun. Probab.}, 20:no. 43, 11, 2015.

\bibitem[HL18]{MR3779690}
M.~Hairer and C.~Labb\'{e}.
\newblock Multiplicative stochastic heat equations on the whole space.
\newblock {\em J. Eur. Math. Soc. (JEMS)}, 20(4):1005--1054, 2018.

\bibitem[HM18]{HaiMat}
M.~Hairer and K.~Matetski.
\newblock Discretisations of rough stochastic {PDE}s.
\newblock {\em Ann. Probab.}, 46(3):1651--1709, 2018.

\bibitem[Hos18]{MR3769661}
M.~Hoshino.
\newblock Paracontrolled calculus and {F}unaki-{Q}uastel approximation for the
  {KPZ} equation.
\newblock {\em Stochastic Process. Appl.}, 128(4):1238--1293, 2018.

\bibitem[HP15]{MR3417505}
M.~Hairer and E.~Pardoux.
\newblock A {W}ong-{Z}akai theorem for stochastic {PDE}s.
\newblock {\em J. Math. Soc. Japan}, 67(4):1551--1604, 2015.

\bibitem[HQ18]{KPZJeremy}
M.~Hairer and J.~Quastel.
\newblock A class of growth models rescaling to {KPZ}.
\newblock {\em Forum Math. Pi}, 6:e3, 112, 2018.

\bibitem[HS16]{MR3452276}
M.~Hairer and H.~Shen.
\newblock The dynamical sine-{G}ordon model.
\newblock {\em Comm. Math. Phys.}, 341(3):933--989, 2016.

\bibitem[HS17]{HaiShen}
M.~Hairer and H.~Shen.
\newblock A central limit theorem for the {KPZ} equation.
\newblock {\em Ann. Probab.}, 45(6B):4167--4221, 2017.

\bibitem[HX18a]{MR3772400}
M.~Hairer and W.~Xu.
\newblock Large-scale behavior of three-dimensional continuous phase
  coexistence models.
\newblock {\em Comm. Pure Appl. Math.}, 71(4):688--746, 2018.

\bibitem[HX18b]{HairerXu2018}
M.~Hairer and W.~Xu.
\newblock Large-scale limit of interface fluctuation models.
\newblock {\em arXiv:1802.08192}, 2018.

\bibitem[Jaf00]{MR1773042}
A.~Jaffe.
\newblock Constructive quantum field theory.
\newblock In {\em Mathematical physics 2000}, pages 111--127. Imp. Coll. Press,
  London, 2000.

\bibitem[Jan97]{janson_1997}
S.~Janson.
\newblock {\em Gaussian Hilbert Spaces}.
\newblock Cambridge Tracts in Mathematics. Cambridge University Press, 1997.

\bibitem[JLM85]{MR815192}
G.~Jona-Lasinio and P.~K. Mitter.
\newblock On the stochastic quantization of field theory.
\newblock {\em Comm. Math. Phys.}, 101(3):409--436, 1985.

\bibitem[Jos13]{jos201340}
J.~V. Jos.
\newblock {\em 40 years of {Berezinskii-Kosterlitz-Thouless} theory}.
\newblock World Scientific, 2013.

\bibitem[K\"16]{MR3526112}
W.~K\"{o}nig.
\newblock {\em The parabolic {A}nderson model: Random walk in random
  potential}.
\newblock Pathways in Mathematics. Birkh\"{a}user/Springer, [Cham], 2016.

\bibitem[KM17]{MR3607594}
A.~Kupiainen and M.~Marcozzi.
\newblock Renormalization of generalized {KPZ} equation.
\newblock {\em J. Stat. Phys.}, 166(3-4):876--902, 2017.

\bibitem[KPZ86]{KPZ86}
M.~Kardar, G.~Parisi, and Y.-C. Zhang.
\newblock Dynamic scaling of growing interfaces.
\newblock {\em Phys. Rev. Lett.}, 56:889--892, 1986.

\bibitem[KR82]{MR683274}
N.~V. Krylov and B.~L. Rozovski\u{\i}.
\newblock Stochastic partial differential equations and diffusion processes.
\newblock {\em Uspekhi Mat. Nauk}, 37(6(228)):75--95, 1982.

\bibitem[KS91a]{MR1121940}
I.~Karatzas and S.~E. Shreve.
\newblock {\em Brownian motion and stochastic calculus}, volume 113 of {\em
  Graduate Texts in Mathematics}.
\newblock Springer-Verlag, New York, second edition, 1991.

\bibitem[KS91b]{krug1991kinetic}
J.~Krug and H.~Spohn.
\newblock Kinetic roughening of growing surfaces in {``Solids Far From
  Equilibrium: Growth, Morphology and Defects (C. Godreche, ed.)''}, 1991.

\bibitem[KS12]{MR3443633}
S.~Kuksin and A.~Shirikyan.
\newblock {\em Mathematics of two-dimensional turbulence}, volume 194 of {\em
  Cambridge Tracts in Mathematics}.
\newblock Cambridge University Press, Cambridge, 2012.

\bibitem[KT73]{kosterlitz1973ordering}
J.~M. Kosterlitz and D.~J. Thouless.
\newblock Ordering, metastability and phase transitions in two-dimensional
  systems.
\newblock {\em Journal of Physics C: Solid State Physics}, 6(7):1181, 1973.

\bibitem[Kup10]{kupiainen2010ergodicity}
A.~Kupiainen.
\newblock Ergodicity of two dimensional turbulence.
\newblock {\em arXiv:1005.0587}, 2010.

\bibitem[Kup16]{Antti}
A.~Kupiainen.
\newblock Renormalization group and stochastic {PDE}s.
\newblock {\em Ann. H. Poincare}, 17:497--535, 2016.

\bibitem[Lab17]{Labbe16b}
C.~Labb\'e.
\newblock Weakly asymmetric bridges and the {KPZ} equation.
\newblock {\em Commun. Math. Phys.}, 353(3):1261--1298, 2017.

\bibitem[LR15]{MR3410409}
W.~Liu and M.~R\"{o}ckner.
\newblock {\em Stochastic partial differential equations: an introduction}.
\newblock Universitext. Springer, Cham, 2015.

\bibitem[Lyo98]{MR1654527}
T.~J. Lyons.
\newblock Differential equations driven by rough signals.
\newblock {\em Rev. Mat. Iberoamericana}, 14(2):215--310, 1998.

\bibitem[Mat03]{MR2050597}
J.~C. Mattingly.
\newblock On recent progress for the stochastic {N}avier {S}tokes equations.
\newblock In {\em Journ\'{e}es ``\'{E}quations aux {D}\'{e}riv\'{e}es
  {P}artielles''}, pages Exp. No. XI, 52. Univ. Nantes, Nantes, 2003.

\bibitem[Mat18]{Matetski2018martingale}
K.~Matetski.
\newblock Martingale-driven approximations of singular stochastic {PDEs}.
\newblock {\em arXiv:1808.09429}, 2018.

\bibitem[McK95]{MR1328250}
H.~P. McKean.
\newblock Statistical mechanics of nonlinear wave equations. {IV}. {C}ubic
  {S}chr\"{o}dinger.
\newblock {\em Comm. Math. Phys.}, 168(3):479--491, 1995.

\bibitem[MP17]{martin2017paracontrolled}
J.~Martin and N.~Perkowski.
\newblock Paracontrolled distributions on bravais lattices and weak
  universality of the 2d parabolic anderson model.
\newblock {\em arXiv:1704.08653}, 2017.

\bibitem[MS77]{MR0441179}
O.~A. McBryan and T.~Spencer.
\newblock On the decay of correlations in {${\rm SO}(n)$}-symmetric
  ferromagnets.
\newblock {\em Comm. Math. Phys.}, 53(3):299--302, 1977.

\bibitem[MU18]{MR3790153}
J.~Magnen and J.~Unterberger.
\newblock The scaling limit of the {KPZ} equation in space dimension 3 and
  higher.
\newblock {\em J. Stat. Phys.}, 171(4):543--598, 2018.

\bibitem[Mue91]{Mueller91}
C.~Mueller.
\newblock On the support of solutions to the heat equation with noise.
\newblock {\em Stoch. \& Stoch. Rep.}, 37(4):225--245, 1991.

\bibitem[MW17a]{MR3628883}
J.-C. Mourrat and H.~Weber.
\newblock Convergence of the two-dimensional dynamic {I}sing-{K}ac model to
  {$\Phi^4_2$}.
\newblock {\em Comm. Pure Appl. Math.}, 70(4):717--812, 2017.

\bibitem[MW17b]{MR3719541}
J.-C. Mourrat and H.~Weber.
\newblock The dynamic {$\Phi^4_3$} model comes down from infinity.
\newblock {\em Comm. Math. Phys.}, 356(3):673--753, 2017.

\bibitem[MW18]{moinat2018space}
A.~Moinat and H.~Weber.
\newblock Space-time localisation for the dynamic $\phi^4_3 $ model.
\newblock {\em arXiv:1811.05764}, 2018.

\bibitem[OR98]{MR1640497}
M.~Oberguggenberger and F.~Russo.
\newblock Nonlinear stochastic wave equations.
\newblock {\em Integral Transform. Spec. Funct.}, 6(1-4):71--83, 1998.
\newblock Generalized functions---linear and nonlinear problems (Novi Sad,
  1996).

\bibitem[OW18]{otto2018quasilinear}
F.~Otto and H.~Weber.
\newblock Quasilinear {SPDEs} via rough paths.
\newblock {\em Archive for Rational Mechanics and Analysis}, pages 1--78, 2018.

\bibitem[Par17]{Shalin}
S.~Parekh.
\newblock The {KPZ} limit of {ASEP} with boundary.
\newblock {\em arXiv:1711.05297}, 2017.

\bibitem[Pau35]{LPauling35}
L.~Pauling.
\newblock The structure and entropy of ice and of other crystals with some
  randomness of atomic arrangement.
\newblock {\em J. Amer. Chem. Soc.}, 57(12):2680--2684, 1935.

\bibitem[PR07]{MR2329435}
C.~Pr\'{e}v\^{o}t and M.~R\"{o}ckner.
\newblock {\em A concise course on stochastic partial differential equations},
  volume 1905 of {\em Lecture Notes in Mathematics}.
\newblock Springer, Berlin, 2007.

\bibitem[PR18]{perkowski2018line}
N.~Perkowski and T.~C. Rosati.
\newblock The {KPZ} equation on the real line.
\newblock {\em arXiv:1808.00354}, 2018.

\bibitem[PW81]{ParisiWu}
G.~Parisi and Y.~S. Wu.
\newblock Perturbation theory without gauge fixing.
\newblock {\em Sci. Sinica}, 24(4):483--496, 1981.

\bibitem[RWZZ18]{MR3787728}
M.~R\"{o}ckner, B.~Wu, R.~Zhu, and X.~Zhu.
\newblock Stochastic heat equations with values in a {R}iemannian manifold.
\newblock {\em Atti Accad. Naz. Lincei Rend. Lincei Mat. Appl.},
  29(1):205--213, 2018.

\bibitem[She07]{MR2322706}
S.~Sheffield.
\newblock Gaussian free fields for mathematicians.
\newblock {\em Probab. Theory Related Fields}, 139(3-4):521--541, 2007.

\bibitem[She18]{Shen2018Abelian}
H.~Shen.
\newblock Stochastic quantization of an {A}belian gauge theory.
\newblock {\em arXiv:1801.04596}, 2018.

\bibitem[Spo86]{MR826856}
H.~Spohn.
\newblock Equilibrium fluctuations for interacting {B}rownian particles.
\newblock {\em Comm. Math. Phys.}, 103(1):1--33, 1986.

\bibitem[Spo14]{MR3176405}
H.~Spohn.
\newblock Nonlinear fluctuating hydrodynamics for anharmonic chains.
\newblock {\em J. Stat. Phys.}, 154(5):1191--1227, 2014.

\bibitem[ST18]{ST18}
H.~Shen and L.-C. Tsai.
\newblock Stochastic telegraph equation limit for the stochastic six vertex
  model.
\newblock {\em arXiv:1807.04678}, 2018.

\bibitem[SW18]{MR3820327}
H.~Shen and H.~Weber.
\newblock Glauber dynamics of 2{D} {K}ac-{B}lume-{C}apel model and their
  stochastic {PDE} limits.
\newblock {\em J. Funct. Anal.}, 275(6):1321--1367, 2018.

\bibitem[SX16]{ShenXu}
H.~Shen and W.~Xu.
\newblock Weak universality of dynamical $\phi^4_3$: non-{G}aussian noise.
\newblock {\em Stochastics and Partial Differential Equations: Analysis and
  Computations}, 2016.

\bibitem[Wal86]{MR876085}
J.~B. Walsh.
\newblock An introduction to stochastic partial differential equations.
\newblock In {\em \'{E}cole d'\'et\'e de probabilit\'es de {S}aint-{F}lour,
  {XIV}---1984}, volume 1180 of {\em Lecture Notes in Math.}, pages 265--439.
  Springer, Berlin, 1986.

\bibitem[Yan18]{Yang}
K.~Yang.
\newblock The {KPZ} equation, non-equilibrium solutions, and weak universality
  for long-range interactions.
\newblock {\em arXiv:1810.02836}, 2018.

\bibitem[Zhu90]{MR1044937}
M.~Zhu.
\newblock Equilibrium fluctuations for one-dimensional {G}inzburg-{L}andau
  lattice model.
\newblock {\em Nagoya Math. J.}, 117:63--92, 1990.

\bibitem[ZZ15]{Zhu2014NS}
R.~Zhu and X.~Zhu.
\newblock Three-dimensional {N}avier-{S}tokes equations driven by space-time
  white noise.
\newblock {\em J. Differential Equations}, 259(9):4443--4508, 2015.

\bibitem[ZZ18]{MR3758734}
R.~Zhu and X.~Zhu.
\newblock Lattice approximation to the dynamical {$\Phi_3^4$} model.
\newblock {\em Ann. Probab.}, 46(1):397--455, 2018.

\end{thebibliography}

\end{document}